\documentclass[preprint,3p,sort&compress]{elsarticle}

\usepackage{amssymb}

\usepackage{placeins}
\usepackage{tikz}
\usepackage{ams math,cases}
\usepackage{siunitx}

\newcommand{\tabref}[1]{Table~\ref{#1}} 
\newcommand{\figref}[1]{Figure~\ref{#1}} 
\newcommand{\secref}[1]{Section~\ref{#1}} 

\newcommand{\tikzcircle}[2][red,fill=red]{\tikz[baseline=-0.5ex]\draw[#1,radius=#2] (0,0) circle ;}

\definecolor{darkred}{rgb}{0.5,0,0}
\definecolor{darkgreen}{rgb}{0,0.3,0}
\definecolor{darkblue}{rgb}{0,0,0.5}
\definecolor{darkbrown}{rgb}{0.28,0.07,0.07}
\definecolor{black}{rgb}{0,0,0}
\definecolor{plotGreen}{rgb}{0,0.5,0}
\definecolor{plotRed}{rgb}{1.0,0,0}
\definecolor{plotBlue}{rgb}{0,0,1.0}
\definecolor{plotCyan}{rgb}{0,1.0,1.0}
\definecolor{plotGray}{gray}{0.25}

\begin{document}

\begin{frontmatter}



\title{A hybrid Eulerian--Lagrangian flow solver}


\author[cst,awep]{A.~Palha\corref{cor}}
\ead{a.palha@tue.nl}
\author[awep]{L.~Manickathan}
\ead{l.manickathan@student.tudelft.nl}
\author[awep]{C.S.~Ferreira}
\ead{c.j.simaoferreira@tudelft.nl}
\author[awep]{G. van Bussel}
\ead{g.j.w.vanbussel@tudelft.nl}

\cortext[cor]{Corresponding author}

\address[cst]{Eindhoven University of Technology, Department of Mechanical Engineering, section of Control Systems Technology P.O. Box 513, 5600 MB Eindhoven, The Netherlands}

\address[awep]{Delft University of Technology, Faculty of Aerospace Engineering, Department of Aerodynamics, Wind Energy Flight Performance and Propulsion , P.O. Box 5058, 2600~GB Delft, The Netherlands}

\begin{abstract}
	Currently, Eulerian flow solvers are very efficient in accurately resolving flow structures near solid boundaries. On the other hand, they tend to be diffusive and to dampen high-intensity vortical structures after a short distance away from solid boundaries. The use of high order methods and fine grids, although alleviating this problem, gives rise to large systems of equations that are expensive to solve. 
	
	Lagrangian solvers, such as the regularized vortex particle method, have shown to eliminate (in practice) the diffusion in the wake. However, as a drawback, the modelling of solid boundaries is less accurate, more complex and costly than with Eulerian solvers (due to the isotropy of its computational elements). 
	
	Given the drawbacks and advantages of both Eulerian and Lagrangian solvers the combination of both methods, giving rise to a hybrid solver, is advantageous. The main idea behind the hybrid solver presented is the following. In a region close to solid boundaries the flow is solved with an Eulerian solver, where the full Navier-Stokes equations are solved (possibly with an arbitrary turbulence model or DNS, the limitations being the computational power and the physical properties of the flow), outside of that region the flow is solved with a vortex particle method.
	
	In this work we present this hybrid scheme and verify it numerically on known 2D benchmark cases: dipole flow, flow around a cylinder and flow around a stalled airfoil. The success in modelling these flow conditions presents this hybrid approach as a promising alternative, bridging the gap between highly resolved and computationally intensive Eulerian CFD simulations and fast but less resolved Lagrangian simulations.

\end{abstract}


\end{frontmatter}


\section{Introduction}
	It is well known that the numerical simulation of advection dominated external flows is a challenging problem. In Eulerian flow solvers one of the biggest difficulties lies in the correct definition of boundary conditions, arising from the necessity to truncate the unbounded domain, see for example \cite{Sani1994,Heywood1996,Gresho1991}. A usual expedient to circumvent this is to employ a large computational domain, placing the outflow boundaries sufficiently far from the region of interest, minimizing errors. Although work has been done towards the reduction of the size of the computational domain required for numerical stability, see for example the recent work in \cite{Dong2014}, physical accuracy still imposes strong restrictions. Another important aspect of external flows is that vorticity is concentrated along the wake and in thin boundary layers surrounding solid boundaries. Therefore, in order to make the numerical solution accurate and at the same time computationally feasible, adaptive schemes should be employed, e.g. \cite{Howell1997}. In this way, high resolution is used in the support of vorticity whereas in the remaining regions lower resolution may be used, reducing the computational cost. Nonetheless, when several moving objects are present, the continuous generation and adaptation of the underlying mesh is not trivial and poses several challenges, see \cite{Tezduyar2001,Lesoinne1996} and references therein. Conversely, the Lagrangian vortex particle method, introduced in \cite{Chorin1973}, is a mesh free method where the vorticity is carried by particles which are then transported by the flow field. In this way, the compact vortical structures that are characteristic of the flow in wakes are naturally captured by this method, which is intrinsically adaptive since vortex particles are placed only where vorticity exists, \cite{Leonard1985_2003}. Additionally, boundary conditions at infinity are automatically satisfied. Although many developments have made this method a robust alternative to Eulerian flow solvers, some challenges still exist. For example, the enforcement of viscous boundary conditions, introduced in \cite{Koumoutsakos1995} and improved in \cite{Ploumhans2000}, is still a challenge on complex geometries, see \cite{Morgenthal2007,Chatelin2014,Rasmussen2011}.
	
	An alternative approach, followed in this work, is to divide the computational domain into different sub-domains and employ solvers specifically tailored to the flow characteristics present in each region. This gives rise to what is commonly referred to as a \emph{hybrid solver}. Within this framework it is possible to use Eulerian flow solvers in the vicinity of solid boundaries and a Lagrangian vortex particle solver for the wake. In this way, the Eulerian flow solvers are used to efficiently and accurately resolve the regions where viscous effects dominate the convective ones. For instance, the flow in boundary layers, which is predominantly anisotropic and unidirectional, is particularly suited for an Eulerian formulation where elongated cells may be used. On the other hand, the vortex particle method accurately and efficiently captures the wake flow. Moreover, by using a multi-domain formulation other advantages arise. First, since meshes exist only in the vicinity of solid boundaries, the movement of each body can be treated independently of all the others, thus greatly simplifying the implementation of moving bodies. Secondly, despite less favourable results in the past, \cite{Ploumhans2002,Rasmussen2011}, recent results start to establish that fast particle solvers such as the Fast Multipole Method (FMM),  \cite{Barnes1986,Greengard1987a}, can scale to massively parallel computations more efficiently than matrix and Fast Fourier Transform (FFT) solvers, \cite{BarbaFMMexa2013,Yokota2012,Lashuk2012}. By employing a hybrid formulation, large scale simulations involving several objects, may replace a single large Eulerian solver by a set of smaller ones coupled to an efficient mesh free vortex particle solver where velocities are computed by a FMM. Multidomain Eulerian solvers are well established, nevertheless, the authors consider that the prior division of the computational domain into disconnected sub-domains potentially leads to a more efficient solver. These potential improvements in computational efficiency have several applications, namely on the simulation of large wind farms where large numbers of objects move with respect to each other, different flow scales interact and the accurate propagation of the wake is fundamental to determine the optimal placement of wind turbines.
	
	\subsection{Literature review of hybrid grid-particle flow solvers}\label{subsection::literature_review_of_hybrid_grid_particla_flow_solvers}
		Although hybrid grid-particle methods are not a popular approach for solving flow problems, several formulations have been proposed in the past.
		
			Particle-grid methods have been formally introduced for the first time by Cottet, \cite{Cottet1991book}, where a ($\boldsymbol{u}$, $\omega$) formulation of the Navier-Stokes equations is used in the whole domain. In that work, the Eulerian and Lagrangian domains partially overlap and the interface conditions are obtained by an alternating Schwarz method. This approach was extended by Ould-Salihi et al., \cite{Ould-Salihi2001a}, to a ($\boldsymbol{u}$, $p$) formulation of the Navier-Stokes equations. A more detailed description of this hybrid formulation is presented by Cottet and Koumoutsakos, \cite{CottetKoumoutsakos2000}, both for the ($\boldsymbol{u}$, $\omega$) and a ($\boldsymbol{u}$, $p$) forms of the Navier-Stokes equations.
			
			One of the few formulations without overlap is the one by Guermond et al., \cite{Guermond1993,guermond1994}. In this approach the ($\psi$, $\omega$) form of the Navier-Stokes equations is solved in the Eulerian domain by a finite differences method. The transmission condition is obtained by assuming that viscosity is small in comparison to the convective term.
			
			Huberson and Voutsinas present an interesting overview of particle-grid methods in their overview article, \cite{Huberson2002}.
			
			More recently, two approaches related to the work presented here have been proposed that do not rely on the classical alternating Schwartz procedure. Daeninck, \cite{daeninckThesis}, presents a hybrid formulation where the particles cover the whole fluid domain and completely overlap the Eulerian one. The Eulerian solution is then used as a correction to the particle solver in the vicinity of the solid boundaries. Stock, Gharakhani and Stone, \cite{Stock2010}, proposed an improved version of this method in order to avoid interpolation errors within the viscous boundary layer, which for high Reynolds numbers can be large.
			
			Oxley, \cite{OxleyThesis}, and more recently Papadakis and Voutsinas, \cite{Papadakis2014}, present a 2D hybrid compressible Euler method for transonic rotorcraft applications using a complete overlap and an iteration procedure to eliminate projection errors between the particles and the grid and ensure compatibility between the solutions in the two computational domains.
			
			With a focus on computer graphics animation, Golas et al., \cite{Golas2012}, present a partially overlapping domain decomposition approach based on a ($\boldsymbol{u}$, $\omega$) formulation of the Navier-Stokes equations. 
		
	\subsection{Outline}\label{subsection::outline}
		As mentioned before, the goal of this work is to present a hybrid grid-particle flow solver. The main principle behind this hybrid approach is to decompose the fluid domain into a set of sub-domains where the most suitable flow solver is applied. In the vicinity of the boundary regions an Eulerian solver will be used, whereas the wake will be solved by Lagrangian particles.
		
		Therefore, the outline of this paper is as follows. In \secref{section::the_hybrid_flow_solver} we present the hybrid grid-particle solver, starting with a brief introduction to its key ingredients in \secref{subsection::hybrid_solver_intro}. This is followed by a discussion of the Lagrangian and Eulerian solvers in \secref{subsection::lagrangian_solver} and \secref{subsection::eulerian_solver}, respectively. In \secref{subsection::hybrid_solver}, we present in detail the coupling strategy between the two flow solvers. After introducing our approach, we apply it to different test cases in \secref{section::numerical_benchmark_cases}.  Finally, in \secref{section::conclusions_and_outlook} the conclusions and further outlook of this work are discussed.
		
\section{The hybrid flow solver} \label{section::the_hybrid_flow_solver}
	\subsection{Hybrid solver introduction}\label{subsection::hybrid_solver_intro}
\begin{figure}[!ht]
			\center
			\includegraphics[width=0.35\textwidth]{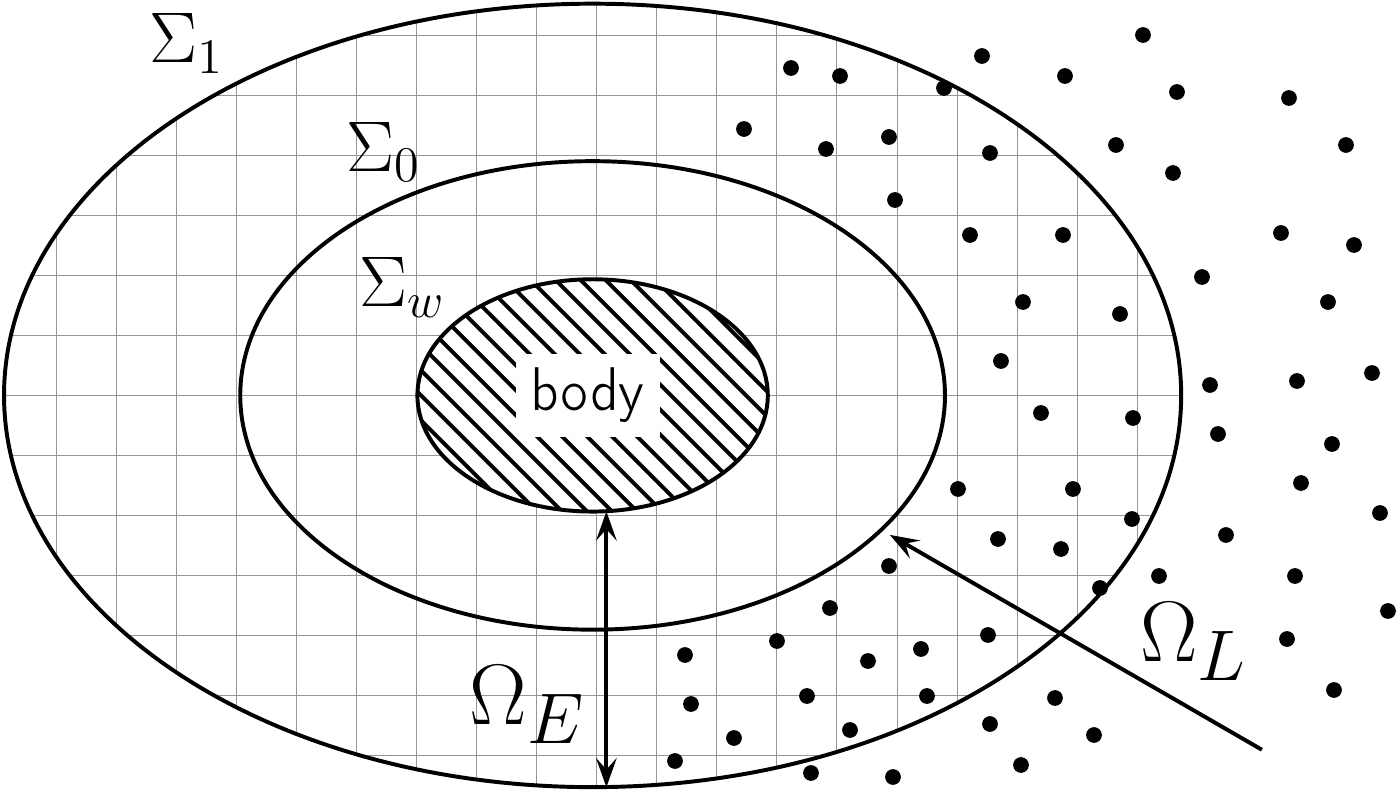}
			\caption{Standard geometrical decomposition of the flow domain with partial overlap. The Eulerian domain, $\Omega_{E}$, is a neighbourhood of the solid boundaries and the Lagrangian domain, $\Omega_{L}$, is a neighbourhood of infinity. Both domains overlap, such that $\Omega_{E}\cap\Omega_{L}\neq\emptyset$.}
			\label{fig::standard_domain_decomposition_overlap}
		\end{figure}

		As mentioned before, the objective of a hybrid solver is to use the most suitable numerical method in each region. In the case of the grid-particle approach discussed in this article, Eulerian grid solvers are used to discretize the flow equations in the vicinity of solid boundaries and a Lagrangian vortex particle formulation approximates the flow in the wake region. Most of the hybrid approaches discussed in \secref{subsection::literature_review_of_hybrid_grid_particla_flow_solvers} are based on domain decomposition with partial overlap, as depicted in \figref{fig::standard_domain_decomposition_overlap}. Another common characteristic of most of these formulations is the use of Schwarz's alternating method to ensure compatibility between the two solvers.
		
		Although the Schwarz alternating method is an excellent technique to obtain a very accurate match between the solutions of the two solvers in the overlap region, $\Omega_{E}\cap\Omega_{L}$, it requires iterations. The need to iterate makes this approach computationally expensive. For this reason, in this work we extend the ideas pioneered by Daeninck, \cite{daeninckThesis}, Stock, Gharakhani and Stone, \cite{Stock2010}, and propose a vorticity conserving hybrid grid-particle solver that does not require an iterative process to ensure compatibility.
		
\begin{figure}[!ht]
			\center
			\includegraphics[scale=.65]{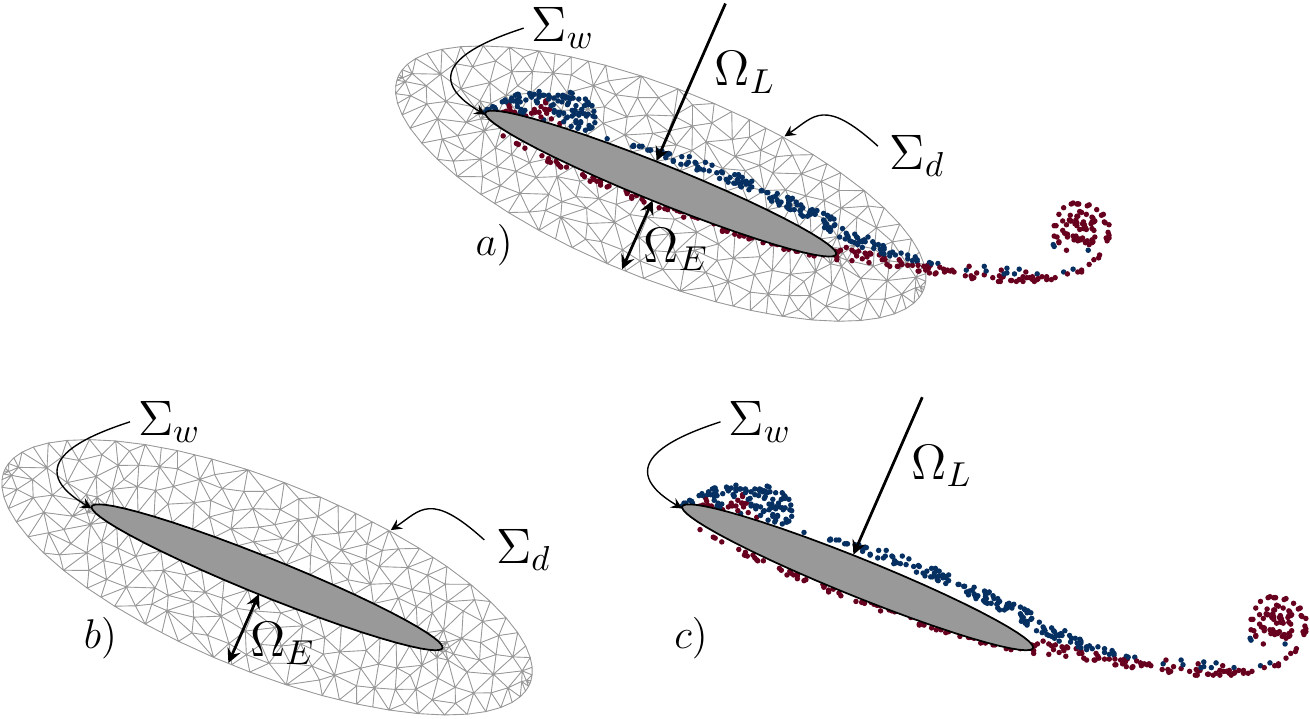}
			\caption{Representation of the different domains present in the domain decomposition strategy: presented in this work (a) Hybrid domain, (b) Eulerian domain and (c) Lagrangian domain. The Lagrangian domain solves the fluid flow in the whole domain and the Eulerian one plays the role of a near-wall correction.}
			\label{fig::pHyFlow_domain_decomposition_overlap_overview}
		\end{figure}

		In this work, the Lagrangian domain, $\Omega_{L}$, completely overlaps the Eulerian one, $\Omega_{E}$, as can be seen in \figref{fig::pHyFlow_domain_decomposition_overlap_overview}a. Instead of focusing on exactly matching the solution of both solvers in the overlap region, this approach follows another route. The Lagrangian vortex particle method is used to obtain the flow solution in the whole domain, \figref{fig::pHyFlow_domain_decomposition_overlap_overview}c, without resolving the boundary layers near the solid walls. The Eulerian solver is then used to correct it in the near-wall region, \figref{fig::pHyFlow_domain_decomposition_overlap_overview}b.
		
\begin{figure}[!ht]
			\center
			\includegraphics[scale=.85]{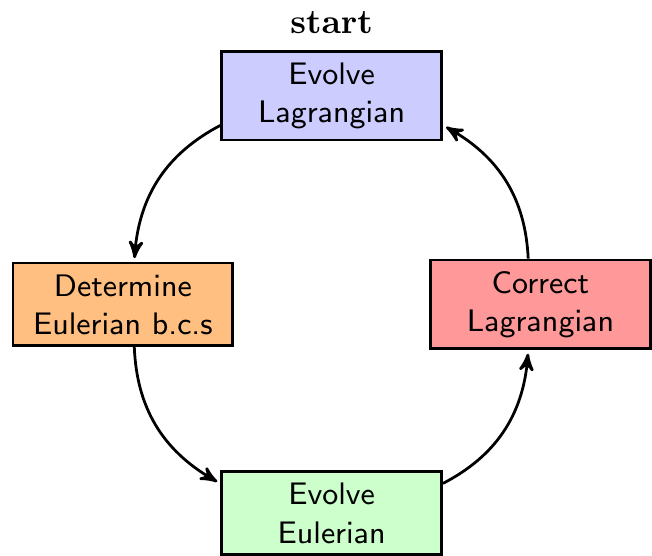}
%
%
			\caption{Flowchart of the coupling procedure to evolve the hybrid grid-particle solver from time instant $t_{n}$ to $t_{n+1}$.}
			\label{fig::pHyFlow_hybrid_flowchart}
		\end{figure}

		To evolve from time instant $t_{n}$ to $t_{n+1}$, the coupling procedure of the hybrid grid-particle solver discussed in this work can be summarised in the following steps:
		
		\begin{enumerate}
			\item \textbf{Evolve Lagrangian solution:} Evolve the vortex particles from time instant $t_{n}$ to $t_{n+1}$ neglecting vortex generation at the solid boundaries, following a standard procedure as presented, for example, by Cottet and Koumoutsakos in \cite{CottetKoumoutsakos2000}.
			\item \textbf{Determine Eulerian boundary conditions:} Use the Lagrangian solution at time instant $t_{n+1}$ to compute the Dirichlet boundary conditions for the velocity field in the Eulerian domain.
			\item \textbf{Evolve Eulerian solution:} Evolve the Eulerian solution from time instant $t_{n}$ to $t_{n+1}$ (possibly using $k_{E}$ sub steps in order to ensure stability of the solution), using any standard Navier-Stokes grid solver based on the $(\boldsymbol{u},p)$ formulation.
			\item \textbf{Correct Lagrangian solution:} Use the Eulerian solution at time instant $t_{n+1}$ to correct the Lagrangian one in the near wall region, $\Omega_{E}\cap\Omega_{L}$.
		\end{enumerate}
		These steps are performed in a loop, as depicted in the flowchart presented in \figref{fig::pHyFlow_hybrid_flowchart}.
		
		In order to introduce the hybrid grid-particle solver we will first discuss the Lagrangian vortex particle method used in this work, followed by the presentation of the finite element discretization employed in the Eulerian domain.
		
		\FloatBarrier

	\subsection{Lagrangian particle solver}\label{subsection::lagrangian_solver}
		Since this work is intended for a wider audience, potentially unfamiliar with the vortex particle method (VPM), in this section we give a brief introduction to this topic, emphasising the specific details of the implementation used in this work. For a more detailed overview of vortex methods we recommend the extensive reference book by Cottet and Koumoutsakos, \cite{CottetKoumoutsakos2000}, and the briefer reviews by Raviart, \cite{RaviartVortex1985}, and Winckelmans, \cite{WinckelmansGeneralVortex2004}. Since this work focuses on two dimensional flows, we will restrict our analysis to two dimensional domains.
		
		On a simply connected two dimensional unbounded domain, $\Omega=\mathbb{R}^{2}$, and for time interval $t\in (0,T]$, the system of the Navier-Stokes (NS) equations describing incompressible fluid dynamics in the $(\boldsymbol{u},\omega)$ formulation in the absence of external forces is given by:
		\begin{numcases}{}
				\frac{\partial\omega}{\partial t} + \left(\boldsymbol{u}\cdot\nabla\right)\omega -\nu\nabla^{2}\omega = 0\,, & in $\Omega$\,, \label{eq::ns_momentum_vorticity} \\
				\nabla\cdot\boldsymbol{u} = 0\,, & in $\Omega$\,, \label{eq::ns_incompressibility_2} \\
				\nabla\times\boldsymbol{u} = \omega\,. & in $\Omega$\,, \label{eq::ns_definition_vorticity} \\
				\omega(\boldsymbol{x},t) = \omega_{0}(\boldsymbol{x})\,, & in $\Omega$ and for $t=0$\,, \\
				\lim_{\|\boldsymbol{x}\|\rightarrow\infty}\boldsymbol{u}(\boldsymbol{x},t) = \boldsymbol{U}_{\infty}\,. & in $\Omega$ and for $t\in]0,T]$\,, \label{eq::velocity_bc_infinity_2} \\
				\lim_{\|\boldsymbol{x}\|\rightarrow\infty}\omega(\boldsymbol{x},t) = 0\,. & in $\Omega$ and for $t\in]0,T]$\,, \label{eq::vorticity_bc_infinity}
		\end{numcases}
		
\begin{figure}[!ht]
			\center
			\includegraphics[scale=1.0]{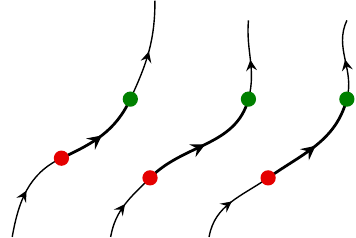}
			\caption{Fluid parcels containing vorticity are traced along field lines. Where the positions of fluid parcels at time instant $t=t_{0}$ are represented by [\tikzcircle[red!90!black,fill=red!90!black]{1.5pt}] and the positions of fluid parcels at a time instant $t>t_{0}$ are represented by [\tikzcircle[green!50!black,fill=green!50!black]{1.5pt}].}
			\label{fig::vorticity_parcels_trace}
		\end{figure}

		
		If we consider inviscid flow, $\nu=0$, \eqref{eq::ns_momentum_vorticity} reduces to:
		\begin{equation}
			\frac{D\omega}{D t} := \frac{\partial \omega}{\partial t} + \left(\boldsymbol{u}\cdot\nabla\right)\omega = 0\,,
		\end{equation}
		where $\frac{D}{D t}$ is the material derivative. This equation simply states that vorticity is transported along the velocity field lines, see \figref{fig::vorticity_parcels_trace}. This points to a Lagrangian formulation where fluid parcels containing vorticity are traced along field lines.

		The vortex particle method stems from this Lagrangian formulation. For incompressible and inviscid flows, the vortex particle method assumes that fluid parcels convect without deformation, therefore the position, $\boldsymbol{x}_{p}$, of each particle, $p$, evolves according to the following ordinary differential equation:		
		\begin{equation}
			\frac{\mathrm{d} \boldsymbol{x}_{p}}{\mathrm{d}t} = \boldsymbol{u}(\boldsymbol{x}_{p},t)\,,  \label{eq::vortex_equations_inviscid}
		\end{equation}
		where $\boldsymbol{u}(\boldsymbol{x},t)$ is the velocity field associated with the vorticity distribution $\omega(\boldsymbol{x},t)$ due to the vortex particles. Given the incompressibility constraint, \eqref{eq::ns_incompressibility_2}, and the definition of vorticity, \eqref{eq::ns_definition_vorticity}, it is possible to show that the velocity, $\boldsymbol{u}$, and vorticity, $\omega$, are related by the following Poisson's equation:
		\begin{equation}
			\nabla^{2}\boldsymbol{u} = \nabla\times\left(\omega\boldsymbol{e}_{z}\right)\,, \label{eq::poisson_velocity}
		\end{equation}
		with boundary conditions as given in \eqref{eq::velocity_bc_infinity_2}.
		
		Green's function formulation can be used to directly compute the velocity field from the vorticity field by noting that the velocity, $\boldsymbol{u}_{\delta}$, associated with a Dirac-$\delta$ distribution of vorticity, $\omega_{\delta}$, located at $\boldsymbol{x}_{p}$ satisfies \eqref{eq::poisson_velocity},
		\[
			\nabla^{2}\boldsymbol{u}_{\delta} = \nabla\times\left(\omega_{\delta}\boldsymbol{e}_{z}\right) :=  \nabla\times\left(\delta(\boldsymbol{x}-\boldsymbol{x_{p}})\boldsymbol{e}_{z}\right)\,.
		\]
		For homogeneous boundary conditions, $\boldsymbol{u}_{\delta}$ is given by the expression:
		\[
			\boldsymbol{u}_{\delta} = -\frac{1}{2\pi}\frac{\boldsymbol{x}-\boldsymbol{x}_{p}}{\|\boldsymbol{x}-\boldsymbol{x}_{p}\|^{2}}\times\boldsymbol{e}_{z}\,.
		\]
		Due to the linearity of Poisson's equation, a given vorticity distribution, $\omega(\boldsymbol{x})$, together with the boundary conditions \eqref{eq::velocity_bc_infinity_2} generates the following velocity field:
		\begin{equation}
			\boldsymbol{u}(\boldsymbol{x}) = -\frac{1}{2\pi} \int_{\Omega} \frac{\boldsymbol{x}-\tilde{\boldsymbol{x}}}{\|\boldsymbol{x}-\tilde{\boldsymbol{x}}\|^{2}}\times\boldsymbol{e}_{z} \,\omega(\tilde{\boldsymbol{x}})\mathrm{d}\tilde{\boldsymbol{x}} + \boldsymbol{U}_{\infty}\,. \label{eq::greens_velocity_general}
		\end{equation}
		If we associate to each vortex particle, $p$, a vorticity distribution, $\omega_{p}(\boldsymbol{x})$, the total vorticity distribution, $\omega(\boldsymbol{x})$, will be given by
		\begin{equation}
			\omega(\boldsymbol{x}) = \sum_{p} \omega_{p}(\boldsymbol{x})\,,\label{eq::vorticity_blob}
		\end{equation}
		and the total velocity field will be, replacing \eqref{eq::vorticity_blob} into \eqref{eq::greens_velocity_general},
		\begin{equation}
			\boldsymbol{u}(\boldsymbol{x}) = -\frac{1}{2\pi}\sum_{p} \int_{\Omega} \frac{\boldsymbol{x}-\tilde{\boldsymbol{x}}}{\|\boldsymbol{x}-\tilde{\boldsymbol{x}}\|^{2}}\times\boldsymbol{e}_{z} \,\omega_{p}(\tilde{\boldsymbol{x}})\mathrm{d}\tilde{\boldsymbol{x}} + \boldsymbol{U}_{\infty}\,. \label{eq::greens_velocity_particles}
		\end{equation}
		In the vortex particle method the vorticity distribution associated with each particle, $\omega_{p}(\boldsymbol{x})$, is given by
		\begin{equation}
			\omega_{p}(\boldsymbol{x}) := \zeta_{\sigma}(\|\boldsymbol{x}-\boldsymbol{x}_{p}\|)\,\Gamma_{p}\,,\label{eq::vorticity_blobs}
		\end{equation}
		where $\Gamma_{p}:=\int_{\Omega_{p}}\omega(\boldsymbol{x})\mathrm{d}\boldsymbol{x}$ is the circulation contained in the vortical fluid element $\Omega_{p}$, associated with each vortex particle, and $\zeta_{\sigma}(r):=\frac{1}{2\pi\sigma^{2}}\zeta(\frac{r}{\sigma})$ is referred to in the literature as \emph{regularised kernel}, see for example \cite{CottetKoumoutsakos2000,WinckelmansGeneralVortex2004}. The parameter $\sigma$ is the \emph{core size} associated to the size of each vortex particle. Typically $\Gamma_{p}$ is computed by numerical integration using some kind of quadrature rule. The most common one, and used in this work, is to use lowest order Gauss-Lobatto quadrature,
		\begin{equation}
			\Gamma_{p} = \omega(\boldsymbol{x}_{p}) h^{2}\,. \label{eq::blob_circulation_quadrature}
		\end{equation}
		
		Under these conditions, if we substitute \eqref{eq::vorticity_blobs} into \eqref{eq::greens_velocity_particles} we get
		\begin{equation}
			\boldsymbol{u}(\boldsymbol{x}) = -\frac{1}{2\pi}\sum_{p}\frac{g_{\sigma}(\|\boldsymbol{x}-\boldsymbol{x}_{p}\|)}{\|\boldsymbol{x}-\boldsymbol{x}_{p}\|^{2}}(\boldsymbol{x}-\boldsymbol{x}_{p})\times\boldsymbol{e}_{z}\,\Gamma_{p}\,, \label{eq::induced_velocity_field}
		\end{equation}
		where $g_{\sigma}(r):=g(\frac{r}{\sigma})$. Several options exist for the pair of functions $(g,\zeta)$, see for example \cite{CottetKoumoutsakos2000,WinckelmansGeneralVortex2004} for the more standard choices, or \cite{winckelmans1993,Beale1985} for a more extensive discussion of Gaussian kernels and \cite{winckelmansThesis,Rosenhead1930,speckThesis} for arbitrary order algebraic kernels. In this work we use the more standard Gaussian kernels of order two:
		\begin{equation}
			g(\rho) := 1 - e^{-\frac{\rho^{2}}{2}} \qquad \text{and} \qquad \zeta(\rho) := e^{-\frac{\rho^{2}}{2}}\,.
		\end{equation}
		
		For the evaluation of the velocity field by \eqref{eq::induced_velocity_field} we use a Fast Multipole Method for a fast and grid-free computation, see Goude and Engblom, \cite{Goude2012}, and Engblom, \cite{Engblom2011}, for further details.
		
		The convergence proof for vortex particle methods, first introduced by Hald in \cite{Hald1979}, extended simultaneously by Beale and Majda in \cite{Beale1982} and Cottet in \cite{Cottet1982}, is presented in full detail by Cottet and Koumoutsakos in \cite{CottetKoumoutsakos2000}. This proof establishes the mathematical justification for the success of this method. On the other hand, it also establishes an important upper bound on the inter-particle distance, $h$, required for the convergence of the method, namely
		\begin{equation}
			h \leq C\, \sigma ^{1+s} \quad \text{with}\quad s,C \geq 0\,. \label{eq::particle_distance_bound}
		\end{equation}
		Due to the nature of fluid flow, it is known that there exists a time $T$ after which \eqref{eq::particle_distance_bound} is no longer satisfied and the method no longer converges. Several techniques have been introduced to preclude the inter-particle distance from increasing above the bound \eqref{eq::particle_distance_bound}. \emph{Iterative circulation processing}, \cite{Beale1988}, \emph{rezoning}, \cite{Beale1985a,Nordmark1991}, and \emph{remeshing}, \cite{Koumoutsakos1997a}, have been three of the most popular approaches. In this work we opted for remeshing for its computational efficiency. This method introduces an underlying structured mesh with spacing equal to a target inter-particle distance, $h$, satisfying the overlap criterion \eqref{eq::particle_distance_bound}. At regular time steps the vortex particles are reinitialised into this grid. This reinitialisation consists of generating a new set of particles on the regular mesh by means of conservative \emph{interpolatory kernels} with compact support, $W$, according to the expression
		\begin{equation}
			\Gamma^{\text{new}}_{q} = \sum_{p\in Q_{q}} W(x^{\text{new}}_{q}-x^{\text{old}}_{p}, y^{\text{new}}_{q}-y^{\text{old}}_{p})\,\Gamma^{\text{old}}_{p}\,, \label{eq::redistribution_general_expression}
		\end{equation}
		where the particle, $p$, located at $\boldsymbol{x} = [x_{p}, y_{p}]$ has an associated circulation $\Gamma_{p}$ and $Q_{q}$ is the set of indices associated to the particles that lie in the support of $W(x^{\mathrm{new}}_{q}-x,y^{\mathrm{new}}_{q}-y)$. Typically, this higher dimensional interpolatory kernels are constructed as the tensor product of one dimensional ones,
		\[
			W(\xi,\eta) := W(\xi) W(\eta)\,.
		\] 
		Many options exist for the definition of $W$, each with different orders of accuracy and properties, see for example \cite{CottetKoumoutsakos2000, WinckelmansGeneralVortex2004} for a detailed discussion. In this work we use the popular $M'_{4}$ interpolating kernel introduced by Monaghan in \cite{Monaghan1985} and given by the following expression
		\begin{equation}
			M'_{4}(\rho):=\left\{\begin{array}{ll} 1-\frac{5}{2}\rho^{2} + \frac{3}{2}|\rho|^{3} & \text{if } |\rho|<1 \\ \frac{1}{2}\left(2-|\rho|\right)^{2}\left(1-|\rho|\right) & \text{if }1\leq |\rho| < 2 \\ 0 & \text{if } |\rho|\geq 2\,,\end{array}\right. \label{eq::redistribution_m4}
		\end{equation}
		and $W(\xi) = M'_{4}(\frac{\xi}{h})$. This kernel is known for its good combination of smoothness and accuracy, as pointed out by Cottet et al. in \cite{Cottet1999}.
		
		For viscous flows, $\nu>0$, viscosity is not only transported along the velocity field lines but also diffused. Therefore \eqref{eq::vortex_equations_inviscid} is extended into
		\begin{numcases}{}
				\frac{\mathrm{d} \boldsymbol{x}_{p}}{\mathrm{d}t} = \boldsymbol{u}(\boldsymbol{x}_{p},t)\,, \label{eq::convection_part_vortex}\\
				\frac{\mathrm{d}\omega}{\mathrm{d}t} = \nu\nabla^{2}\omega\,. \label{eq::diffusion_part_vortex}
		\end{numcases}
		This characteristic of fluid flow led  to the introduction of the method of \emph{viscous splitting} by Chorin in \cite{Chorin1973}. Instead of solving the coupled system of convection \eqref{eq::convection_part_vortex} and diffusion \eqref{eq::diffusion_part_vortex}, this method decouples the two equations. This is done by introducing two sub-steps. In the first one, particles move along the field lines, according to \eqref{eq::convection_part_vortex}, as discussed before. In the second one, the diffusion effects are computed according to a discrete version of \eqref{eq::diffusion_part_vortex}. The most popular approaches for the discretization of the diffusion term, \eqref{eq::diffusion_part_vortex}, are the \emph{particle strength exchange} (PSE) method, \cite{Degond1989}, the \emph{vortex redistribution method} (VRM), \cite{Shankar1996}, and the \emph{random vortex method} (RVM), \cite{Chorin1973}. For a detailed discussion of the different particle diffusion schemes we recommend the shorter overview by Barba, \cite{Barba2005}, or the more extensive works by Barba, \cite{BarbaPhD}, and Cottet and Koumoutsakos, \cite{CottetKoumoutsakos2000}. Although the VRM has the advantage of maintaining the Lagrangian character of the vortex particle method, since it does not require an underlying mesh, it comes at the cost of solving at each time step an underdetermined system of equations per particle. In order to improve the efficiency of this method, Wee and Ghoniem, \cite{Wee2006}, proposed to redistribute the particle strengths into target particles placed on a regular mesh. Although very efficient, this method imposes restrictive bounds on the size of the time steps, particularly on their minimum size. To circumvent this limitation, we opted to employ the method of Wee and Ghoniem within the context of the Particle-Mesh approach used by Sbalzarini et al. \cite{Sbalzarini2006}. This formulation consists of two steps:
		\begin{enumerate}
			\item Remesh the particles into a regular mesh, using \eqref{eq::redistribution_general_expression} together with the $M'_{4}$ interpolation kernel, \eqref{eq::diffusion_part_vortex}.
			\item Diffuse the vorticity by redistributing the particle's circulation according to
		\end{enumerate} 
		
		\begin{align}
			\Gamma_{i,j}^{k+1} & = \Gamma_{i,j}^{k}\left(1-\frac{4\nu\Delta t}{h^{2}} + \frac{4\nu^{2}\Delta t^{2}}{h^{4}}\right) + \nonumber \\ 
			                               & \quad + \left(\Gamma_{i-1,j}^{k} + \Gamma_{i+1,j}^{k} + \Gamma_{i,j-1}^{k} + \Gamma_{i,j+1}^{k}\right)\left(\frac{\nu\Delta t}{h^{2}} - \frac{2\nu^{2}\Delta t^{2}}{h^{4}}\right) + \nonumber \\
			                               & \quad +\left(\Gamma_{i-1,j-1}^{k} + \Gamma_{i-1,j+1}^{k} + \Gamma_{i+1,j-1}^{k} + \Gamma_{i+1,j+1}^{k}\right)\frac{\nu^{2}\Delta t^{2}}{h^{4}}\,. \label{eq::wee_vrm_mesh}
		\end{align}
		In this expression, as before, $h$ is the inter-particle distance in the regular mesh, the indices $i,j$ correspond to the particles' locations on the mesh, $\boldsymbol{x}_{i,j} = [i h, j h]$, $\Delta t$ is the time step size and the index $k$ corresponds to the time instant $t_{k} = k\Delta t$. This particle circulation update corresponds to the diffusive redistribution scheme presented in \cite{Wee2006} for the special case of regularly distributed particles. It is important to note that this update also satisfies the redistribution formulae given in \cite{Shankar1996} and is equivalent to a PSE method, \cite{CottetKoumoutsakos2000}. Moreover, \eqref{eq::wee_vrm_mesh}, corresponds to a finite difference in time solution of the diffusion equation on a mesh, similar to what is done in \cite{Sbalzarini2006}.
		
		For viscous flows in the presence of solid boundaries, two additional issues need to be taken into account: (i) the enforcement of the no-slip boundary conditions and (ii) the vorticity flux at the solid boundaries. Currently, the standard approach is based on the vorticity boundary conditions pioneered by Lighthill, \cite{Lighthill1963}, and then extended to the vortex particle method by Cottet, \cite{Cottet1994}, and Koumoutsakos et al. in \cite{Koumoutsakos1994,Koumoutsakos1996,Koumoutsakos1995}. For a more detailed discussion of this and other approaches, the reader is directed to the book by Cottet and Koumoutsakos, \cite{CottetKoumoutsakos2000}. On what follows, we give a brief description of this method and its differences in the context of the hybrid flow solver. 
		
\begin{figure}[!ht]
			\center
			\includegraphics[scale=.85]{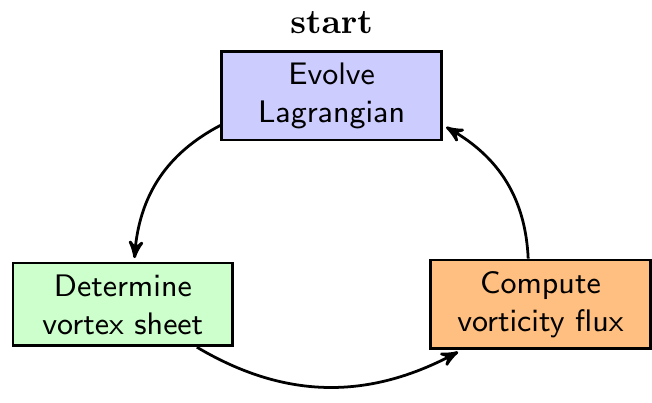}
			\hspace{1.5cm}
			\includegraphics[scale=.85]{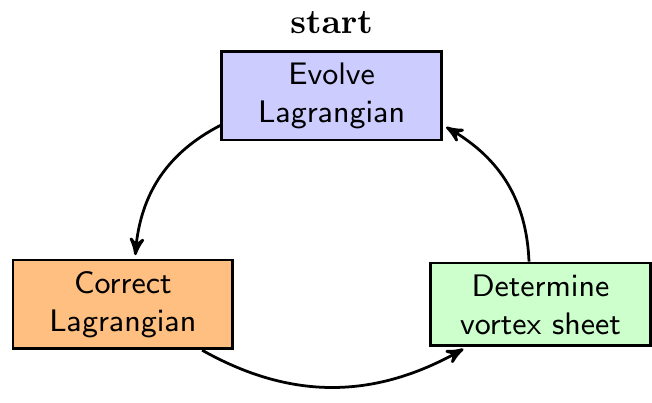}
			\caption{Flowcharts of the viscous Lagrangian vortex particle method (left) and hybrid method (right) in the presence of solid boundaries.}
			\label{fig::lagrangian_flowchart_solid_boundaries}
		\end{figure}

		In summary, this formulation consists of combining the vortex particle method with the boundary element method (vortex sheet panels). Initially, at the beginning of time instant $t_{n}=n\Delta t$, it is assumed that a vorticity field satisfying the no-slip boundary condition exists. The first sub step of this algorithm consists in computing the advection-diffusion evolution of the vortex particles, \eqref{eq::convection_part_vortex} and \eqref{eq::diffusion_part_vortex}, as explained above, therefore without enforcing the boundary conditions. By the end of this first sub step, since no boundary conditions are enforced, a spurious slip velocity exists. In the second sub step the vortex sheet strength that cancels this slip velocity is computed. This in turn is related to the vorticity flux at the solid boundary. This method is represented graphically in \figref{fig::lagrangian_flowchart_solid_boundaries}, left.
		
		In the hybrid method presented in this work, the computation of the vorticity flux is intrinsically taken into account in the Eulerian sub-domain and transferred to the Lagrangian one in the correction sub-step (see \figref{fig::pHyFlow_hybrid_flowchart}). Therefore, in this work, the computation of the vorticity flux is replaced by the correction step and the  computation of the vortex sheet strength, $\gamma(s)$,  is performed at the end in order to enforce the no-slip velocity boundary conditions, as depicted in \figref{fig::lagrangian_flowchart_solid_boundaries}, right.
		
		In this work we focus only on enforcing the no-slip boundary conditions since, due to the \emph{linked boundary conditions}, there is a direct relation between the no-slip and the no-through boundary conditions. This relation is detailed by Koumoutsakos, \cite{KoumoutsakosPhD}, and Shiels, \cite{ShielsPhD}.
		
		The vortex sheet strength is obtained from the following Fredholm equation of the second kind:
		\begin{equation}
			\frac{\gamma(s)}{2} -\frac{1}{2\pi} \oint_{\Sigma_{w}}\frac{\partial}{\partial n}\left[\log(|\boldsymbol{x}(s)-\boldsymbol{x}(s')|)\right] \,\gamma(s')\mathrm{d}s' = \boldsymbol{u}_{slip}\cdot\hat{\boldsymbol{s}}\,,
		\end{equation}
		where $\Sigma_{w}$ is the solid boundary, $\boldsymbol{u}_{slip}$ is the slip velocity and $\hat{\boldsymbol{s}}$ is a unit vector tangent to the solid boundary in the direction of integration. These equations were discretized using the boundary element method (panel method), resulting in:
		\begin{equation}
			\frac{\gamma_{k}(s)}{2} -\frac{1}{2\pi}\sum_{k=1}^{N} \int_{\sigma_{k}}\frac{\partial}{\partial n}\left[\log(|\boldsymbol{x}(s)-\boldsymbol{x}(s')|)\right] \,\gamma_{k}(s')\mathrm{d}s' = \boldsymbol{u}_{slip}\cdot\hat{\boldsymbol{s}}\quad \text{for} \quad k=1,\dots,N\,, \label{eq::panels_vortex_sheet}
		\end{equation}
		where $\gamma_{k}$ is the vortex sheet distribution associated to the panel $k$ and $\sigma_{k}$ is the boundary segment associated to the panel $k$. The drawback of this formulation is that the system of equations \eqref{eq::panels_vortex_sheet} is singular, and therefore allows an infinite number of solutions. An additional equation that limits the solution space to a unique solution is one prescribing the total strength of the vortex sheet:
		\begin{equation}
			\sum_{k=1}^{N}\int_{\sigma_{k}}\gamma_{k}(s)\mathrm{d}s = \Gamma_{\gamma}\,. \label{eq::vortex_sheet_total_circulation}
		\end{equation}
		The total strength, $\Gamma_{\gamma}$, can be obtained from Kelvin's circulation theorem which states that the total circulation must be conserved in time. As will be seen in more detail in \secref{subsection::hybrid_solver}, the total circulation of the vortex sheet can be determined from the coupling between the Eulerian and the Lagrangian sub-domains.
		
		By adding this additional equation, \eqref{eq::vortex_sheet_total_circulation}, the resulting system of equations will be overdetermined since there will be $N$ unknowns but $N+1$ equations. Several options exist to solve this problem, see for example Cottet and Koumoutsakos, \cite{CottetKoumoutsakos2000}. In this work we solve this overdetermined system of equations using a Least-Squares approach.
		
		It is worth to mention a more robust formulation proposed by Koumoutsakos and Leonard, \cite{KoumoutsakosBoundary1993}, consisting of eliminating the first term in the spectral decomposition of the kernel of the system of equations. By doing this, the singularity is removed and the system of equations is well-posed, in the future we intend to explore this route.
		
		By following the steps outlined above, at each time step the velocity field in the Lagrangian domain will be given by:
		\begin{equation}
			\boldsymbol{u} = \boldsymbol{U}_{\infty} + \boldsymbol{u}_{\omega} + \boldsymbol{u}_{\gamma}\,, \label{eq::velocity_field_lagrangian}
		\end{equation}
		where $\boldsymbol{U}_{\infty}$ is the free stream velocity, $\boldsymbol{u}_{\omega}$ is the velocity associated to the vortex particle distribution and $\boldsymbol{u}_{\gamma}$ is the velocity induced by the vortex sheet at the solid boundary.
	
	\subsection{Eulerian solver} \label{subsection::eulerian_solver}
		In this section we introduce the Eulerian flow solver used in the hybrid method presented in this work. The most common approaches to discretize the Navier-Stokes equations are the Finite Volume Method, the Finite Differences Method and the Finite Element method. We have opted for a Finite Element discretization for its geometrical flexibility and $hp$-refinement capabilities, which we intend to use in future developments of this work.
		
		Here, as opposed to \secref{subsection::lagrangian_solver}, we consider the  $(\boldsymbol{u},p)$ formulation of the Navier-Stokes equations in a two dimensional bounded domain, $\Omega \subset \mathbb{R}^{2}$, in the absence of external forces,
		\begin{numcases}{}
				\frac{\partial\boldsymbol{u}}{\partial t} +\left(\boldsymbol{u}\cdot\nabla\right)\boldsymbol{u} -\nabla\cdot\sigma = 0\,, & in $\Omega$\,, \label{eq::ns_momentum_u_p} \\
				\nabla\cdot\boldsymbol{u} = 0\,, & in $\Omega$\,, \label{eq::ns_incompressibility} \\
				\boldsymbol{u}(\boldsymbol{x},t) = \boldsymbol{u}_{0}(\boldsymbol{x})\,, & in $\Omega$ and for $t=0$\,, \\
				\boldsymbol{u}(\boldsymbol{x},t) = \boldsymbol{u}_{b}(\boldsymbol{x},t)\,. & in $\partial\Omega$ and for $t\in]0,T]$\,, \label{eq::velocity_bc}
		\end{numcases}
		where $\boldsymbol{u}$ denotes the velocity, $\boldsymbol{u}_{0}$ corresponds to the initial velocity condition, $\boldsymbol{u}_{b}$ is the velocity at the boundary and $\sigma$ is the Cauchy stress tensor defined as
		\[
			\sigma(\boldsymbol{u},p) := 2\nu\epsilon(\boldsymbol{u}) - p\,\mathsf{I}\,,
		\]
		where $\nu$ is the kinematic viscosity, $p$ is the pressure, $\mathsf{I}$ is the identity matrix and $\epsilon(\boldsymbol{u}):=\frac{1}{2}(\nabla\boldsymbol{u}+\nabla\boldsymbol{u}^{t})$.
		
		To solve the Navier-Stokes equations with the Finite Element Method (FEM) it is necessary to construct the associated weak form and to select the appropriate function spaces for the velocity, $\boldsymbol{u}$, and pressure, $p$, in order to ensure stability.
		
		It is well known that suitable function spaces must satisfy the Ladyzhenskaya-Babu\v{s}ka-Brezzi (LBB) compatibility condition, see Brezzi and Fortin, \cite{brezzi1991mixed}. Such a pair of function spaces is the Taylor-Hood family, \cite{TaylorHood1973}. This family of function spaces consists in $C^{0}$ Lagrange interpolants of polynomial order two ($V_{h}$) for the velocity and polynomial order one ($Q_{h}$) for pressure.
		
		For the construction of the weak form we use the \emph{Incremental Pressure Correction Scheme} (IPCS). This method, introduced by Goda, \cite{Goda1979}, is an improvement to Chorin's projection scheme, \cite{Chorin1968}, and consists of computing a tentative velocity, $\boldsymbol{u}^{\star}$, by neglecting the pressure in \eqref{eq::ns_momentum_u_p} and then correcting it by determining the pressure field that gives rise to a divergence free velocity field. To advance the solution in time we use a forward Euler scheme. This algorithm can be summarised in the following three steps:
		
		\begin{enumerate}
			\item \textbf{Compute the tentative velocity:}  At time instant $t_{n}=n\Delta t$, find the approximate tentative velocity $\boldsymbol{u}_{h}^{\star}\in V_{h}$ such that:
			\begin{equation}
				\begin{split}
		\langle \frac{\boldsymbol{u}_{h}^{\star}-\boldsymbol{u}_{h}^{n-1}}{\Delta t}, \boldsymbol{v} \rangle &+ \langle \boldsymbol{u}_{h}^{n-1}\cdot\nabla\boldsymbol{u}_{h}^{n-1},\boldsymbol{v}\rangle + \langle \sigma(\boldsymbol{u}_{h}^{n-\frac{1}{2}},p_{h}^{n-1}), \epsilon(\boldsymbol{v}) \rangle \quad \\ &\quad+ \langle p_{h}^{n-1}\hat{\boldsymbol{n}},\boldsymbol{v}\rangle_{\partial \Omega} - \langle\nu\hat{\boldsymbol{n}} \cdot (\nabla \boldsymbol{u}_{h}^{n-\frac{1}{2}} )^{t},\boldsymbol{v} \rangle_{\partial \Omega} = 0 \qquad \forall \boldsymbol{v}\in V_{h}\,,
				\end{split}
				\label{eq:tentativeVel}
			\end{equation}
			where $\hat{\boldsymbol{n}}$ is the unit vector normal to the boundary, $\partial\Omega$, and $\boldsymbol{u}_{h}^{n-\frac{1}{2}}:=\frac{\boldsymbol{u}_{h}^{\star}-\boldsymbol{u}_{h}^{n-1}}{2}$. The Dirichlet boundary conditions for $\boldsymbol{u}_{h}$, \eqref{eq::velocity_bc}, are also applied to $\boldsymbol{u}_{h}^{\star}$ in this step. In the hybrid method, these correspond to no-slip conditions at the solid boundaries and at the external boundary of the Eulerian sub-domain, $\Sigma_{d}$ in \figref{fig::pHyFlow_domain_decomposition_overlap_overview}, the velocity is obtained from the Lagrangian sub-domain.
			
			\item \textbf{Compute the pressure:} The pressure is obtained by finding $p^{n}_{h}\in Q_{h}$ such that:
			\begin{equation}
		\langle \nabla p_{h}^n, \nabla q \rangle = \langle \nabla p_{h}^{n-1}, \nabla q\rangle - \frac{1}{\Delta t}\langle \nabla \cdot \boldsymbol{u}_{h}^{\star}, q \rangle \qquad \forall q\in{Q_{h}}\,.
			\label{eq:pressureCorrection}
			\end{equation}
			Where weak homogeneous Neumann boundary conditions are used for the pressure, see Guermond et al. \cite{Guermond2006} for a discussion on pressure boundary conditions in the context of projection schemes.
			\item \textbf{Determine the corrected velocity:} The corrected velocity field is obtained by finding $\boldsymbol{u}_{h}^{n}\in V_{h}$ such that:
			\begin{equation}
				\langle \boldsymbol{u}_{h}^n, \boldsymbol{v}\rangle = \langle \boldsymbol{u}_{h}^{\star},\boldsymbol{v} \rangle - \Delta t \,\langle \nabla(p^{n}_{h} - p_{h}^{n-1}),\boldsymbol{v} \rangle \qquad \boldsymbol{v}\in V_{h}.
				\label{eq:velocityCorrection}	
			\end{equation}
			The Dirichlet boundary conditions for velocity are taken into account at in step 1.
		\end{enumerate}
		
		This solver was implemented using the finite element library FEniCS, see Logg and Wells, \cite{Logg2010}, for more details.

	\subsection{Hybrid solver}\label{subsection::hybrid_solver}
		As seen in the introduction to the hybrid flow solver, \secref{subsection::hybrid_solver_intro}, the main idea behind the hybrid solver discussed in this work is the following. Divide the computational domain into two sub-domains: one Eulerian and one Lagrangian, see \figref{fig::pHyFlow_domain_decomposition_overlap_overview}. The Lagrangian solver covers the whole computational domain and the Eulerian one exists only in the vicinity of solid walls, being capable of capturing near-wall effects and serves as a correction for the vortex particles in that region. This, in turn, results in an evolution algorithm comprised of four steps: \emph{evolve Lagrangian}, \emph{determine Eulerian boundary conditions}, \emph{evolve Eulerian} and \emph{correct Lagrangian}, see \figref{fig::pHyFlow_hybrid_flowchart}.  Therefore, in this section, we will discuss each step of the evolution algorithm of the hybrid method.
		
		\subsubsection{Evolve Lagrangian}
			The first step of the algorithm to advance the hybrid solver solution from time instant $t_{n}$ to $t_{n+1} = t_{n} + \Delta t_{L}$ is the evolution of the Lagrangian sub-domain. We assume that at the start of the algorithm the velocity field in the overlap region is nearly identical in the two sub-domains. Therefore, at the start of this step the velocity field in the Lagrangian sub-domain, \eqref{eq::velocity_field_lagrangian}, satisfies the no-slip boundary conditions. After evolving the Lagrangian sub-domain the vortex particles are remeshed into a regular grid, using \eqref{eq::redistribution_general_expression}, as discussed before.
			
			Under these conditions, the Lagrangian sub-domain is evolved as described in \secref{subsection::lagrangian_solver}. Once more we note that, in the context of the hybrid solver, the evolution of the Lagrangian solver does not explicitly include the generation of vorticity at the solid boundary. For this reason, at the end of this step, the Lagrangian solution will not fully resolve the near-wall region. Nevertheless, as discussed by Daeninck, \cite{daeninckThesis}, and Stock, Gharakhani and Stone, \cite{Stock2010}, this error will not affect the accuracy of the induced velocity field in the outer boundary of the Eulerian sub-domain, $\Sigma_{d}$, \figref{fig::hybrid_interpolation_domain_and_boundaries} right. 
			 
		\subsubsection{Determine Eulerian boundary conditions}
			Once the Lagrangian solution is evolved, it is possible to use \eqref{eq::velocity_field_lagrangian} to compute the velocity field, $\boldsymbol{u}_{b}$, at the outer boundary of the Eulerian sub-domain, $\Sigma_{d},$ at time instant $t_{n+1} = t_{n} + \Delta t_{L}$. If the velocity field at the outer boundary is required at any other time instant $t\in [t_{n},t_{n+1}]$ a linear interpolation in time is used,
			\begin{equation}
				\boldsymbol{u}_{b}(t) = \boldsymbol{u}_{b}(t_{n}) + \frac{t-t_{n}}{\Delta t_{L}}\left[\boldsymbol{u}_{b}(t_{n+1}) - \boldsymbol{u}_{b}(t_{n})\right] \qquad t\in[t_{n},t_{n+1}]\,.\label{eq::velocity_boundary_intermediate}
			\end{equation}
			
			Note that higher order interpolation could be used.
			
		\subsubsection{Evolve Eulerian}
\begin{figure}[!ht]
			\center
			\includegraphics[width=0.5\textwidth]{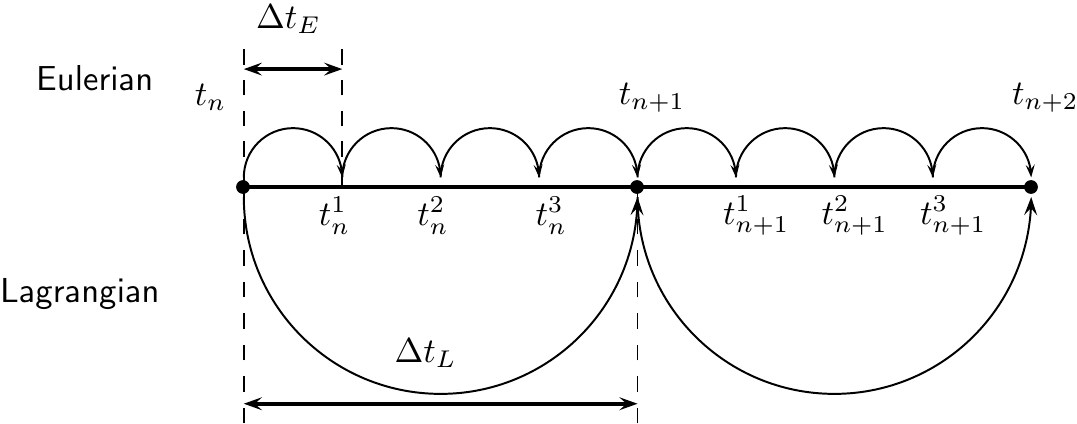}
			\caption{Eulerian multi-stepping to match the Lagrangian time step. The figure shows $\Delta t_{L} = k_{E}\Delta t_{E}$ with $(k_{E} = 4)$ Eulerian sub-steps to time march from $t_{n}$ to $t_{n+1}$.}
			\label{fig::hybrid_eulerian_multistep}
		\end{figure}

			Due to the different nature of the solvers in the two sub-domains, different time step constraints are required. Typically, the time step required in the Lagrangian sub-domain, $\Delta t_{L}$, will be larger than the one in the Eulerian sub-domain, $\Delta t_{E}\leq \Delta t_{L}$. For this reason, and in order to improve the computational efficiency, each solver is allowed to advance in time according to different time steps, subject to: $\Delta t_{L} = k_{E}\Delta t_{E}$ with $k_{E}\in\mathbb{N}$, see \figref{fig::hybrid_eulerian_multistep}.
			
			Therefore, to advance the Eulerian sub-domain from time instant $t_{n}$ to $t_{n+1}$, $k_{E}$ Eulerian evolution sub-steps associated to the time instants
			\[
				t_{n}^{k} := t_{n} + k\Delta t_{E} \qquad k=0,\dots,k_{E}
			\]
			are performed. The velocity boundary conditions at the intermediate time instants $t_{n}^{k}$ are given by \eqref{eq::velocity_boundary_intermediate}. These Eulerian evolution sub-steps are performed according to what was described in \secref{subsection::eulerian_solver}.
			
			\FloatBarrier
			
		\subsubsection{Correct Lagrangian}
			The final step of the algorithm for evolution of the hybrid solver is the correction step. As mentioned before, in this step the vorticity field computed in the Eulerian sub-domain is transferred to the Lagrangian one. This step replaces the vorticity flux step in pure Lagrangian vortex particle solvers. This step is fundamental for an accurate coupling between the two sub-domains, therefore it is crucial that this transfer of vorticity satisfies conservation of circulation.
			
\begin{figure}[!ht]
        				\center
        				\includegraphics[width=0.3\textwidth]{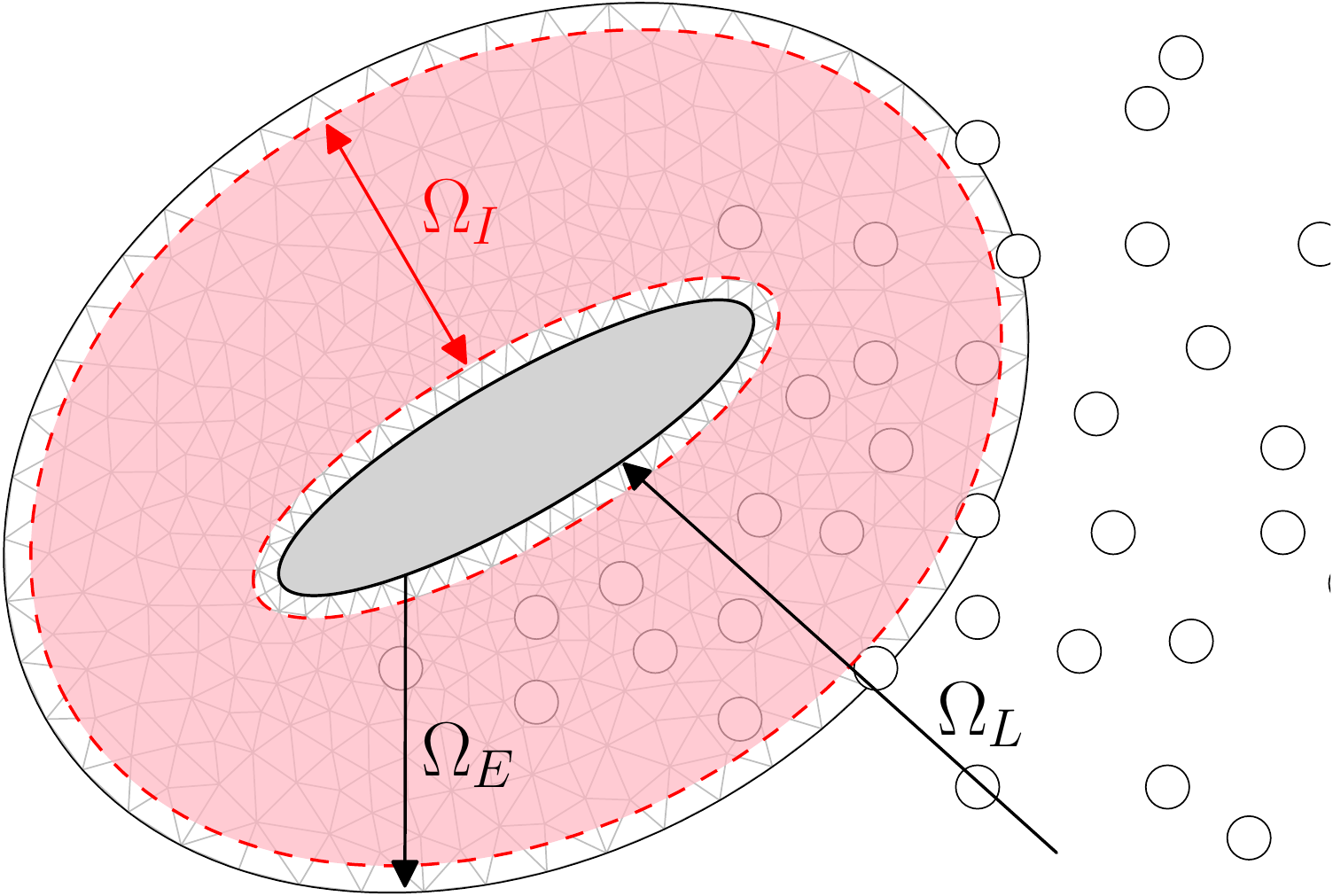}
        				\hspace{1.5cm}
        				\includegraphics[width=.3\textwidth]{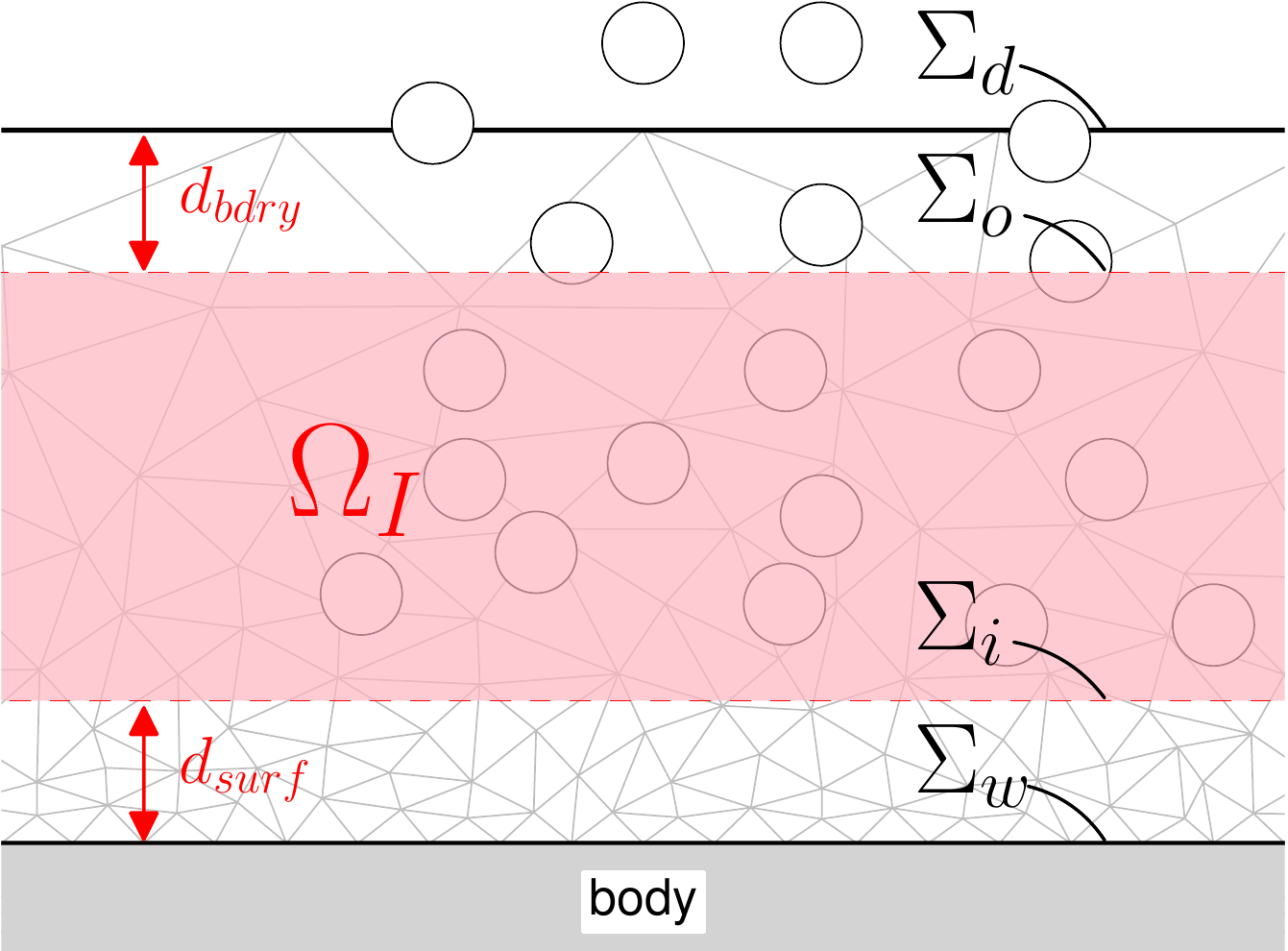}
        				\caption{Definition of the Eulerian, $\Omega_{E}$, and Lagrangian, $\Omega_{L}$, sub-domains and interpolation domain $\Omega_{I}$ where the Lagrangian solution is corrected, with boundary $\partial\Omega_{I}=\Sigma_{i}\cup\Sigma_{o}$. Left: Representation of the whole domain. Right: Detail of the sub-domains close to the boundary.}
        				\label{fig::hybrid_interpolation_domain_and_boundaries}
        			\end{figure}

			In this work, as in the work of Stock, Gharakhani and Stone, \cite{Stock2010}, the interpolation (or correction) region, $\Omega_{I}$, is a subset of the Eulerian sub-domain, $\Omega_{I}\subset \Omega_{E}$, see \figref{fig::hybrid_interpolation_domain_and_boundaries}, left. The interpolation region is defined such that its outer boundary, $\Sigma_{o}$, is at a distance $d_{bdry}$ of the outer boundary of the Eulerian domain, $\Sigma_{d}$, and its inner boundary, $\Sigma_{i}$, is at a distance $d_{surf}$ of the solid wall, $\Sigma_{w}$, see \figref{fig::hybrid_interpolation_domain_and_boundaries}, right. This is done firstly to ensure that the high gradients of vorticity near the solid boundary do not introduce interpolation errors in the Lagrangian domain and secondly to minimise the errors in the Lagrangian sub-domain due to small discrepancies in the boundary conditions of the Eulerian sub-domain.
		
			The algorithm used to correct the Lagrangian solution in the vicinity of the solid walls is summarised as: \emph{generate new particles}, \emph{assign circulations} and \emph{ensure no-slip conditions}.
			
\begin{figure}[!ht]
			\center
			\includegraphics[width=0.25\textwidth]{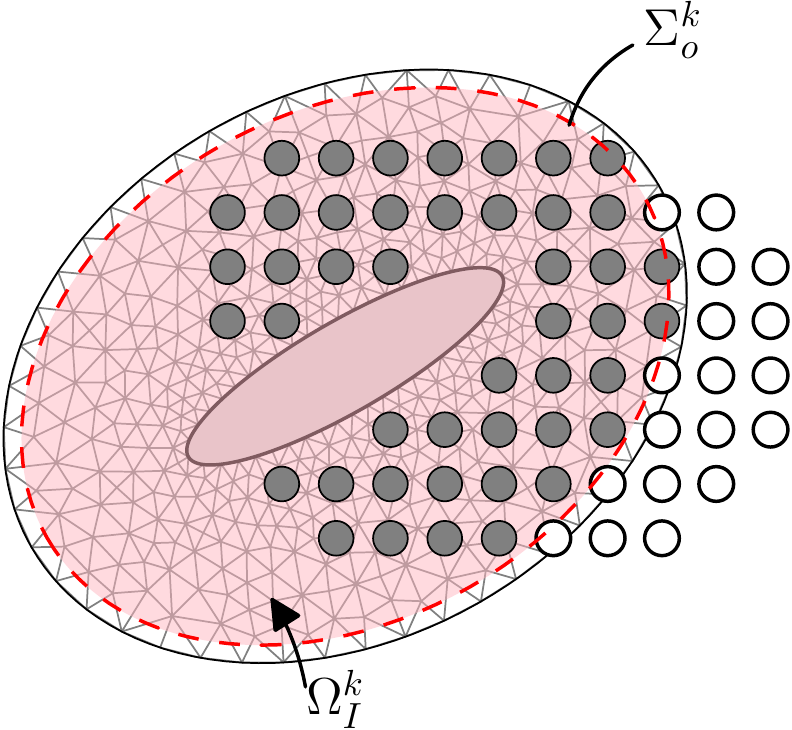}
			\hspace{1.5cm}
			\includegraphics[width=.25\textwidth]{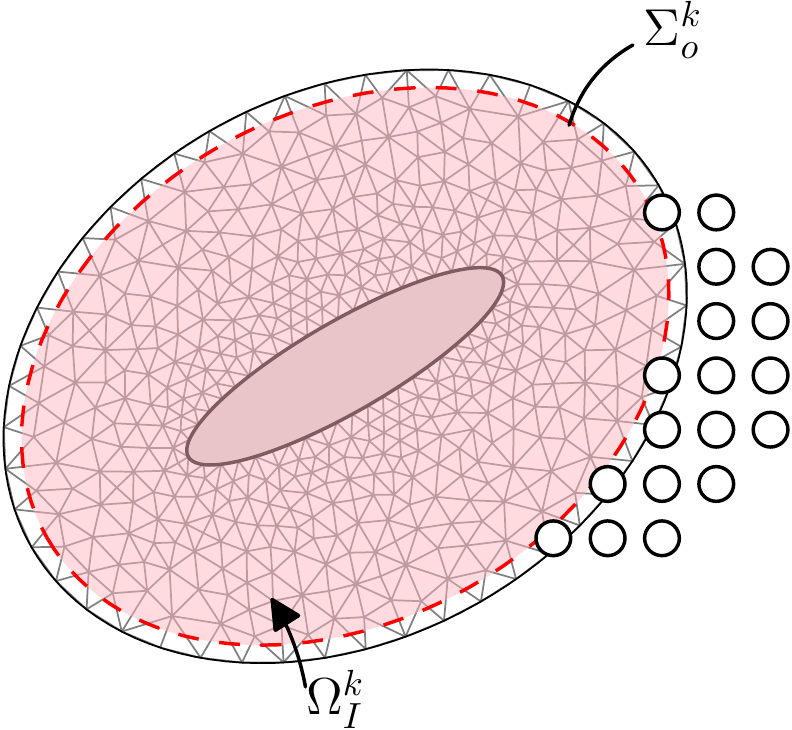}
			\hspace{1.5cm}
			\includegraphics[width=.25\textwidth]{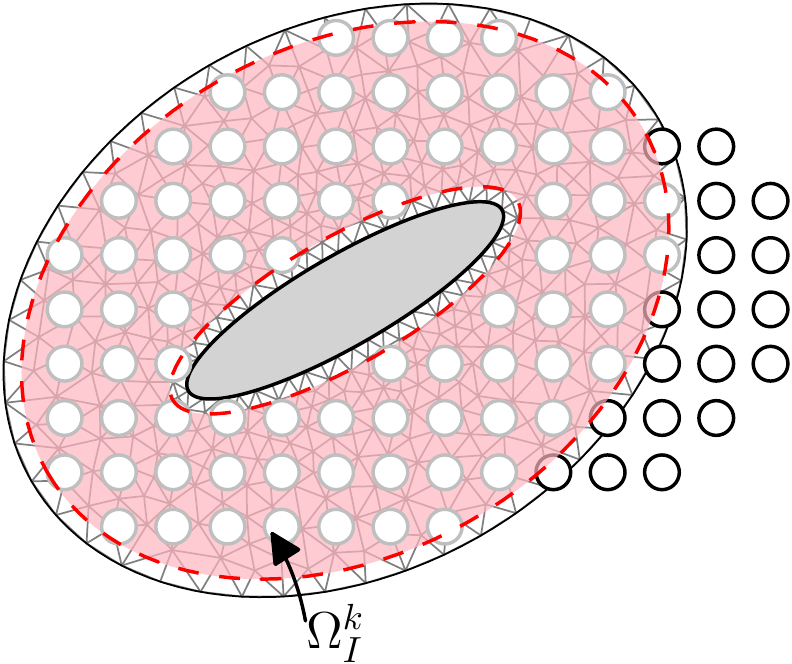}
			\caption{Correction of Lagrangian solution in the vicinity of a solid object $k$: generation of new particles. Left: Particles in the region enclosed by the boundary $\Sigma_{o}^{k}$ are selected. centre: These particles are removed. Right: The interpolation region, $\Omega_{I}$, is covered with new particles with zero circulation and regularly distributed.}
			\label{fig::hybrid_generation_of_new_particles}
		\end{figure}

			The first step in the correction of the Lagrangian solution is the generation of new particles in the interpolation region, $\Omega_{I}$. We start by identifying the particles that lie inside the region bounded by the outer boundary of the interpolation region, $\Sigma_{o}$, red region in \figref{fig::hybrid_generation_of_new_particles} left. These particles are then removed, \figref{fig::hybrid_generation_of_new_particles} centre. We do not use the region $\Omega_{I}$ as a selection region since, during the evolution of the Lagrangian solution, some particles may move to the region bounded by $\Sigma_{i}$ and $\Sigma_{w}$ and therefore would not be removed. Finally, the interpolation region is filled with zero circulation particles regularly distributed, \figref{fig::hybrid_generation_of_new_particles} right.
			
		The second step corresponds to the assignment of circulation to the newly generated vortex particles. This follows the same procedure as outlined in \secref{subsection::lagrangian_solver}, namely \eqref{eq::blob_circulation_quadrature}. Now, the vorticity field is the one obtained from the Eulerian finite element solver, $\omega_{h}(\boldsymbol{x},t)=\nabla\times\boldsymbol{u}_{h}(\boldsymbol{x},t)$.  In general the Eulerian sub-domain will be discretized with an unstructured grid. Evaluating the discrete vorticity field at the location of each new vortex particle is computationally expensive when unstructured grids are used. For this reason, we use an intermediate fixed structured grid. The vorticity field is initially interpolated into the structured grid and then from the structured grid into the vortex particles.	
		
\begin{figure}[!ht]
			\center
			\includegraphics[width=0.55\textwidth]{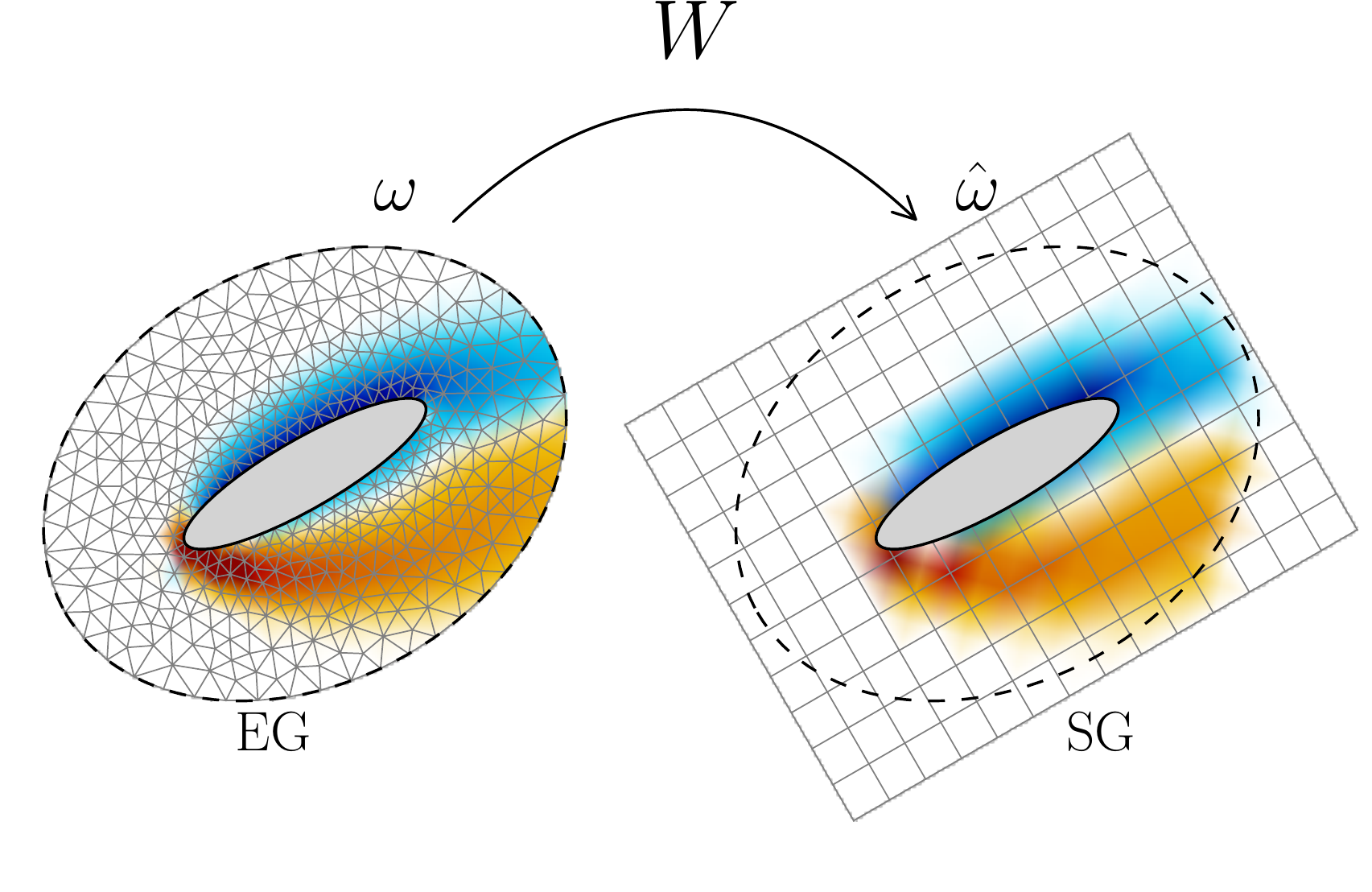}
			\caption{Interpolation of the vorticity $\omega$ from the unstructured Eulerian grid into a fixed structured grid.}
			\label{fig::hybrid_interpolation_unstructured_structured}
		\end{figure}

\begin{figure}[!ht]
			\center
			\includegraphics[width=0.3\textwidth]{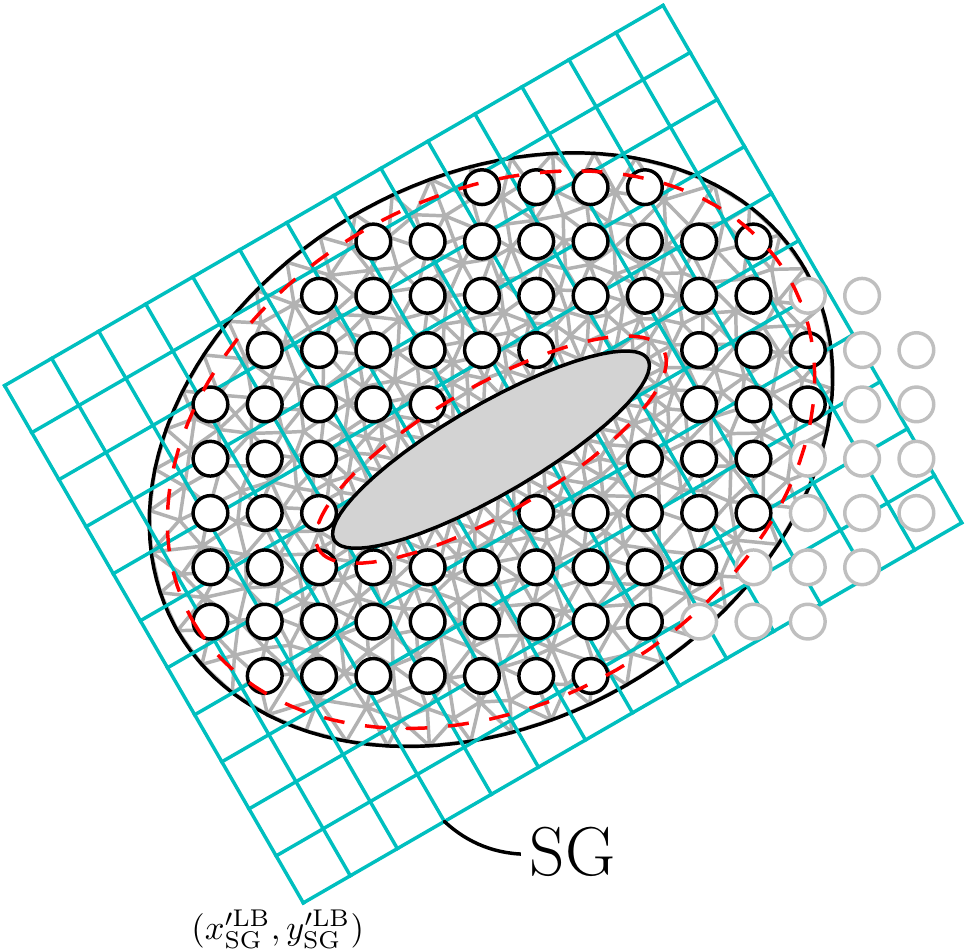}
			\hspace{1cm}
			\includegraphics[width=.3\textwidth]{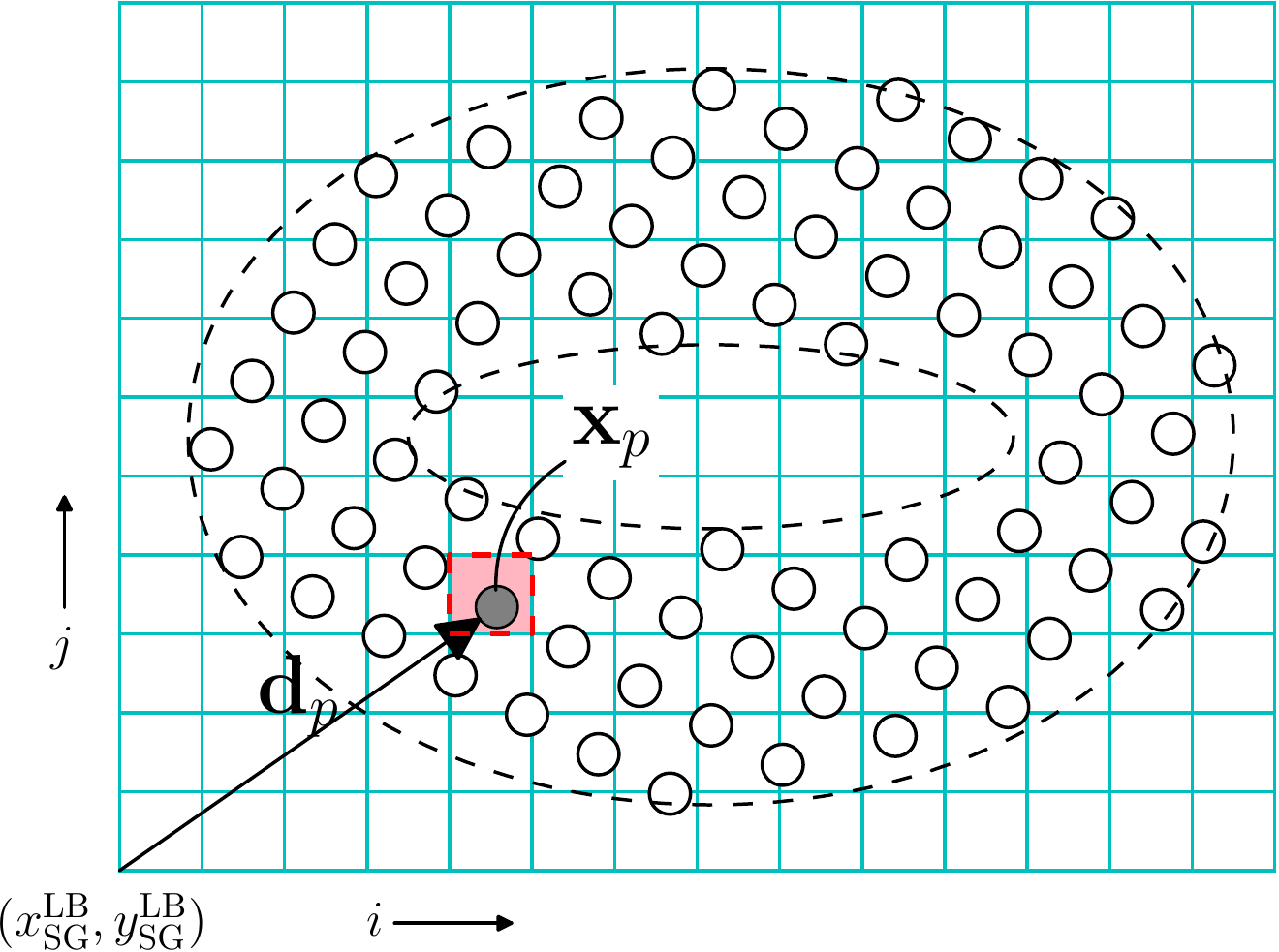}
			\hspace{1cm}
			\includegraphics[width=.25\textwidth]{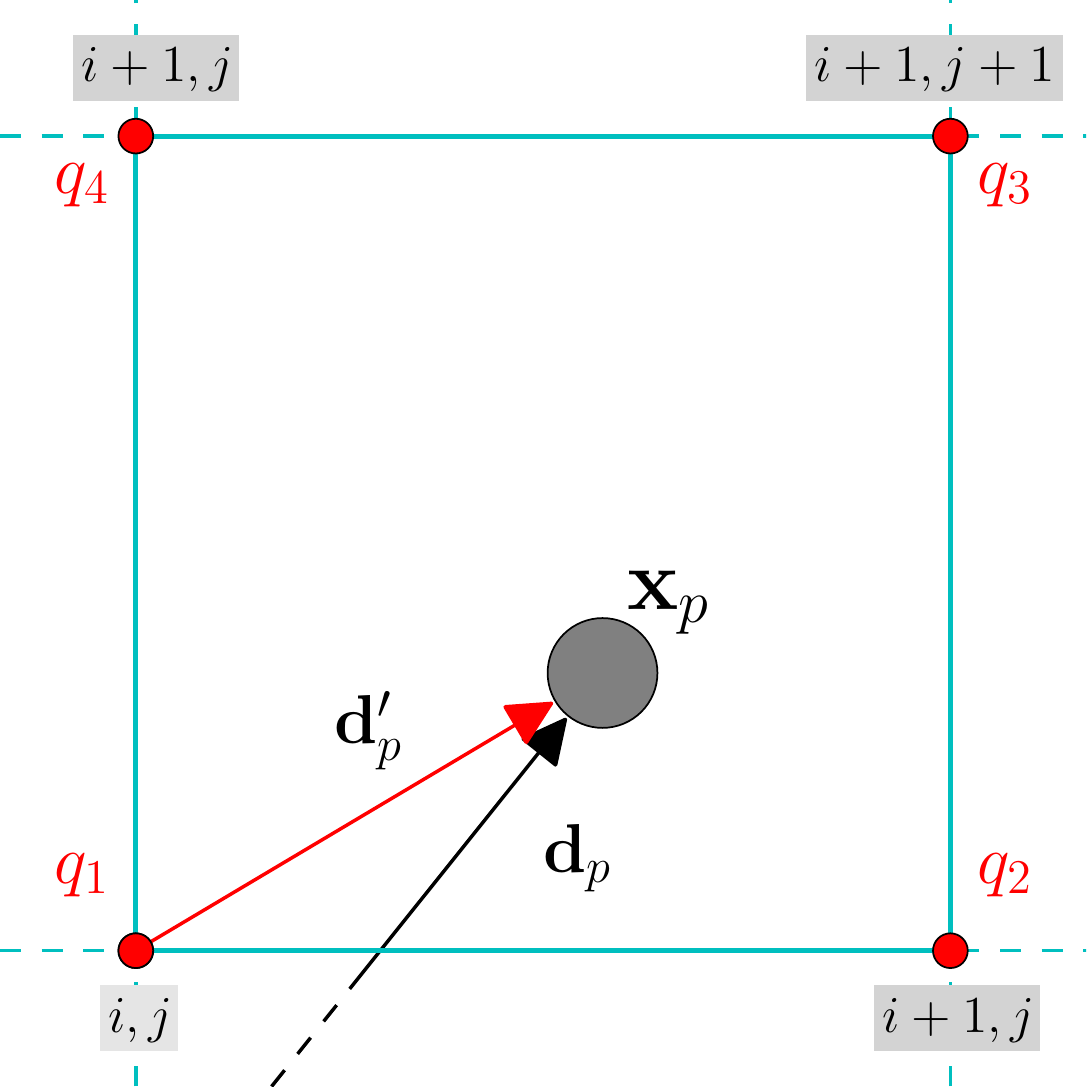}
			\caption{Algorithm to assign circulation strengths to the particles from the vorticity distribution in the structured grid. Left: Vortex particles overlapped to the structured grid and Eulerian unstructured grid. centre: Vortex particle, $p$, in the reference frame of the structured grid is represented by $\boldsymbol{d}_{p}$. Right: Vortex particles in the reference frame of a cell of the structured grid are represented by $\boldsymbol{d}'_{p}$.}
			\label{fig::hybrid_algorithm_circulation_from_vorticity_in_structured_grid}
		\end{figure}

		The structured grid is constructed such that is covers the unstructured one, see \figref{fig::hybrid_algorithm_circulation_from_vorticity_in_structured_grid} left. Since the two grids are fixed to each other, the weights, $W_{ij}$, that define the matrix $W$ representing the interpolation of vorticity from the unstructured grid, $\omega_{h}$, onto the structured one, $\hat{\omega}_{h}$, are constant throughout the simulation and therefore need to be computed only once, reducing considerably the computational cost. We can then compute the vorticity in the structured grid, see \figref{fig::hybrid_interpolation_unstructured_structured}, by using
		\[
			\hat{\omega}_{h,i} = \sum_{i}W_{ij}\omega_{h,j}\,.
		\]
		
		The computation of the vorticity field at the particles' locations is done by first transforming their coordinates into the reference frame of the structured grid, see \figref{fig::hybrid_algorithm_circulation_from_vorticity_in_structured_grid} middle. This allows a fast identification of the cell containing each vortex particle. The computation of the vorticity at each particle is obtained by linear interpolation,
		\[
			\hat{\omega}_{h}(\boldsymbol{x}_{p}) = \sum_{q=1}^{4}\hat{W}_{pq} \,\hat{\omega}_{h,q}\,, 
		\]
		where $q\in\{1,2,3,4\}$ denotes the four nodes of the cell in which the particle $p$ is located, as shown in \figref{fig::hybrid_algorithm_circulation_from_vorticity_in_structured_grid} right, $\hat{\omega}_{h,q}$ represents the associated values of vorticity in those vertices and $\hat{W}_{pq}$ are the weights of the bilinear interpolation given as usual by,
		\begin{align}
			W_{p1} &= \frac{(d'_{x,p}-h)(d'_{y,p}-h)}{h^{2}}\,, \nonumber \\
			W_{p2} &= -\frac{d'_{x,p}(d'_{y,p}-h)}{h^{2}}\,, \nonumber \\
			W_{p3} &= \frac{d'_{x,p}d'_{y,p}}{h^{2}}\,,\nonumber \\
			W_{p4} &= -\frac{d'_{y,p}(d'_{x,p}-h)}{h^{2}}\,, \nonumber 
		\end{align}
		where $\boldsymbol{d}'_{p} = [d'_{x,p},d'_{y,p}]$ is the coordinates of the particle $p$ in the reference frame of the cell, see \figref{fig::hybrid_algorithm_circulation_from_vorticity_in_structured_grid} right, and $h$ is the cell spacing of the structured grid in both $x$ and $y$ directions. We use the same grid spacing in the structured grid as in the remeshing of the vortex particles. Finally, the circulations are assigned to each particle, using \eqref{eq::blob_circulation_quadrature}, see \figref{fig::hybrid_algorithm_circulation_from_vorticity_in_structured_grid}.
		
		\begin{figure}[!ht]
			\center
			\includegraphics[width=0.55\textwidth]{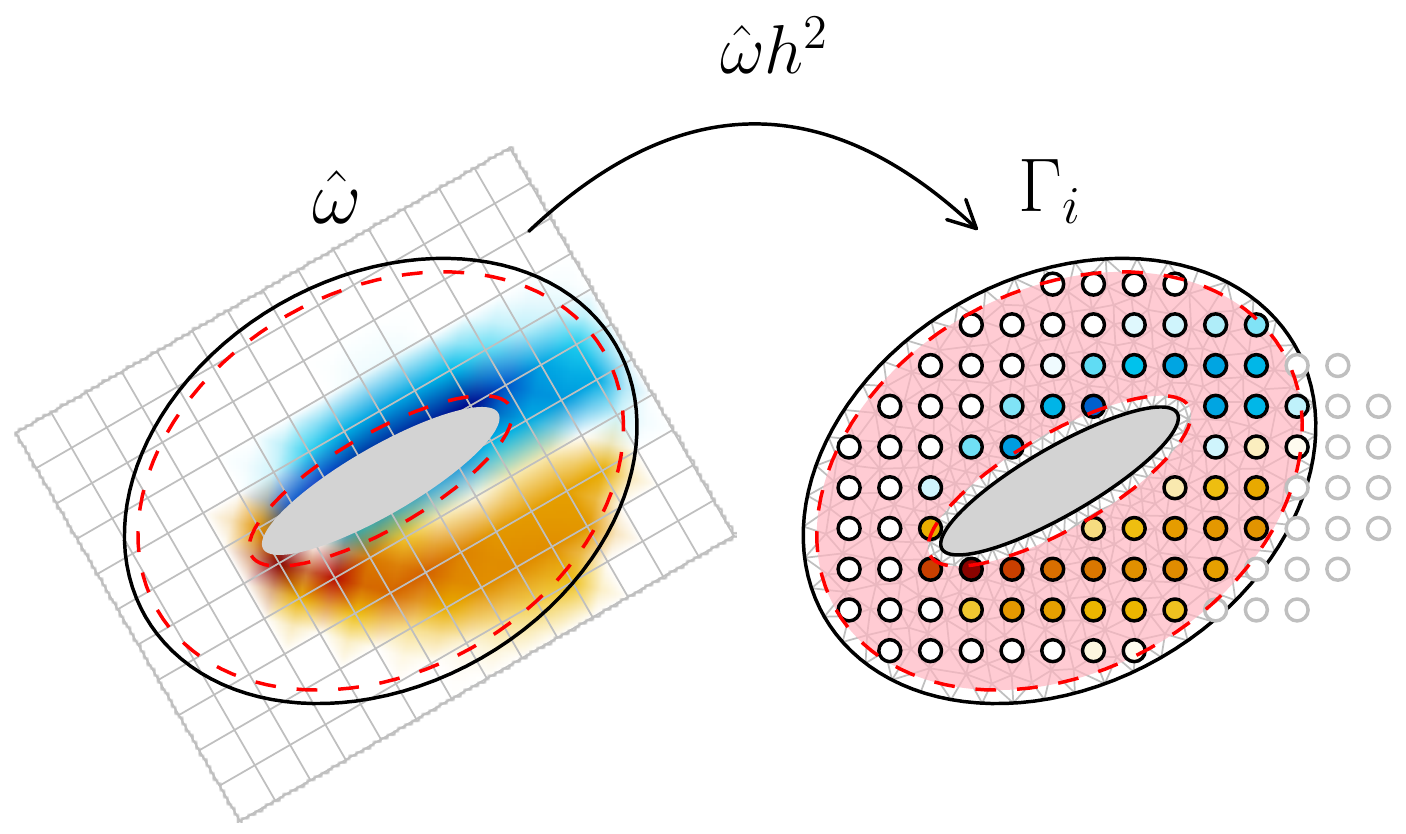}
			\caption{Assignment of circulation into the new blobs in the Lagrangian domain.}
			\label{fig::hybrid_assignment_circulation_blobs}
		\end{figure}
		
		The third and final step consists in ensuring the no-slip boundary condition since, at this stage, there exists a slip velocity at the solid walls. As discussed in \secref{subsection::lagrangian_solver}, we need to compute the new vortex sheet strength that satisfies the no-slip boundary conditions. Because we use a hybrid approach, we do not need to compute the vorticity flux at the boundary, since the assignment of circulations to the new particles automatically takes this into account. It was also seen in \secref{subsection::lagrangian_solver} that the new vortex sheet can be obtained by solving a discrete version of the Fredholm integral equation of the second kind, \eqref{eq::panels_vortex_sheet} and the singularity of this equation circumvented by prescribing the total circulation of the vortex sheet, $\Gamma_{\gamma}$, as an additional equation and solving the overdetermined system in a Least-Squares sense. For bodies rotating at a constant speed, Kelvin's circulation theorem states that the total circulation of our system remains constant in time. Therefore we have:
		\begin{equation}
			\Gamma_{\gamma}(t) + \Gamma_{L}(t) = \mathrm{constant}\,,
		\end{equation}
		where $\Gamma_{L}:=\sum_{p}\Gamma_{p}$, is the total circulation contained in the vortex particles. For the case of bodies with a time varying rotation a similar expression can be deduced, see \cite{CottetKoumoutsakos2000}. Since the Eulerian and the Lagrangian solvers must coincide in the correction region, $\Omega_{I}$, we have that:
		\begin{equation}
			\left.\Gamma_{L}(t)\right|_{\Omega_{I}} = \int_{\Omega_{I}}\omega_{h}(\boldsymbol{x},t)\,\mathrm{d}\boldsymbol{x}\,, \label{eq::circulation_conservation_correction_region}
		\end{equation}
		where $\left.\Gamma_{L}(t)\right|_{\Omega_{I}}$ denotes the expected total circulation contained in the vortex particles that lie in the correction region. Since the step where the circulation is assigned to the particles in the correction region relies on an approximate quadrature formula, \eqref{eq::blob_circulation_quadrature}, there will exist an error, $\delta\Gamma_{L}$, therefore:
		\[
			\left.\Gamma_{L}(t)\right|_{\Omega_{I}} = \left.\tilde{\Gamma}_{L}(t)\right|_{\Omega_{I}}+\delta\Gamma_{L} (t) = \int_{\Omega_{I}}\omega_{h}(\boldsymbol{x},t)\,\mathrm{d}\boldsymbol{x}\,,
		\]
		where  $\left.\tilde{\Gamma}_{L}(t)\right|_{\Omega_{I}}$ is the actual total circulation of the vortex particles that lie in $\Omega_{I}$. In order to exactly conserve circulation, \eqref{eq::circulation_conservation_correction_region}, the circulation of the newly generated vortex particles is corrected using:
		\[
			\Gamma^{\mathrm{corrected}}_{p} = \Gamma_{p} + \frac{\delta\Gamma_{L}}{N_{\Omega_{I}}}\,,
		\]
		where $N_{\Omega_{I}}$ is the number of particles that lie inside the correction region, $\Omega_{I}$.
		
		To finally compute the total circulation of the vortex sheet notice that:
		\[
			\Gamma_{\gamma}(t) + \left.\Gamma_{L}(t)\right|_{\Omega_{E}} = \int_{\Omega_{E}}\omega_{h}(\boldsymbol{x},t)\,\mathrm{d}\boldsymbol{x}\,.
		\]
		Therefore the total circulation of the vortex sheet can be simply obtained by the expression:
		\[
		\Gamma_{\gamma}(t) = \int_{\Omega_{E}}\omega_{h}(\boldsymbol{x},t)\,\mathrm{d}\boldsymbol{x} - \left.\Gamma_{L}(t)\right|_{\Omega_{E}}\,.
		\]
		It is important to note that if the Lagrangian solver could fully capture the boundary layer, then the total circulation of the vortex sheet, $\Gamma_{\gamma}$, would be zero, although each vortex panel would not be zero. We see, then, that the vortex sheet plays not only the fundamental role of ensuring the no-slip boundary conditions, but also of representing the boundary layer. By following this approach, the hybrid solver exactly conserves circulation, greatly improving the results, see \cite{Manickathan2015} for further details. 
		
		\FloatBarrier
		
\section{Numerical benchmark cases} \label{section::numerical_benchmark_cases}
	On what follows we present the results for the solution of different flow problem using the hybrid solver presented in this work. In order to adequately compare the full FEM solution to the hybrid solution we have constructed the meshes such that they are identical in the Eulerian sub-domain.
	
	\subsection{Dipole in unbounded flow} \label{subsection::dipole_in_unbounded_flow}
\begin{figure}[!h]
			\center
			\includegraphics[trim=0cm 2.5cm 0cm 2.5cm, clip, width=0.65\textwidth]{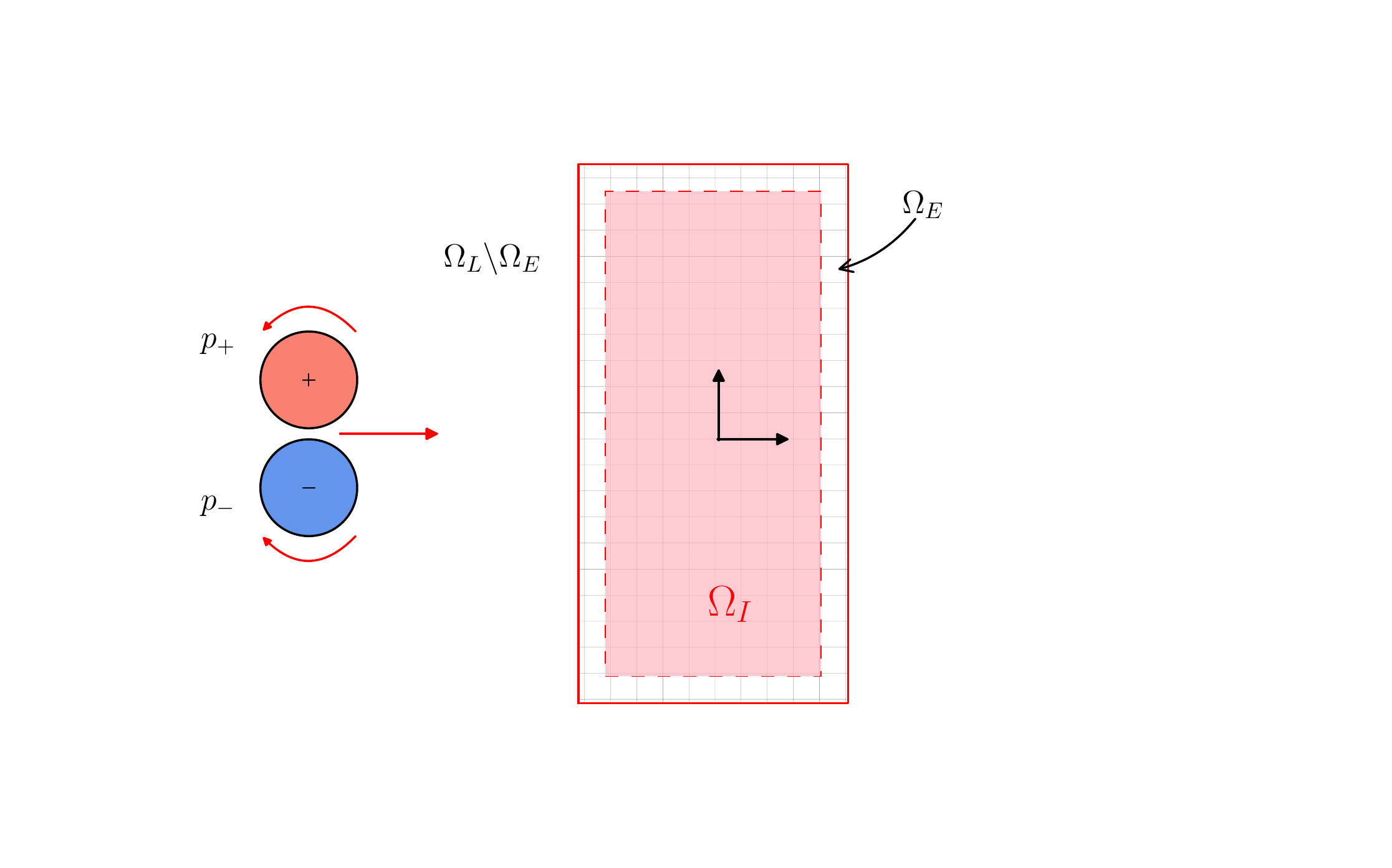}
			\caption{The domain decomposition for the Clercx-Bruneau dipole convection problem, with the positive pole located at $p_{+}=(x_1,y_1) = (-1,0.1)$ and negative pole located at $p_{-}=(x_2,y_2)=(-1,-0.1)$. (\textit{Not to scale})}
			\label{fig::cb_convection_schematics}
		\end{figure}

		The first test case is the evolution of a vortex dipole in an unbounded domain. With this test case we intend to investigate how the flow solution produced by the hybrid solver is perturbed as it traverses the Eulerian sub-domain. In order to do this, we used as initial condition the Clercx-Bruneau dipole, \cite{Clercx2006}, with a positive monopole at $(x_{1}, y_{1}) = (-1.0, 0.1)$ and negative monopole at $(x_{2}, y_{2}) = (-1.0,-0.1)$, each having a core radius $R = 0.1$ and characteristic vorticity magnitude $\omega_{e} = 299.528385375226$ as given by Renac et. al., \cite{Renac2013}, 
		\begin{equation}
			\omega(x,y,0) = \omega_{e}\left(1-\frac{r_{1}^{2}}{R^{2}}\right) e^{-\frac{r_{1}^{2}}{R^{2}}} - \omega_{e}\left(1-\frac{r_{2}^{2}}{R^{2}}\right) e^{-\frac{r_{2}^{2}}{R^{2}}}\,,
		\end{equation}
		where $r_{i}^{2}=\left(x-x_{i}\right)^{2}+\left(y-y_{i}\right)^{2}$. The Eulerian sub-domain is defined as $\Omega_{E}=[-0.25,0.25]\times [-0.5,0.5]$, meaning that the dipole is initially placed to its left, as depicted in \figref{fig::cb_convection_schematics}.
		
\begin{figure}[!h]
			\centering
			\includegraphics[width=0.65\textwidth]{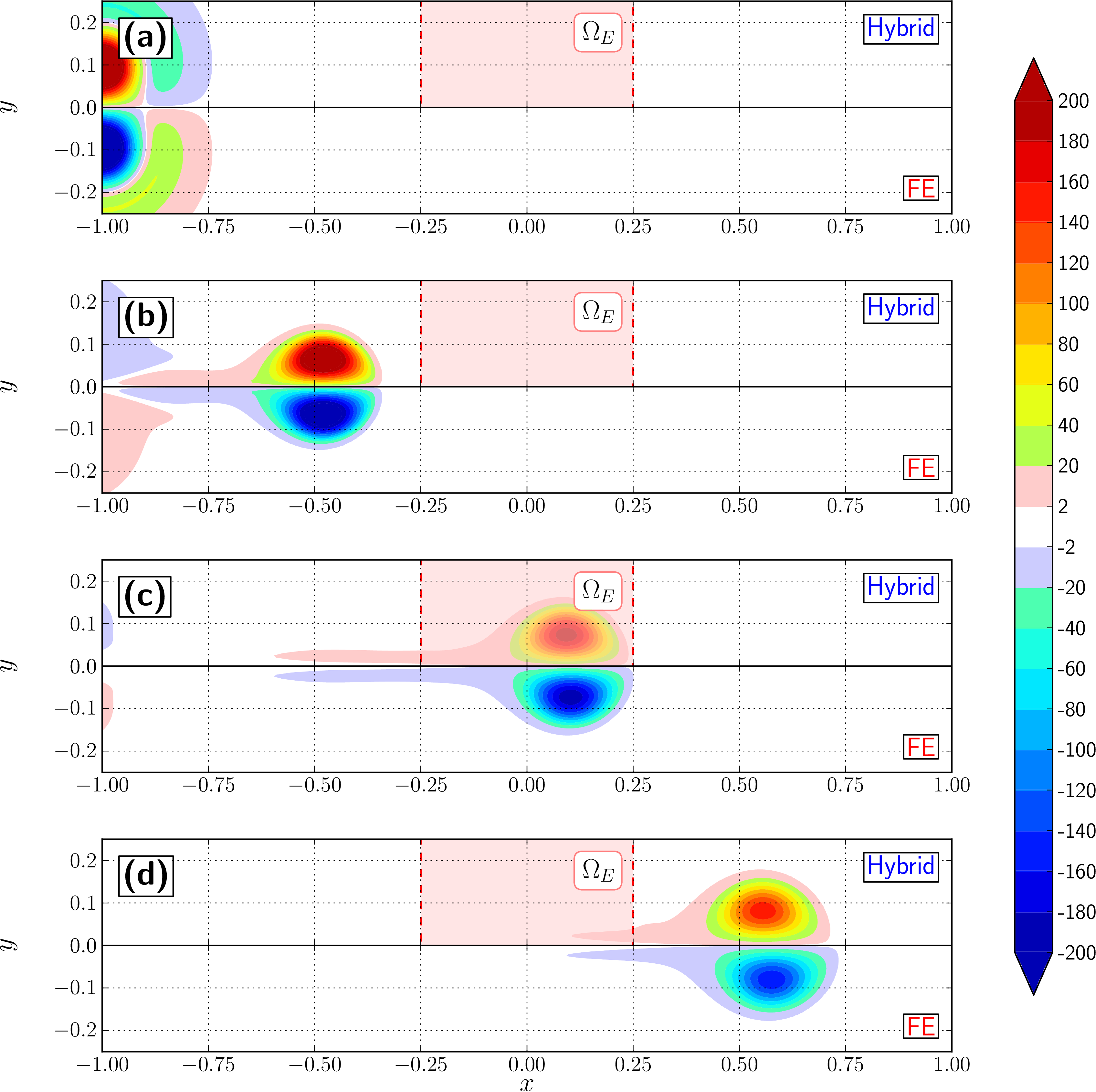}
			\caption{Plot of the Clercx-Bruneau dipole at $t=[0,0.2,0.4,0.7]$. The figure compares the hybrid simulation (top halves) against the FE only simulation (bottom halves).}
			\label{fig:hybrid_doubleMonopoleConvection_contourfPlots}
		\end{figure}

\begin{figure}[!ht]
			\centering
             		\includegraphics[width=0.35\textwidth]{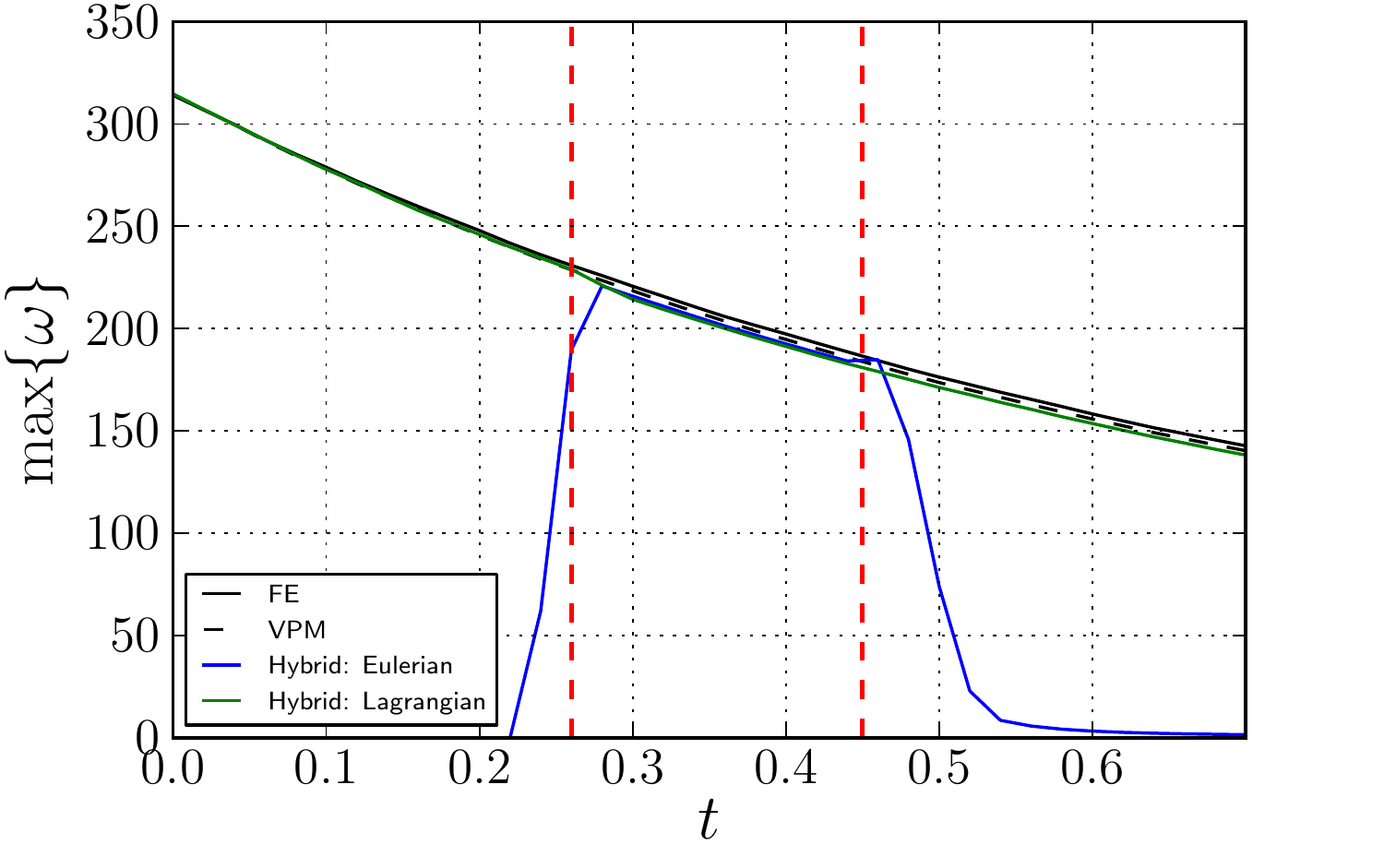}
             		\includegraphics[width=0.35\textwidth]{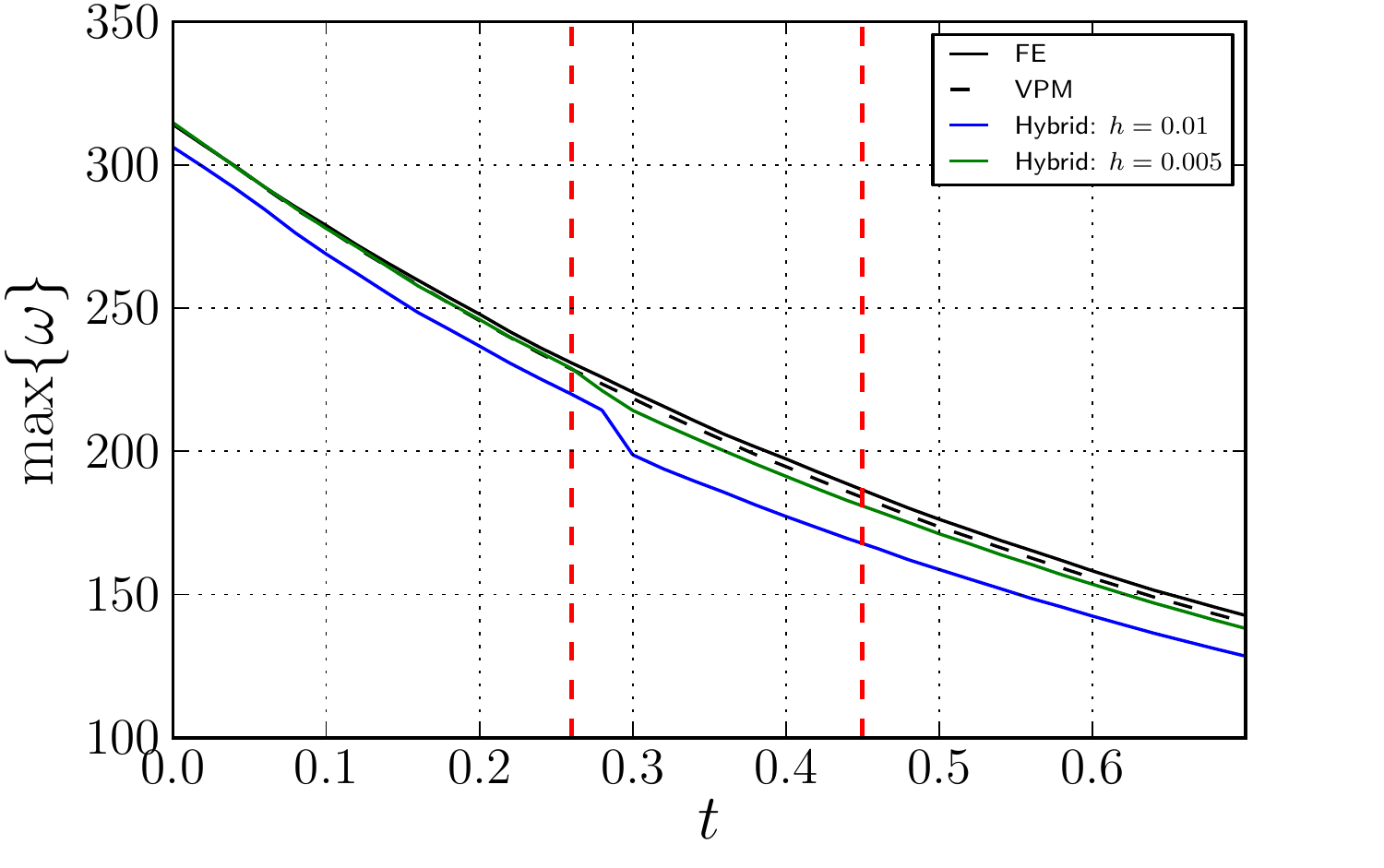}
              		\caption{Evolution of the maximum of vorticity, $\max\{\omega\}$, from $t=0$ to $t=0.7$. The solutions are compared to FE [---, solid black] and  VPM [- -, dashed black] benchmark simulations with characteristics identical to the Eulerian and Lagrangian components of the hybrid simulation. The figure depicts (left) the maximum vorticity in the Eulerian and Lagrangian sub-domains of the hybrid method for a blob spacing $h=0.005$, and (right) the maximum vorticity of the hybrid method with nominal blob spacing $h=0.01$ and $h=0.005$.}
     			\label{fig:hybrid_dipoleConvection_comparison_wMax}
     		\end{figure}

		We ran a simulation with the hybrid solver, using the parameters presented in \tabref{tab:parameters_dipole_unbounded}, and compared it to the results obtained with the FE solver. A contour plot of vorticity, comparing the hybrid results to the FE ones is shown in \figref{fig:hybrid_doubleMonopoleConvection_contourfPlots}. In \figref{fig:hybrid_dipoleConvection_comparison_wMax}, left, we plot the time evolution of the maximum of vorticity, showing the transfer of information between the two sub-domains of the hybrid solver. 
		 
		In \figref{fig:hybrid_dipoleConvection_comparison_wMax}, right, we compare the evolution of the vorticity maximum for different particle sizes, $h=\num{5e-3}$ and $h=\num{1e-2}$. As can be seen, for larger particle sizes it is possible to notice a strong decrease in the vorticity maximum as the dipole enters the Eulerian sub-domain. This produces a small reduction in the mean propagation speed of the dipole, as can be seen in \figref{fig:hybrid_doubleMonopoleConvection_contourfPlots}d where the hybrid solution is slightly lagging behind the FE one.
		
		\begin{table}[htbp]
			\centering
		   	\caption{Summary of the parameters used in the hybrid simulation of the Clercx-Bruneau dipole convection problem.}
		   	\label{tab:parameters_dipole_unbounded}
		   	\begin{tabular}{lcll} 
				\hline\hline
		   		Parameters 					& Value 	& Unit					& Description \\ \hline
				$\nu$						& \num{1.6e-3} & \si{kg.s^{-1}.m^{-1}} & Kinematic viscosity\\
				$\lambda$					& 1 & - & Overlap ratio\\
				$h$							& \num{5e-3} & \si{m} & Nominal blob spacing\\
				$h_{grid}$ 					& $\num{7e-3}$ & \si{m}	& FE cell diameter \\
				$ N_{\mathrm{cells}}$ 		& $\num{4e4}$ 	& -						& Number of mesh cells\\
				$d_{bdry}$					& $2h$ & \si{m} & Interpolation domain, $\Omega_{I}$, offset from $\partial\Omega_{E}$ \\
				$\Delta t_L$				& \num{2.5e-4} & \si{s} & Lagrangian time step size\\
				$\Delta t_E$				& \num{2.5e-5} & \si{s} & Eulerian time step size \\ \hline
		   	\end{tabular}
		\end{table}
		
		\FloatBarrier
		
	\subsection{Dipole in bounded flow} \label{subsections::dipole_in_bounded_flow}
		To further investigate the properties of the hybrid flow solver, namely its ability to capture the generation of vorticity at a solid wall, we applied it to the collision of the Clercx-Bruneau dipole with a solid wall, \cite{Clercx2006}. A FE solution was validated against the results of Clercx and Bruneau, \cite{Clercx2006}, and used as reference.
		
\begin{figure}[!ht]
			\centering
			\includegraphics[width=0.30\textwidth]{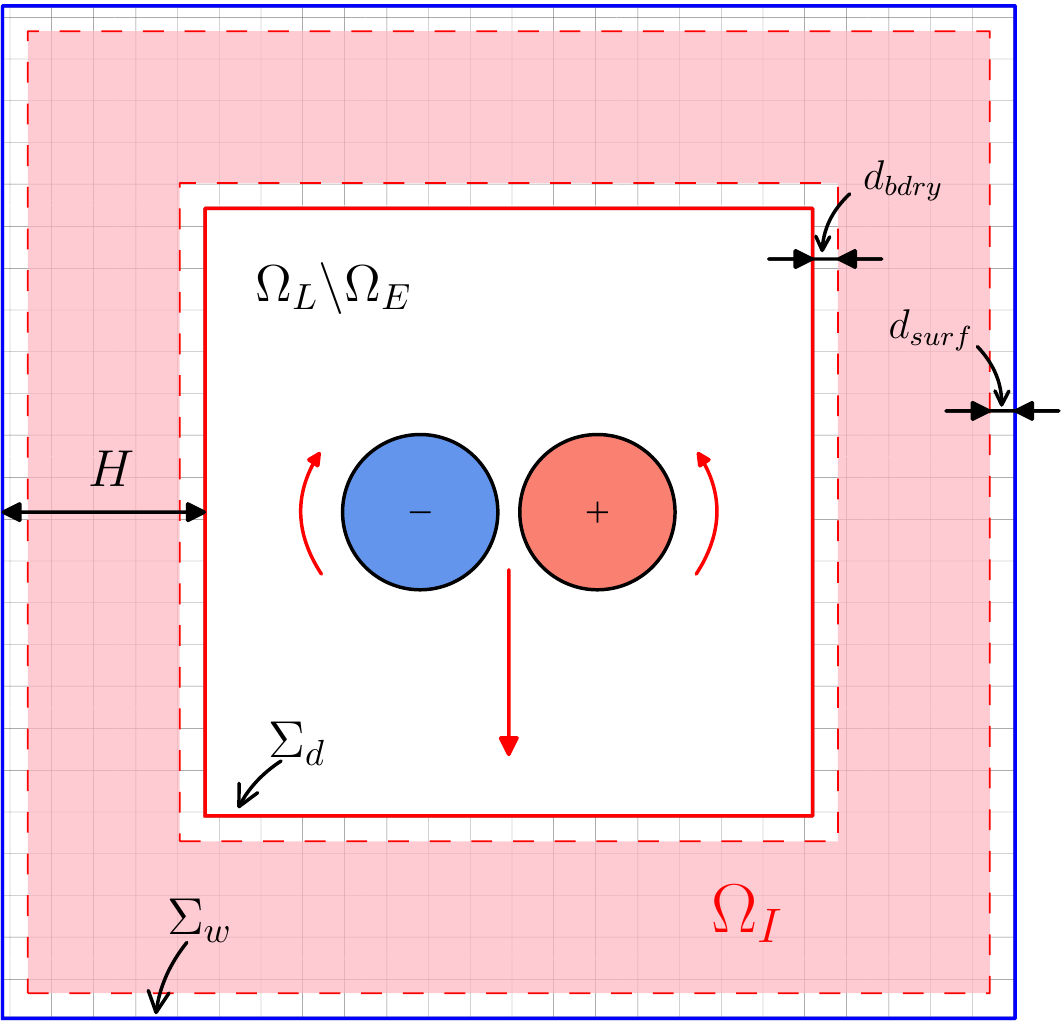}
			\caption{The domain decomposition for the Clercx-Bruneau dipole collision problem, with the positive pole at $p_{+}=(x_1,y_1) = (0.1,0)$ and negative pole at $p_{-}=(x_2,y_2)=(-0.1,0)$. (\textit{Not to scale})}
			\label{fig:hcbdc_dd}
		\end{figure}

		The setup of the hybrid domain is as shown in \figref{fig:hcbdc_dd}. The Eulerian sub-domain, $\Omega_{E}$, covers the near-wall region and the Lagrangian sub-domain domain extends to the complete fluid domain, which is bounded by the no-slip wall $\Sigma_{w}$. The parameters used in this simulation are shown in \tabref{tab:parameters_dipole_bounded}.
		 
		\begin{table}[htbp]
			\centering
		   	\caption{Summary of the parameters used in the hybrid simulation of the Clercx-Bruneau dipole wall collision problem.}
		   	\label{tab:parameters_dipole_bounded}
		   	\begin{tabular}{lcll} 
				\hline\hline
		   		Parameters 					& Value 	& Unit					& Description \\ \hline
				$\nu$						& \num{1.6e-3} & \si{kg.s^{-1}.m^{-1}} & Kinematic viscosity\\
				$\lambda$					& 1 & - & Overlap ratio\\
				$h$							& \num{3e-3} & \si{m} & Nominal blob spacing\\
				$h_{grid}$ 					& \num{5e-3} to \num{1e-2} & \si{m}	& FE cell diameter \\
				$ N_{\mathrm{cells}}$ 		& \num{58272} 	& -						& Number of mesh cells\\
				$d_{bdry}$					& $2h$ & \si{m} & Interpolation domain, $\Omega_{I}$ offset from $\Sigma_{d}$ \\
				$d_{surf}$					& $3h$ & \si{m} &  Interpolation domain, $\Omega_{I}$ offset from $\Sigma_{w}$ \\
				$\Delta t_L$				& \num{2.5e-4} & \si{s} & Lagrangian time step size\\
				$\Delta t_E$				& \num{2.5e-5} & \si{s} & Eulerian time step size \\ \hline
		   	\end{tabular}
		\end{table}

\begin{figure}[!ht]
			\centering
			\includegraphics[width=0.65\textwidth]{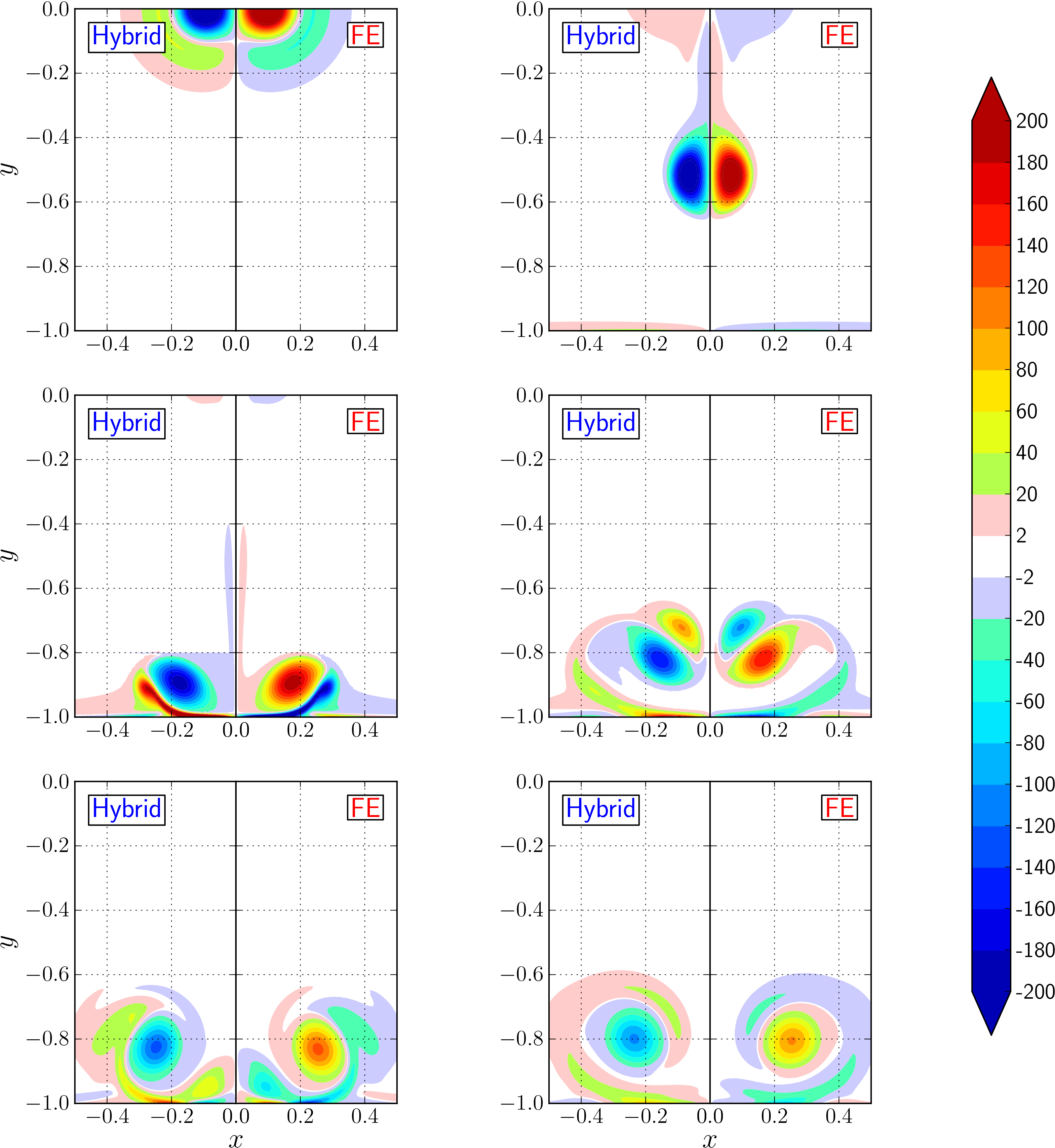}
			\caption{Plot of the dipole at $t = [0, 0.2, 0.4, 0.6, 0.8, 1]$ (from left to right and top to bottom), comparing the hybrid simulation (left half) and FE only simulation (right half).}
			\label{fig:hybrid_doubleMonolope_contourfComparison}
		\end{figure}

\begin{figure}[!ht]
     			\centering
     			\includegraphics[width=0.35\textwidth]{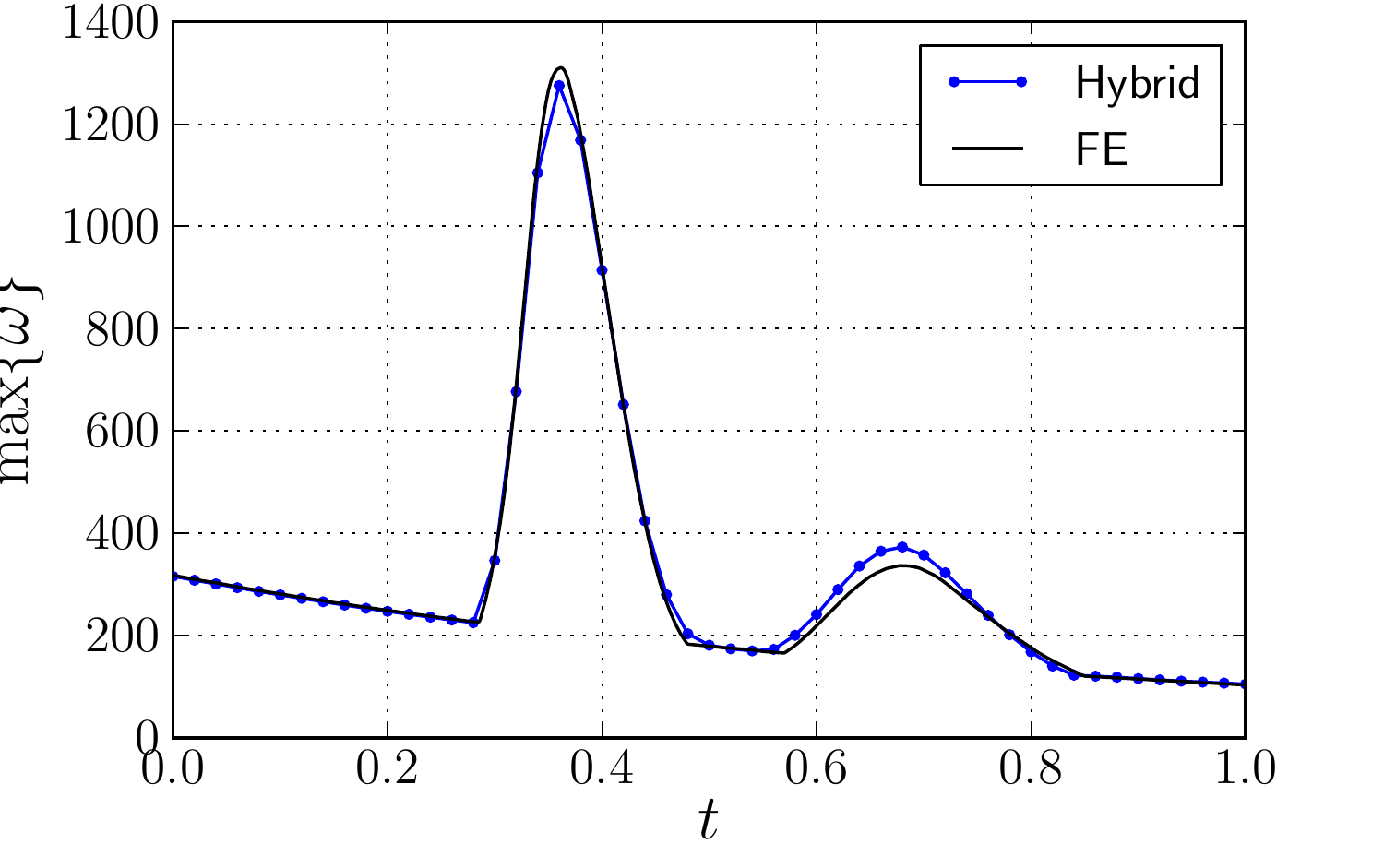}
             		\includegraphics[width=0.35\textwidth]{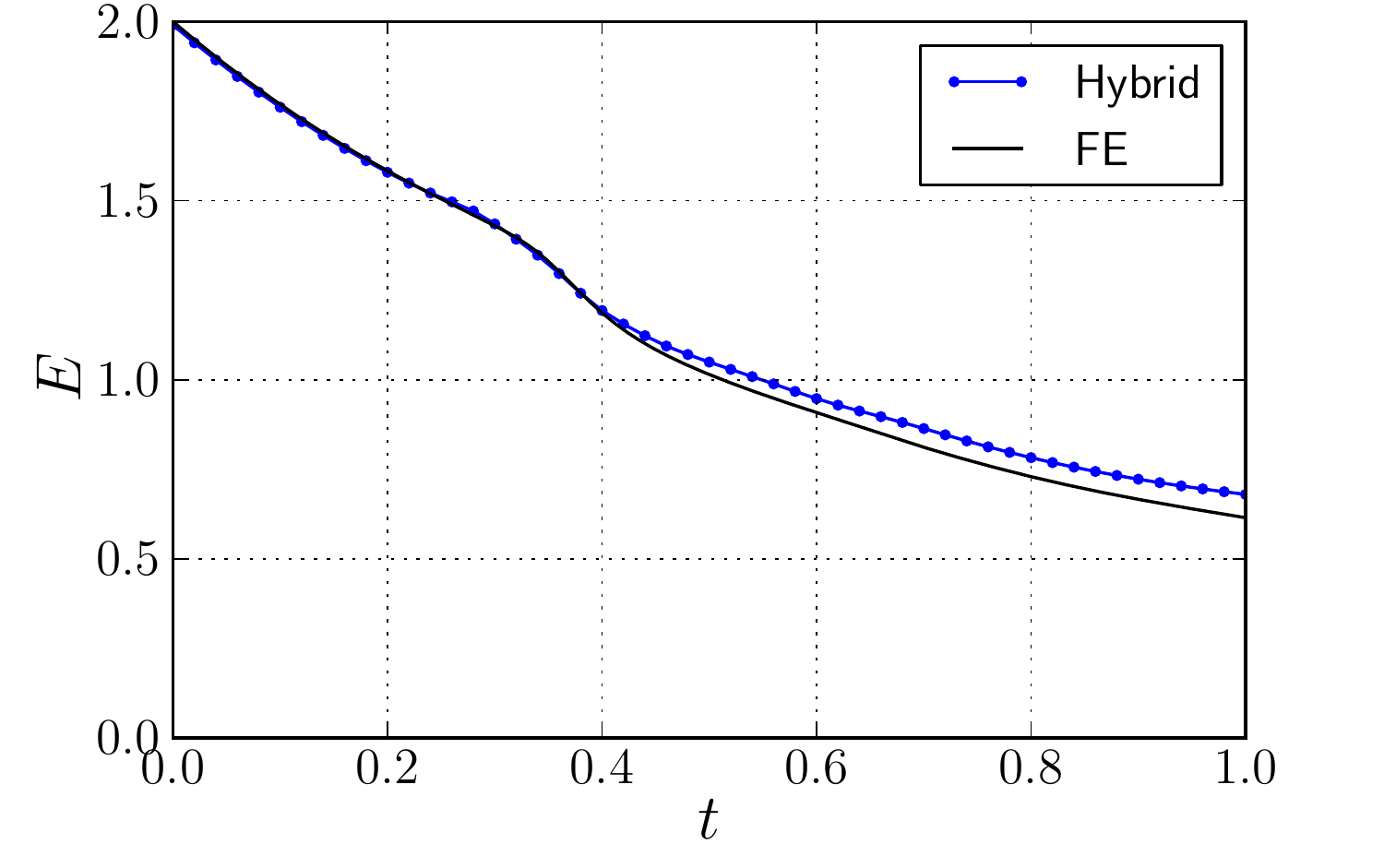}
             		\includegraphics[width=0.35\textwidth]{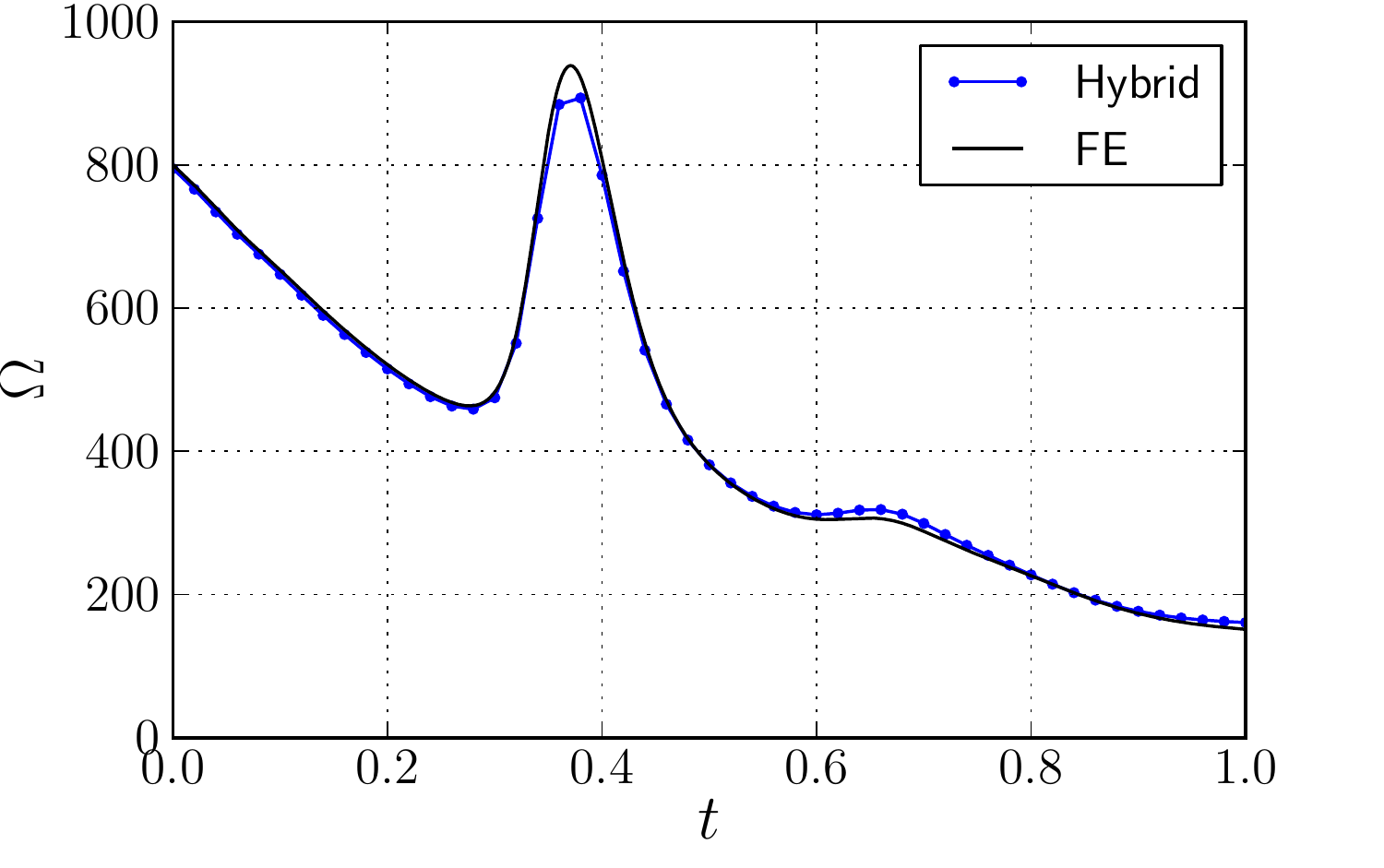}
             		\includegraphics[width=0.35\textwidth]{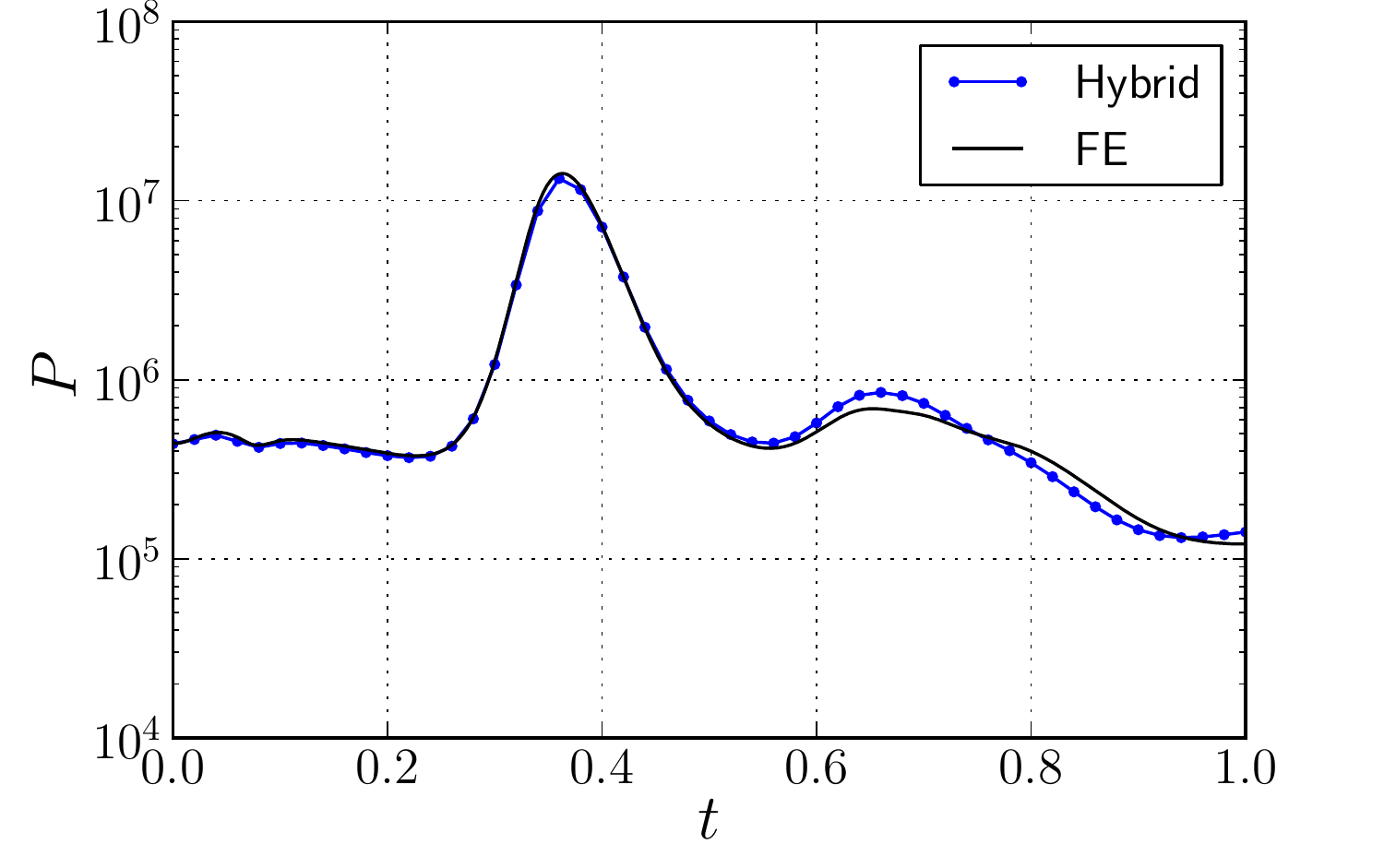}
             		\caption{Variation in the fluid parameters from $t=0$ to $t=1$. The figure compares the hybrid results [{\color{plotBlue}{---}}, solid blue] with the FE only [---, solid black] results.}
     			\label{fig:hybrid_dipole_comparison}
		\end{figure}

		In \figref{fig:hybrid_doubleMonolope_contourfComparison} we show contour plots of vorticity, comparing the hybrid results to the FE ones. We can see a good agreement between the two results. In \figref{fig:hybrid_dipole_comparison} we compare the evolution of vorticity maximum, energy, $E$, enstrophy, $\Omega$, and palinstrophy, $P$. As can be seen, the hybrid solver is capable of reproducing the results of the FE solver, with only small variations on the energy and palinstrophy.
		
		\FloatBarrier
		
	\subsection{Flow around cylinder, $Re=550$} \label{subsection::flow_around_cylinder_re_550}
		An important aspect of flow simulations is the accurate calculation of forces, specifically lift and drag. Therefore, in this section we apply the hybrid solver to the flow around an impulsively started cylinder at $\mathrm{Re}=550$ and determine the forces acting on the cylinder. This test case problem has been extensively analysed in the literature, for example Koumoutaskos and Leonard, \cite{Koumoutsakos1995}, and Rosenfeld et al., \cite{Rosenfeld1991}, and these results will serve as a reference for the assessment of the hybrid solver, since one is a pure vortex particle solver and the other one is a pure Eulerian grid solver.
		
\begin{figure}[!ht]
			\centering
			\includegraphics[width=0.40\textwidth]{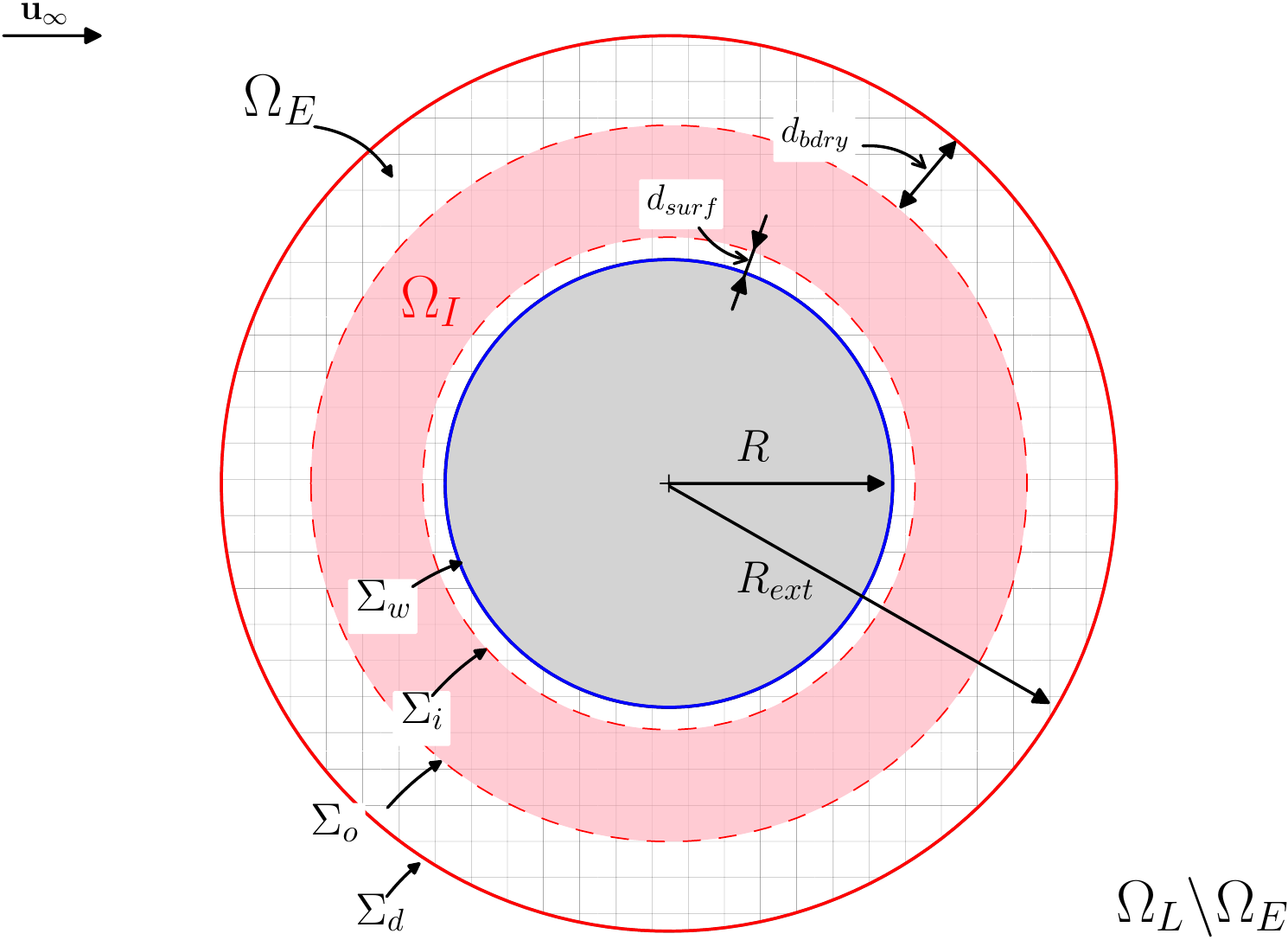}
			\caption{The domain decomposition for the impulsively started cylinder. The parameters of the domain are shown in Table \ref{tab:parameters_cylinder}. (\textit{Not to scale})}
			\label{fig:hisc_dd}
		\end{figure}

		The configuration of the hybrid domain is presented in \figref{fig:hisc_dd}. As before, the Lagrangian domain, $\Omega_{L}$, covers the whole fluid domain and the Eulerian domain, $\Omega_{E}$, is confined to a small region around the cylinder. The parameters used in this simulation are presented in \tabref{tab:parameters_cylinder}.
		
		\begin{table}[htbp]
			\centering
		   	\caption{Summary of the parameters used in the hybrid simulation of the impulsively started cylinder at $Re=550$.}
		   	\label{tab:parameters_cylinder}
		   	\begin{tabular}{lcll} 
				\hline\hline
		   		Parameters 					& Value 	& Unit					& Description \\ \hline
				$\boldsymbol{U}_{\infty}$						& $1.0\boldsymbol{e}_{x}$ & \si{m.s^{-1}} & Free stream velocity\\
				$\nu$						& \num{3.6e-3} & \si{kg.s^{-1}.m^{-1}} & Kinematic viscosity\\
				$\lambda$					& 1 & - & Overlap ratio\\
				$h$							& \num{8e-3} & \si{m} & Nominal blob spacing\\
				$h_{grid}$ 					& \num{8e-3} to \num{4e-2} & \si{m}	& FE cell diameter \\
				$ N_{\mathrm{cells}}$ 		& \num{32138} 	& -						& Number of mesh cells\\
				$R$					& \num{1.0} & \si{m} & Radius of the cylinder \\
				$R_{ext}$					& \num{1.5} & \si{m} & Radius of the external boundary $\Sigma_{d}$ \\
				$d_{bdry}$					& $0.2R$ & \si{m} & Interpolation domain, $\Omega_{I}$ offset from $\Sigma_{d}$ \\
				$d_{surf}$					& $3h$ & \si{m} &  Interpolation domain, $\Omega_{I}$ offset from $\Sigma_{w}$ \\
				$\Delta t_L$				& \num{3.0e-3} & \si{s} & Lagrangian time step size\\
				$\Delta t_E$				& \num{1.0e-3} & \si{s} & Eulerian time step size \\ \hline
		   	\end{tabular}
		\end{table}
				
\begin{figure}[!ht]
			\centering
			\includegraphics[width=0.75\textwidth]{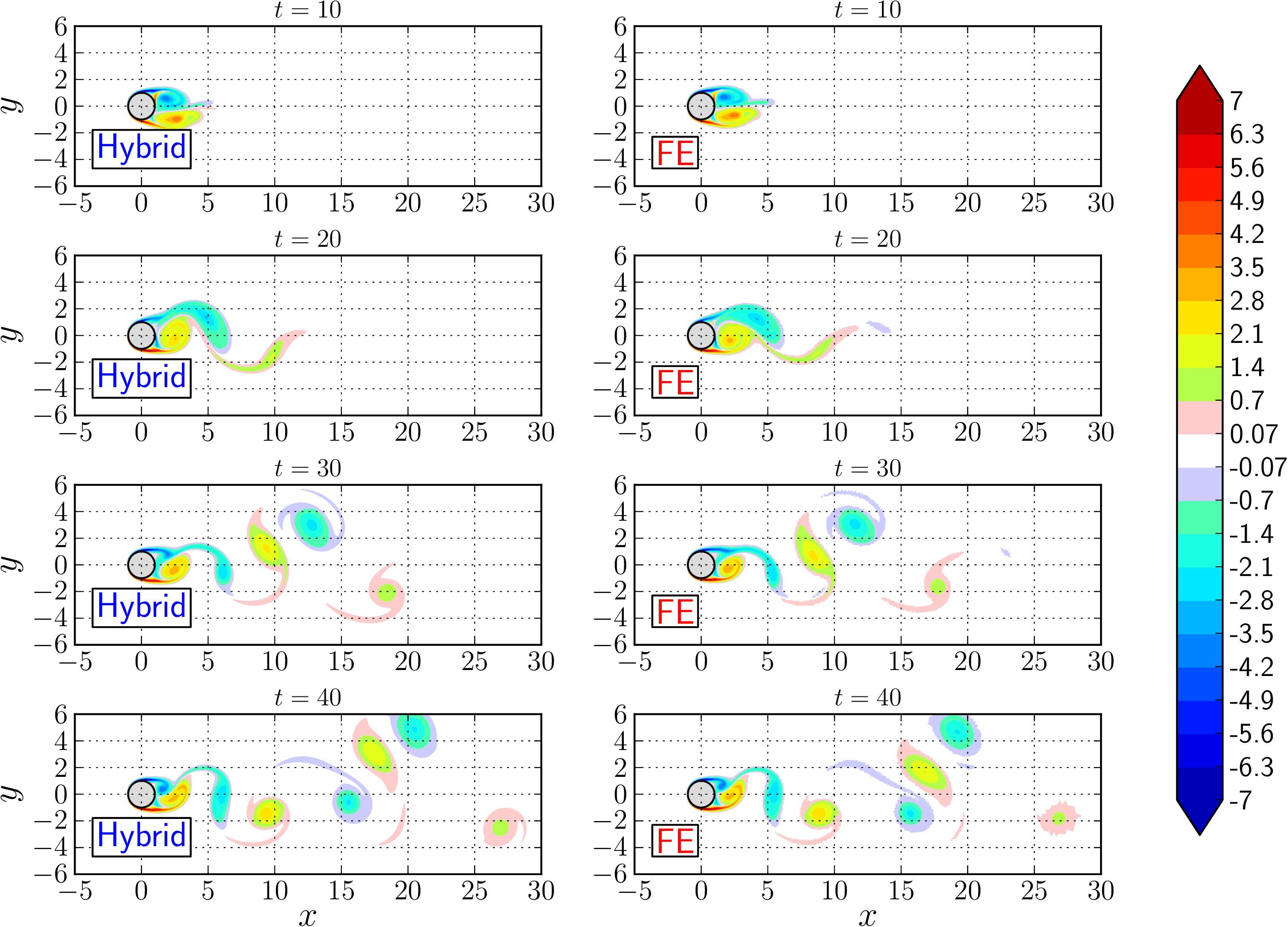}
			\caption{Plot of the vorticity field at $t=[10,20,30,40]$, comparing the hybrid simulation (left) with the FE simulation (right).}
			\label{fig:hybrid_cylinder_LongRun_contourfComparison}
		\end{figure}

\begin{figure}[!ht]
			\centering
			\includegraphics[width=0.35\textwidth]{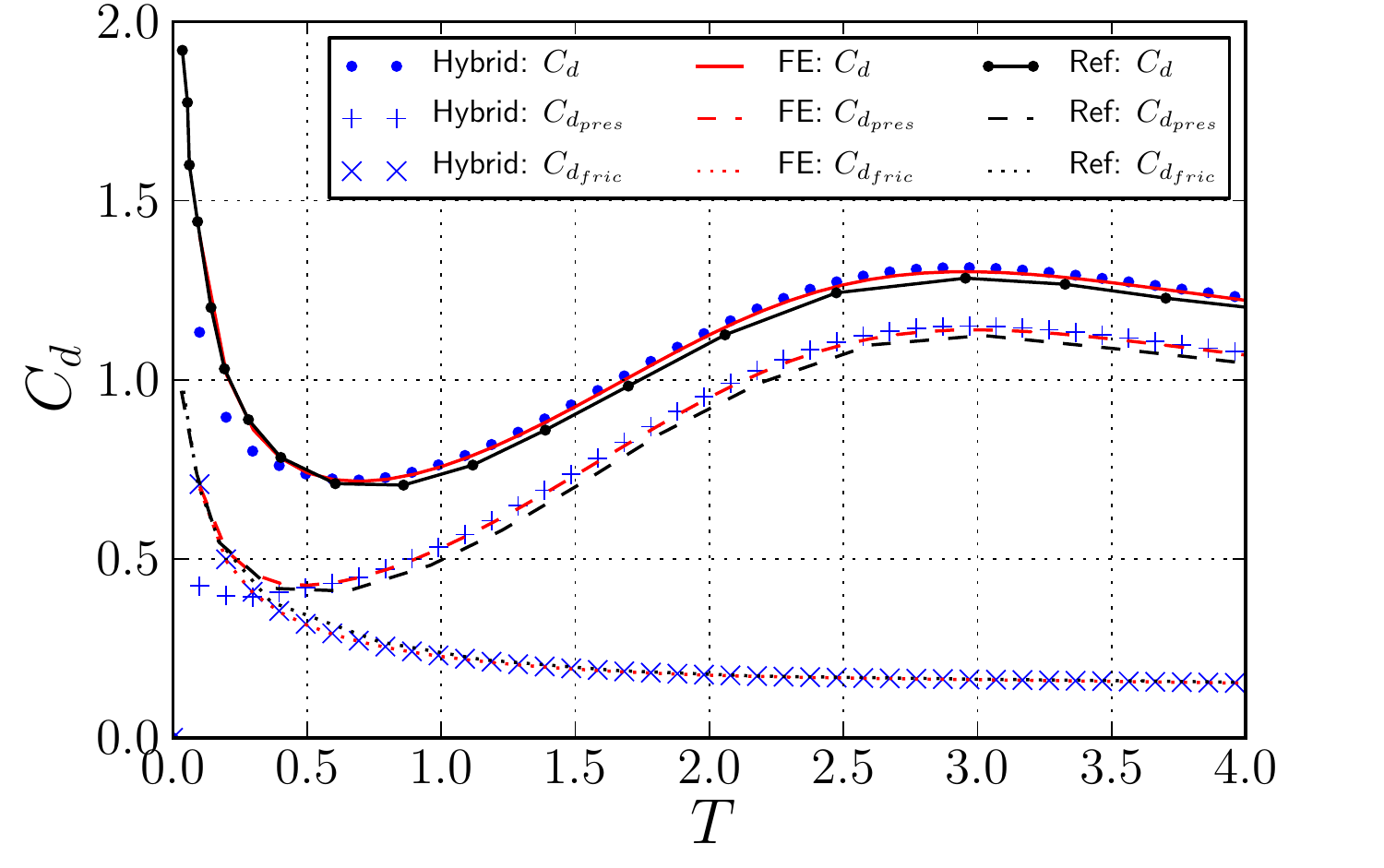}
			\includegraphics[width=0.35\textwidth]{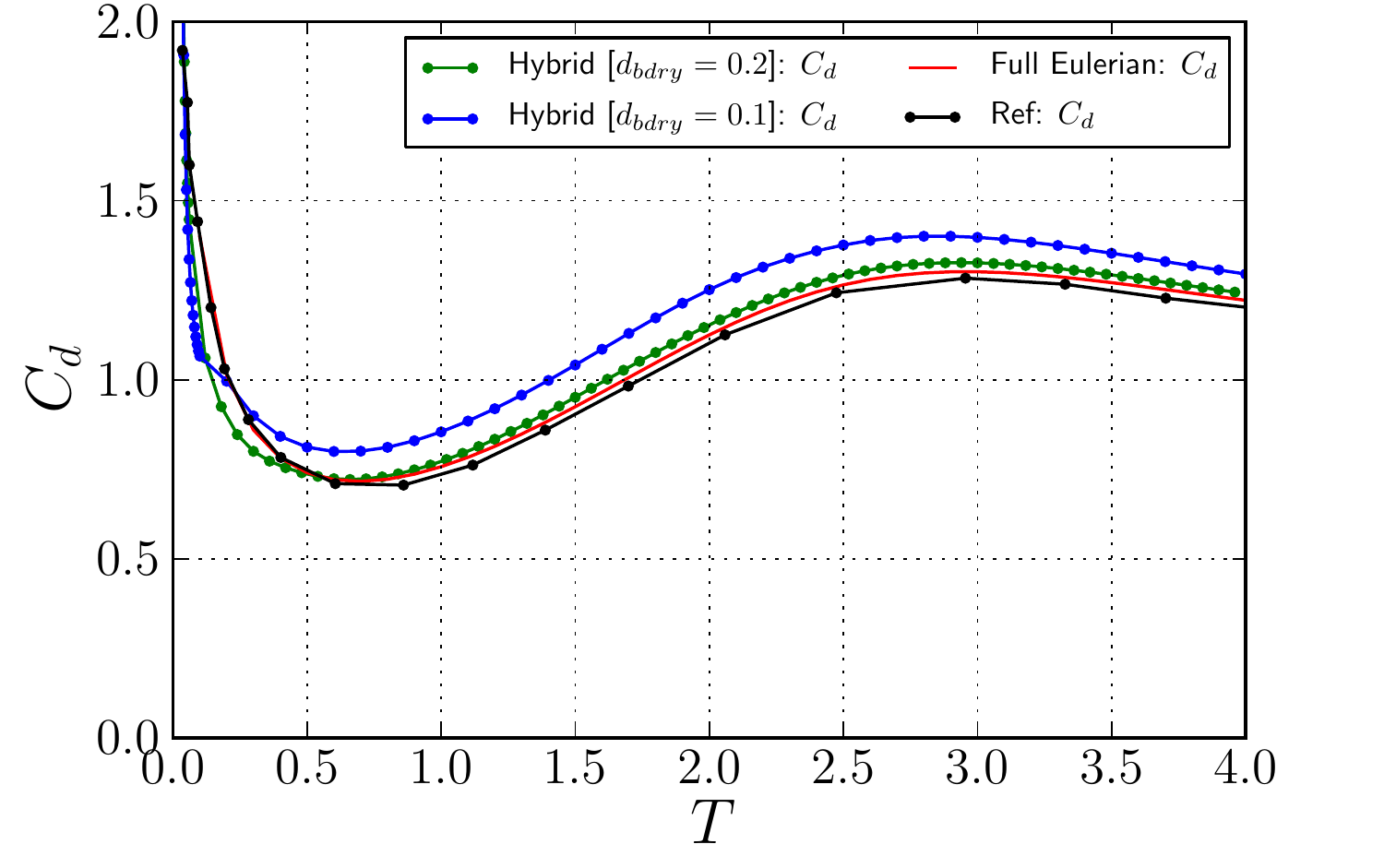}
			\caption{Left: Evolution of the drag coefficient during the initial stages $t=0$ to $t=4$ with total drag coefficient $C_d$ (solid), pressure drag coefficient $C_{d_{pres}}$ (dashed) and friction drag coefficient $C_{d_{fric}}$ (dotted). The figure compares results of hybrid simulation ({\color{plotBlue}{\textbf{blue}}}), FE simulation ({\color{plotRed}{\textbf{red}}}) and reference data (\textbf{black}) of Koumoutsakos and Leonard \cite{Koumoutsakos1995}. Right: Comparison of total drag coefficient $C_{d}$ for  $d_{bdry} \in \{0.1R,0.2R\}$. }
			\label{fig:hybrid_ISC_drag}
		\end{figure}

\begin{figure}[!ht]
			\centering
			\includegraphics[width=0.35\textwidth]{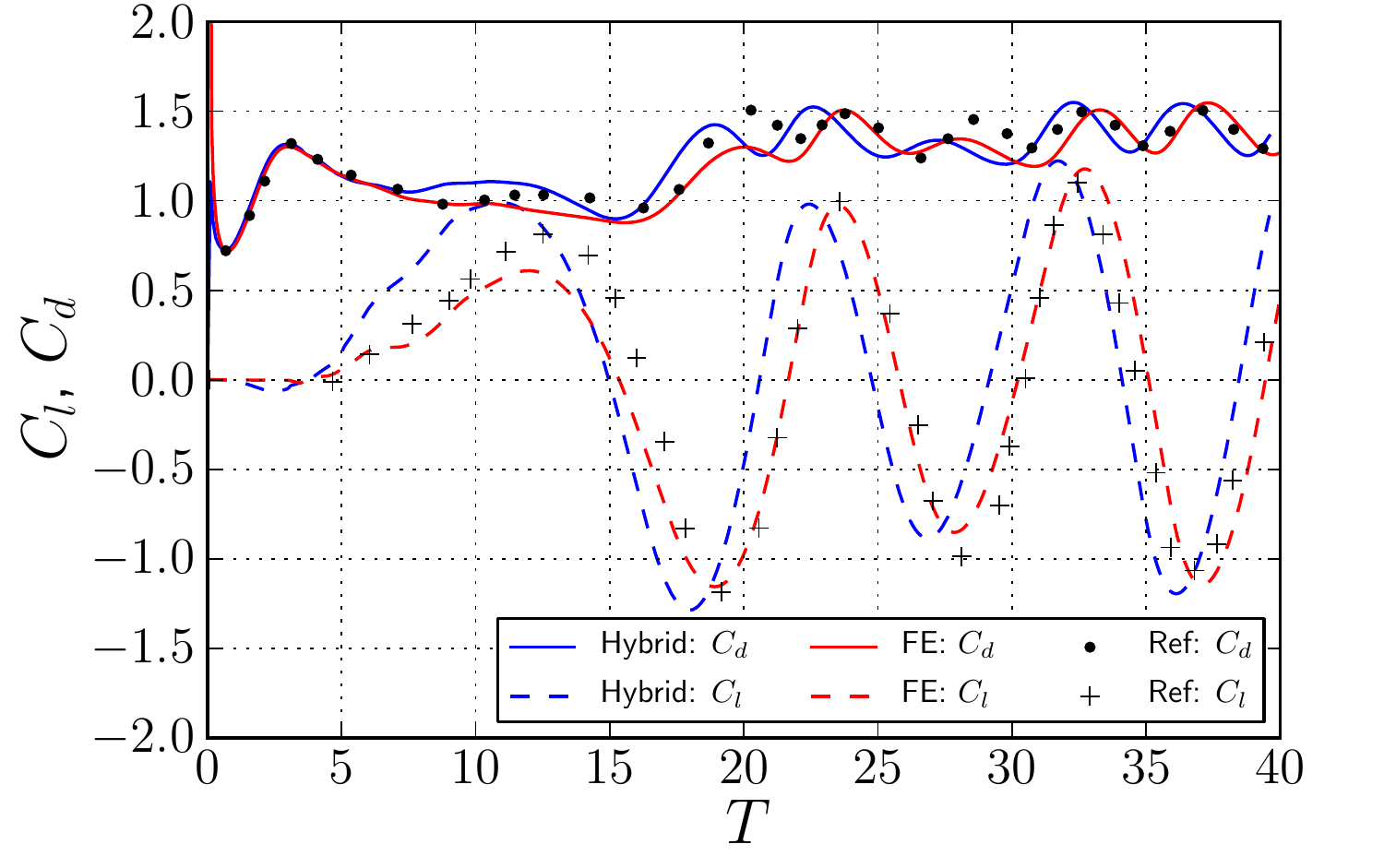}
			\caption{Evolution of the lift and drag coefficients from $t=0$ to $t=40$ with artificial perturbation \cite{Lecointe1984}. The figure compares hybrid ({\color{plotBlue}{\textbf{blue}}}), FE only ({\color{plotRed}{\textbf{red}}}), and the reference data (\textbf{black}) from Rosenfeld et al. \cite{Rosenfeld1991}.}
			\label{fig:hybrid_cylinder_LongRun_liftDrag}
		\end{figure}

		The contour plots of vorticity, comparing the hybrid results to the FE ones, are presented in \figref{fig:hybrid_cylinder_LongRun_contourfComparison}. We can see that the two solvers give very similar results. Regarding the drag and lift, we can say that the hybrid solver is capable of reproducing both the FE results and the results of Koumoutsakos and Leonard, \cite{Koumoutsakos1995}, see \figref{fig:hybrid_ISC_drag}, left. The differences exist mainly in the first instants of the simulation. It is important to highlight that, increasing $h_{bdry}$ from $\num{0.1}R$ to $\num{0.2}R$ improved considerably the results, see  \figref{fig:hybrid_ISC_drag}, right. We have noted that this parameter is important and further research should be done in order to find its optimal value. A longer time simulation, $t\in[0,40]$, was performed and the lift and drag compared to the results of Rosenfeld et al. \cite{Rosenfeld1991}, see \figref{fig:hybrid_cylinder_LongRun_liftDrag}. The hybrid solver follows very well both the reference results of Rosenfeld et al. and of the FE simulation up to the onset of the vortex shedding, $t\approx 5s$. After that, all results stop having a good match but all remain within the same bounds and show similar frequency, as expected.
		
\begin{figure}[!ht]
     			\centering
     			\includegraphics[width=0.395\textwidth]{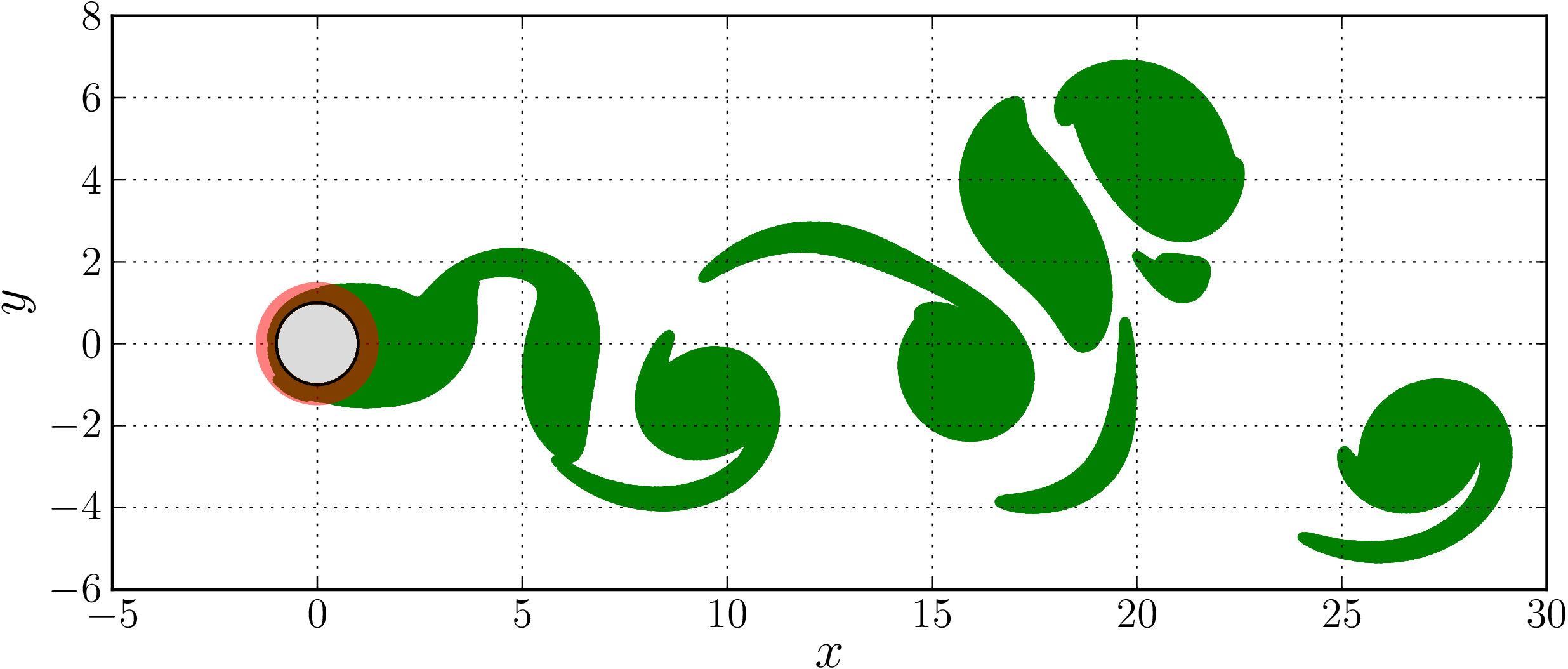}
             		\includegraphics[width=0.45\textwidth]{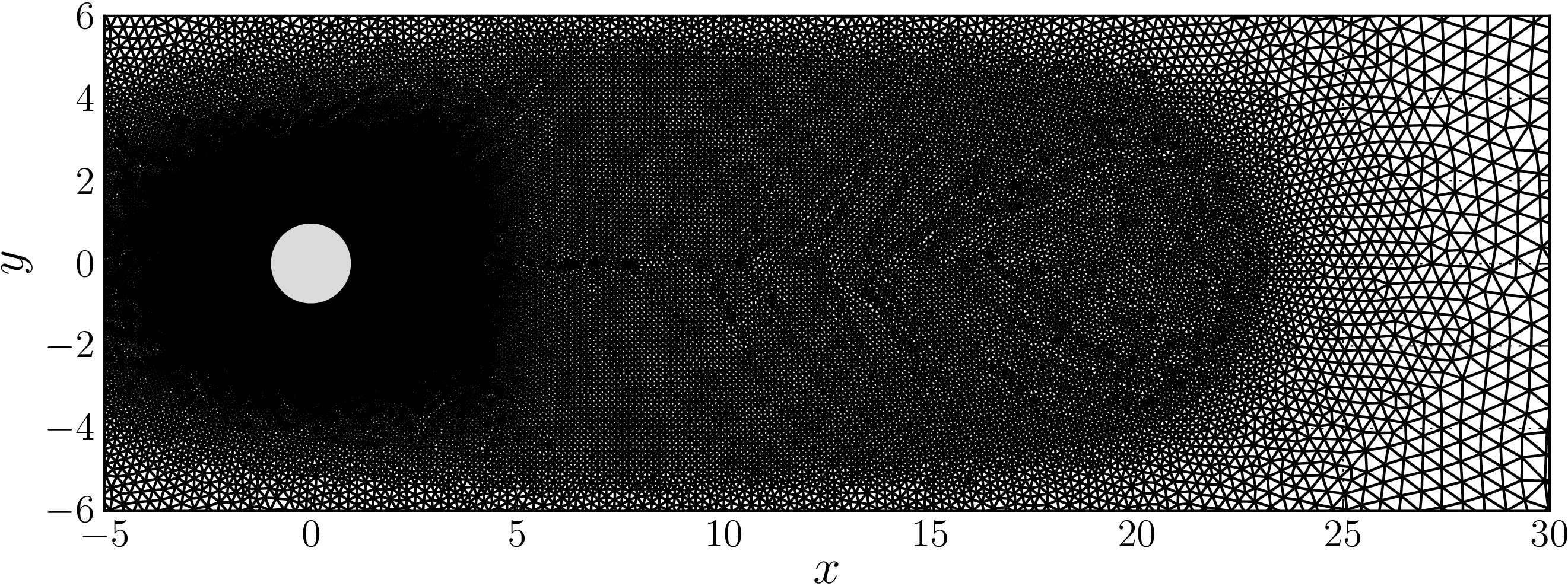}
             		\caption{Spatial distribution of computational elements for the hybrid method (left) and the FE method (right).}
     			\label{fig:hybrid_gridResolution}
		\end{figure}

		Another aspect of the hybrid method, which is inherited from the vortex particle method, is its automatic adaptivity, where computational elements exist only in the support of vorticity, as opposed to standard grid solvers where the computational elements  exist in the whole computational domain, see \figref{fig:hybrid_gridResolution}.
		
		\FloatBarrier

	\subsection{Flow around stalled ellipsoid, $Re=5000$} \label{subsection::flow_around_stalled_ellipsoid_re_5000}
		The final test case consists in the simulation of the flow around a stalled ellipsoid at a Reynolds number $Re=5000$. With this test case we aimed to assess the hybrid method's capability to simulate flows with higher Reynolds number and to evaluate its performance in a situation where vortices exit and re-enter the Eulerian sub-domain.
		
\begin{figure}[!ht]
			\centering
			\includegraphics[width=0.45\textwidth]{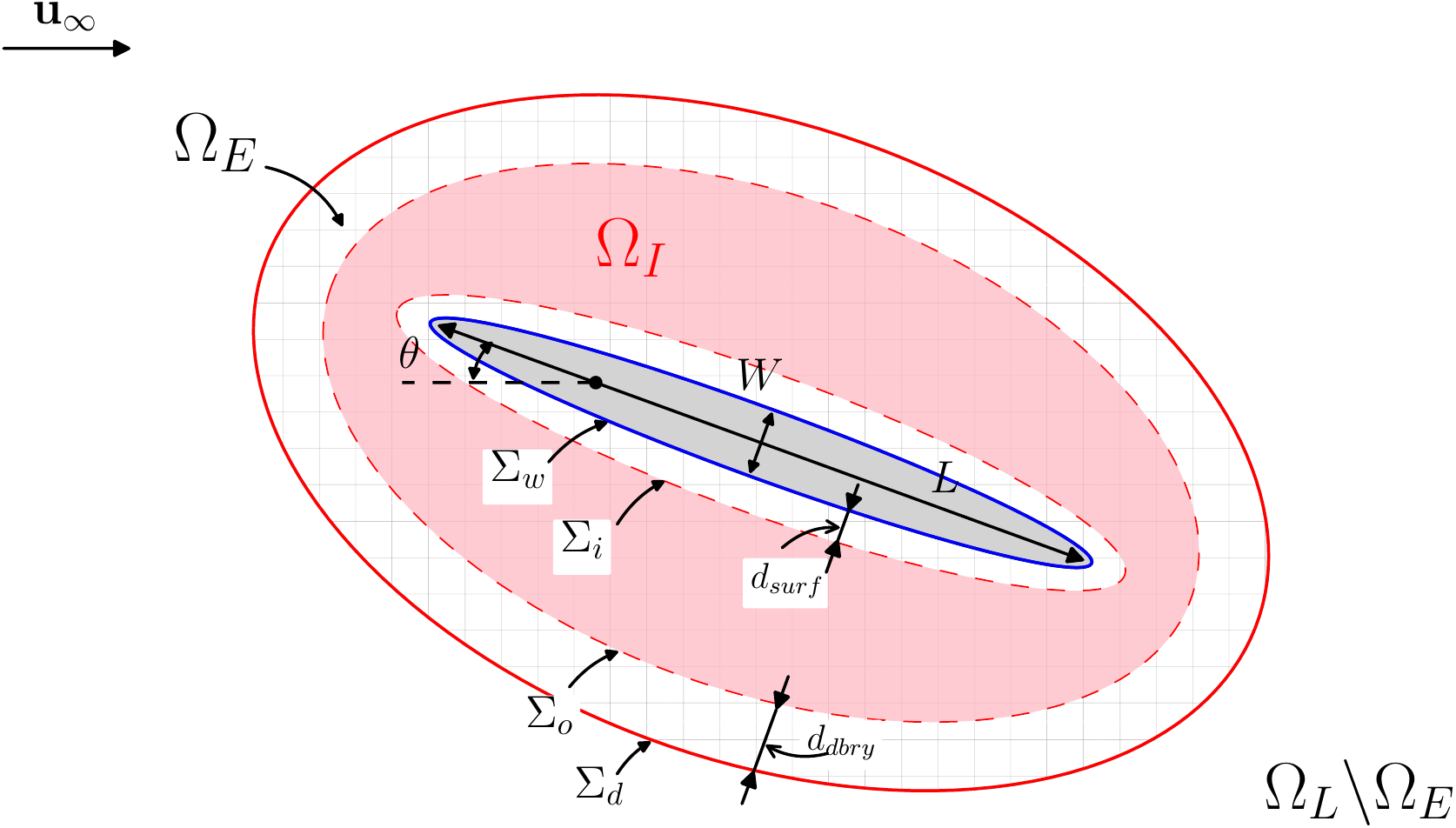}
			\caption{The domain decomposition for the stalled elliptical airfoil test case. The parameters of the domain are tabulated in \tabref{tab:parameters_ellipsoid}. (\textit{Not to scale})}
			\label{fig:hellipticAirfoil_dd-crop}
		\end{figure}

		The configuration of the hybrid domain is presented in \figref{fig:hellipticAirfoil_dd-crop}. Once more, the Lagrangian sub-domain, $\Omega_{L}$, covers the whole fluid domain and the Eulerian domain, $\Omega_{E}$, is restricted to the vicinity of the solid boundary. The parameters used in this simulation are presented in \tabref{tab:parameters_ellipsoid}.
		
		\begin{table}[htbp]
			\centering
		   	\caption{Summary of the parameters used in the hybrid simulation of the stalled ellipsoid at $Re=5000$.}
		   	\label{tab:parameters_ellipsoid}
		   	\begin{tabular}{lcll} 
				\hline\hline
		   		Parameters 					& Value 	& Unit					& Description \\ \hline
				$\boldsymbol{U}_{\infty}$						& $1.0\boldsymbol{e}_{x}$ & \si{m.s^{-1}} & Free stream velocity\\
				$\nu$						& \num{2e-4} & \si{kg.s^{-1}.m^{-1}} & Kinematic viscosity\\
				$\lambda$					& 1 & - & Overlap ratio\\
				$h$							& \num{1.67e-3} & \si{m} & Nominal blob spacing\\
				$h_{grid}$ 					& \num{1.4e-3} to \num{8e-3} & \si{m}	& FE cell diameter \\
				$ N_{\mathrm{cells}}$ 		& \num{118196} 	& -						& Number of mesh cells\\
				$L$					& \num{1.0} & \si{m} & Chord length \\
				$L_{ext}$					& \num{1.0} & \si{m} & Length of the external boundary $\Sigma_{d}$\\
				$W$					& \num{0.12} & \si{m} & Maximum thickness \\
				$W_{ext}$					& \num{1.0} & \si{m} & Maximum thickness of the external boundary $\Sigma_{d}$\\
				$R_{ext}$					& \num{1.5} & \si{m} & Radius of the external boundary $\Sigma_{d}$ \\
				$d_{bdry}$					& $0.1L$ & \si{m} & Interpolation domain, $\Omega_{I}$ offset from $\Sigma_{d}$ \\
				$d_{surf}$					& $3h$ & \si{m} &  Interpolation domain, $\Omega_{I}$ offset from $\Sigma_{w}$ \\
				$\Delta t_L$				& \num{3.0e-3} & \si{s} & Lagrangian time step size\\
				$\Delta t_E$				& \num{1.0e-3} & \si{s} & Eulerian time step size \\ \hline
		   	\end{tabular}
		\end{table}
		
\begin{figure}[!ht]
			\centering
			\includegraphics[width=0.65\textwidth]{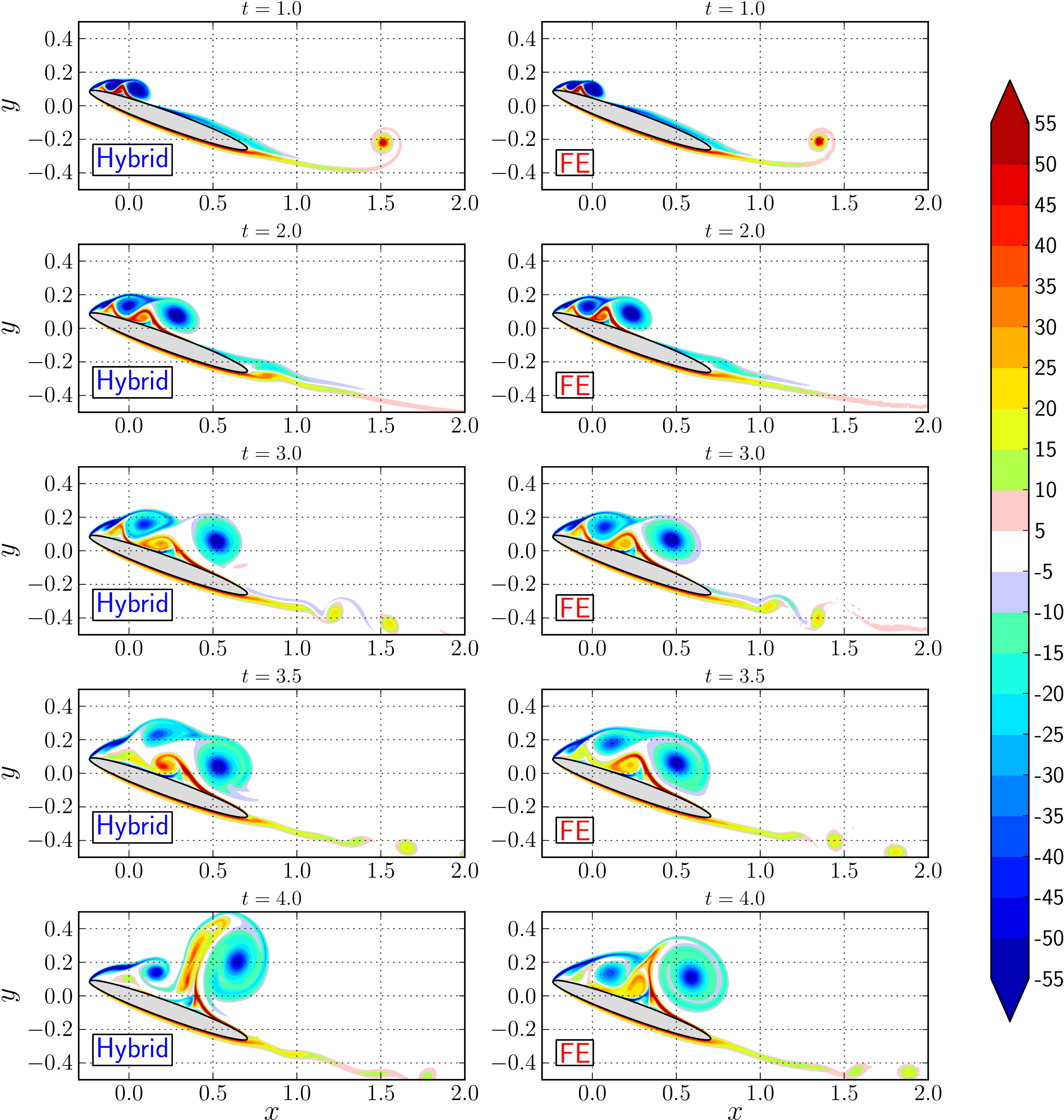}
			\caption{Plot of the vorticity field for $t \in \{1.0, 2.0, 3.0, 3.5, 4.0\}$, comparing the hybrid simulation (left) with the FE simulation (right).}
			\label{fig:hybrid_ellipse_Hybrid_contours}
		\end{figure}

\begin{figure}[!ht]
     			\centering
     			\includegraphics[width=0.35\textwidth]{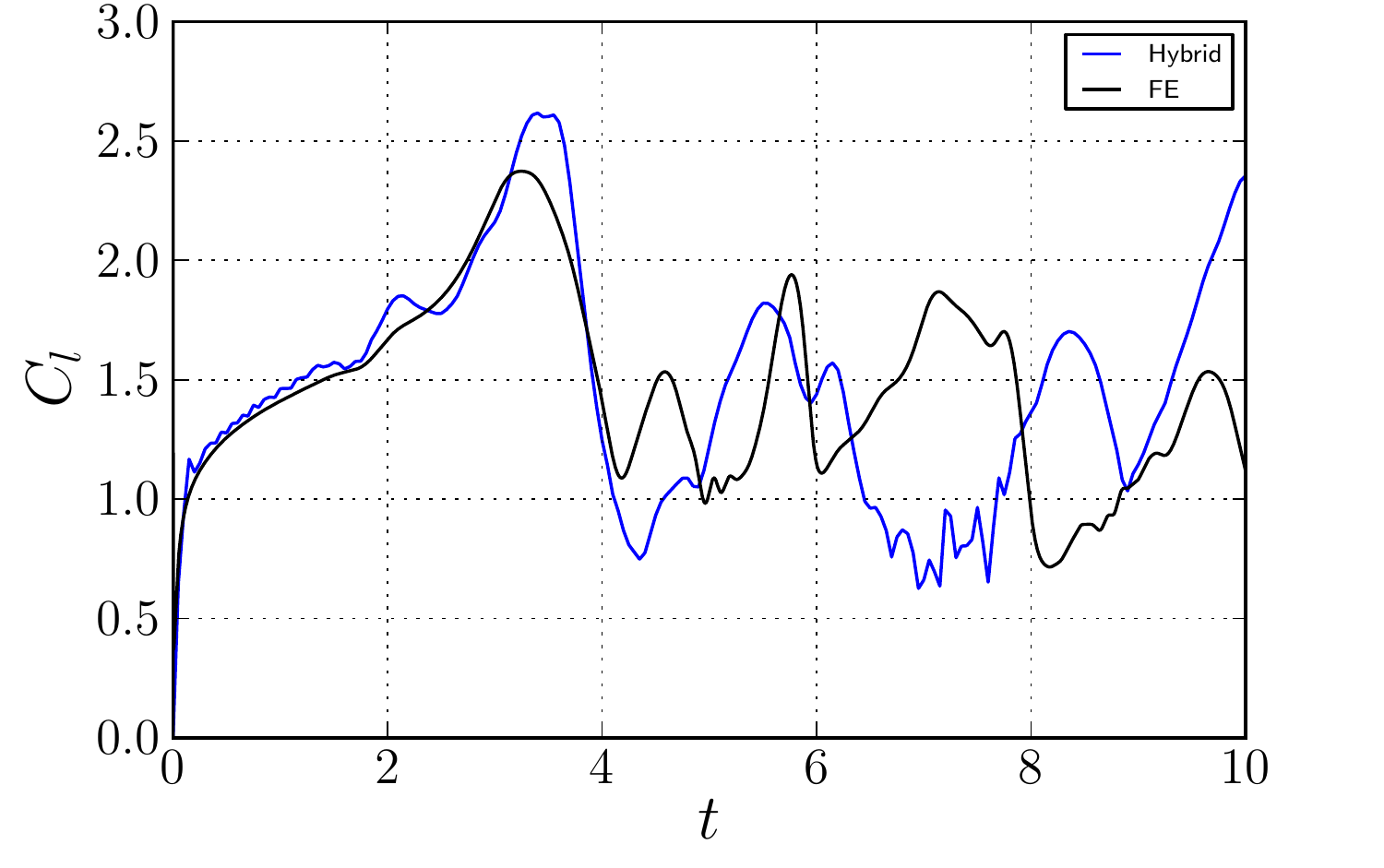}
             		\includegraphics[width=0.35\textwidth]{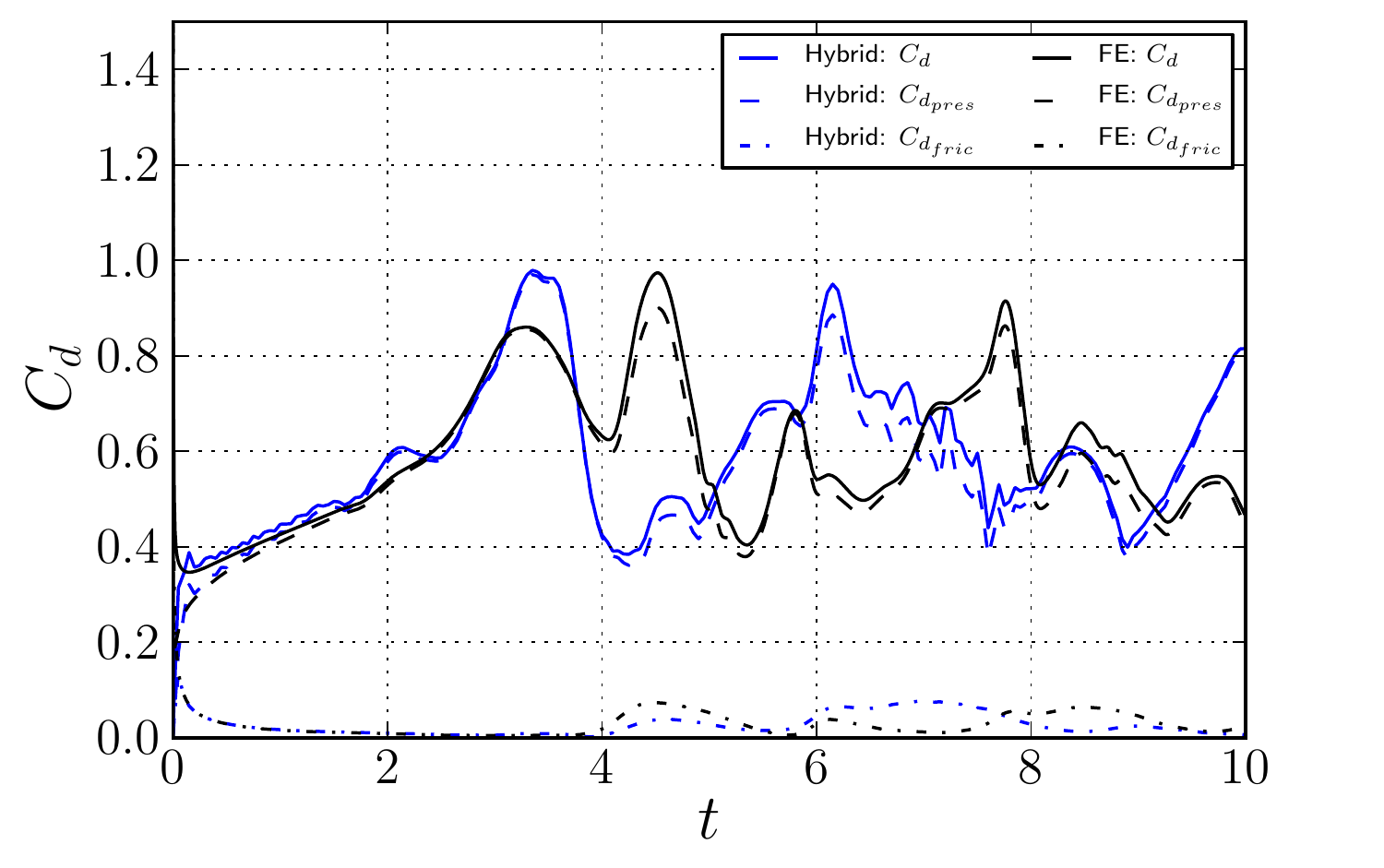}
             		\caption{Evolution of the lift and drag coefficient from $t=0$ up to $t=10$. The figure compares the hybrid results with FE simulation.}
     			\label{fig:hybrid_ellipseForce}
		\end{figure}

		The contour plots of vorticity, comparing the hybrid results to the FE ones, are presented in \figref{fig:hybrid_ellipse_Hybrid_contours}. We can see that the two solvers give very similar results up to $t=3s$ and after that they start to diverge. This behaviour is expected due to the non-linear nature of the problem. Nevertheless we can see a comparable vortical structure. Regarding the drag and lift, we can observe that the hybrid solver is capable of reproducing very well the FE results up to $t=2s$, see \figref{fig:hybrid_ellipseForce}, left. After that instant, which corresponds to the start of the vortex shedding, the results stop having a good match, but remain within the same bounds and show similar behaviour.
		
		\FloatBarrier

%
%
%
		
\section{Conclusions and outlook} \label{section::conclusions_and_outlook}
	In this work we have presented an efficient hybrid Eulerian-Lagrangian flow solver that does not rely on the Schwarz iteration method and is capable of exactly preserving the total circulation, improving the results of similar approaches.
	
	By dividing the computational domain into two types of sub-domains, Lagrangian based and Eulerian based, it is possible to use the most suitable solver in each region of the flow domain. This formulation allows for the use of small and highly resolved Eulerian solvers in the near wall regions. These solvers can efficiently capture the generation of vorticity at the solid boundaries, potentially using wall functions and other advanced methods. In the wake region, the Lagrangian solver accurately models the development of the wake due to its practically absent numerical diffusion.
	
	This approach paves the way to the development of complex simulations of wind turbines since each solid object can be individually modelled by a disjoint Eulerian subdomain embedded in the Lagrangian one. This greatly simplifies the parallelizability  of large flow simulations.
	
	The future developments of this work are the extension to moving and deformable objects, extension to large Reynolds number flows (including the coupling of turbulence between the two sub-domains) and extension to three dimensional flows.
	
\section*{Acknowledgements}
	This work is made possible in part due to funds from Technology Foundation STW - Veni grant, and
Future Emerging Technologies FP7-Deepwind project and the U.S. Department of Energy.



\bibliographystyle{elsart-num-sort} 
\def\url#1{}
\bibliography{library}

\begin{thebibliography}{10}
\expandafter\ifx\csname url\endcsname\relax
  \def\url#1{\texttt{#1}}\fi
\expandafter\ifx\csname urlprefix\endcsname\relax\def\urlprefix{URL }\fi

\bibitem{BarbaPhD}
L.~A. Barba, {Vortex method for computing high-Reynolds number flows: Increased
  accuracy with a fully mesh-less formulation}, Ph.D. thesis, California
  Institute of Technology (2004).

\bibitem{Barba2005}
L.~A. Barba, A.~Leonard, C.~B. Allen, {Advances in viscous vortex
  methods--meshless spatial adaption based on radial basis function
  interpolation}, International Journal for Numerical Methods in Fluids 47~(5)
  (2005) 387--421.

\bibitem{BarbaFMMexa2013}
L.~A. Barba, R.~Yokota, {How Will the Fast Multipole Method Fare in the
  Exascale Era?}, SIAM News 46~(6).

\bibitem{Barnes1986}
J.~Barnes, P.~Hut, A hierarchical $\mathcal{O}(n \log n)$ force-calculation
  algorithm, Nature 324~(6096) (1986) 446--449.

\bibitem{Beale1988}
J.~Beale, {On the Accuracy of Vortex Methods at Large Times}, in: Computational
  Fluid Dynamics and Reacting Gas Flows SE - 2, vol.~12, Springer New York,
  1988, pp. 19--32.

\bibitem{Beale1985}
J.~Beale, A.~Majda, {High order accurate vortex methods with explicit velocity
  kernels}, Journal of Computational Physics 58~(2) (1985) 188--208.

\bibitem{Beale1985a}
J.~Beale, A.~Majda, {High order accurate vortex methods with explicit velocity
  kernels}, Journal of Computational Physics 58~(2) (1985) 188--208.

\bibitem{Beale1982}
J.~T. Beale, A.~Majda, {Vortex methods. II. Higher order accuracy in two and
  three dimensions}, Mathematics of Computation 39~(159) (1982) 29--52.

\bibitem{brezzi1991mixed}
F.~Brezzi, M.~Fortin, {Mixed and hybrid finite element methods}, vol.~15,
  Springer-Verlag, New York, 1991.

\bibitem{Chatelin2014}
R.~Chatelin, P.~Poncet, {Hybrid grid-particle methods and Penalization: A
  Sherman-Morrison-Woodbury approach to compute 3D viscous flows using FFT},
  Journal of Computational Physics 269 (2014) 314--328.

\bibitem{Chorin1968}
A.~J. Chorin, {Numerical solution of the Navier-Stokes equations}, Mathematics
  of Computation 22~(104) (1968) 745--762.

\bibitem{Chorin1973}
A.~J. Chorin, {Numerical study of slightly viscous flow}, Journal of Fluid
  Mechanics 57~(04) (1973) 785--796.

\bibitem{Clercx2006}
H.~Clercx, C.-H. Bruneau, {The normal and oblique collision of a dipole with a
  no-slip boundary}, Computers \& Fluids 35~(3) (2006) 245--279.

\bibitem{Cottet1982}
G.~H. Cottet, {M\'{e}thodes particulaires pour l'\'{e}quation d'Euler dans le
  plan}, Ph.D. thesis, Universit\'{e} Paris VI, Paris (1982).

\bibitem{Cottet1991book}
G.~H. Cottet, {Particle-grid domain decomposition methods for the Navier-Stokes
  equations in exterior domains}, in: Vortex Dynamics and Vortex Methods,
  American Mathematical Society, 1991, pp. 103--117.

\bibitem{Cottet1994}
G.-H. Cottet, {A vorticity creation algorithm for the Navier-Stokes equations
  in arbitrary domain}, in: Navier-Stokes equations and related non-linear
  problems, 1994.

\bibitem{CottetKoumoutsakos2000}
G.-H. Cottet, P.~D. Koumoutsakos, {Vortex Methods: Theory and Practice},
  Cambridge University Press, 2000.

\bibitem{Cottet1999}
G.-H. Cottet, M.-L. {Ould Salihi}, M.~{El Hamroui}, {Multi-purpose regridding
  in vortex methods}, ESAIM: Proceedings 7 (1999) 94--103.

\bibitem{daeninckThesis}
G.~Daeninck, {Developments in hybrid approaches : Vortex method with known
  separation location}, Ph.D. thesis, Universit\'{e} Catolique de Louvain
  (2006).

\bibitem{Degond1989}
P.~Degond, S.~Mas-Gallic, {The Weighted Particle Method for
  Convection-Diffusion Equations Part 1: The Case of an Isotropic Viscosity},
  Mathematics of Computation 53 (1989) 485--507.

\bibitem{Dong2014}
S.~Dong, G.~E. Karniadakis, C.~Chryssostomidis, {A robust and accurate outflow
  boundary condition for incompressible flow simulations on severely-truncated
  unbounded domains}, Journal of Computational Physics 261 (2014) 83--105.

\bibitem{Engblom2011}
S.~Engblom, {On well-separated sets and fast multipole methods}, Applied
  Numerical Mathematics 61~(10) (2011) 1096--1102.

\bibitem{Goda1979}
K.~Goda, {A multistep technique with implicit difference schemes for
  calculating two- or three-dimensional cavity flows}, Journal of Computational
  Physics 30~(1) (1979) 76--95.

\bibitem{Golas2012}
A.~Golas, R.~Narain, J.~Sewall, P.~Krajcevski, P.~Dubey, M.~Lin, {Large-scale
  fluid simulation using velocity-vorticity domain decomposition}, ACM
  Transactions on Graphics 31~(6) (2012) 1.

\bibitem{Goude2012}
A.~Goude, S.~Engblom, {Adaptive fast multipole methods on the GPU}, The Journal
  of Supercomputing 63~(3) (2012) 897--918.

\bibitem{Greengard1987a}
L.~Greengard, V.~Rokhlin, {A fast algorithm for particle simulations}, Journal
  of Computational Physics 73~(2) (1987) 325--348.

\bibitem{Gresho1991}
P.~M. Gresho, {Incompressible fluid dynamics: Some fundamental formulation
  issues}, Annual Review of Fluid Mechanics 23~(1) (1991) 413--453.

\bibitem{Guermond2006}
J.~Guermond, P.~Minev, J.~Shen, {An overview of projection methods for
  incompressible flows}, Computer Methods in Applied Mechanics and Engineering
  195~(44-47) (2006) 6011--6045.

\bibitem{Guermond1993}
J.-L. Guermond, S.~Huberson, W.-Z. Shen, {Simulation of 2D External Viscous
  Flows by Means of a Domain Decomposition Method}, Journal of Computational
  Physics 108~(2) (1993) 343--352.

\bibitem{guermond1994}
J.-L. Guermond, W.-Z. Shen, {A domain decomposition method for simulating 2D
  external viscous flows}, in: Domain Decomposition Methods in Science and
  Engineering: The Sixth International Conference on Domain Decomposition,
  American Mathematical Soc., 1994, pp. 463--467.

\bibitem{Hald1979}
O.~H. Hald, {Convergence of Vortex Methods for Euler’s Equations. II}, SIAM
  Journal on Numerical Analysis 16~(5) (1979) 726--755.

\bibitem{Heywood1996}
J.~G. Heywood, R.~Rannacher, S.~Turek, {Artificial boundaries and flux and
  pressure conditions for the incompressible Navier-Stokes equations},
  International Journal for Numerical Methods in Fluids 22~(5) (1996) 325--352.

\bibitem{Howell1997}
L.~H. Howell, J.~B. Bell, {An adaptive mesh projection method for viscous
  incompressible flow}, SIAM Journal on Scientific Computing 18~(4) (1997)
  996--1013.

\bibitem{Huberson2002}
S.~G. Huberson, S.~G. Voutsinas, {Particles and grid}, Computers \& Fluids
  31~(4-7) (2002) 607--625.

\bibitem{KoumoutsakosPhD}
P.~Koumoutsakos, {Direct numerical simulations of unsteady separated flows
  using vortex methods}, Ph.D. thesis, California Institute of Technology
  (1993).

\bibitem{Koumoutsakos1997a}
P.~Koumoutsakos, {Inviscid Axisymmetrization of an Elliptical Vortex}, Journal
  of Computational Physics 138~(2) (1997) 821--857.

\bibitem{KoumoutsakosBoundary1993}
P.~Koumoutsakos, A.~Leonard, {Improved boundary integral method for inviscid
  boundary condition applications}, AIAA Journal 31~(2) (1993) 401--404.

\bibitem{Koumoutsakos1995}
P.~Koumoutsakos, A.~Leonard, {High-resolution simulations of the flow around an
  impulsively started cylinder using vortex methods}, Journal of Fluid
  Mechanics 296~(1) (1995) 1--38.

\bibitem{Koumoutsakos1994}
P.~Koumoutsakos, A.~Leonard, F.~P\'{e}pin, {Boundary Conditions for Viscous
  Vortex Methods}, Journal of Computational Physics 113~(1) (1994) 52--61.

\bibitem{Koumoutsakos1996}
P.~Koumoutsakos, D.~Shiels, {Simulations of the viscous flow normal to an
  impulsively started and uniformly accelerated flat plate}, Journal of Fluid
  Mechanics 328 (1996) 177--227.

\bibitem{Lashuk2012}
I.~Lashuk, G.~Biros, A.~Chandramowlishwaran, H.~Langston, T.-A. Nguyen,
  R.~Sampath, A.~Shringarpure, R.~Vuduc, L.~Ying, D.~Zorin, {A massively
  parallel adaptive fast multipole method on heterogeneous architectures},
  Communications of the ACM 55~(5) (2012) 101.

\bibitem{Lecointe1984}
Y.~Lecointe, J.~Piquet, {On the use of several compact methods for the study of
  unsteady incompressible viscous flow round a circular cylinder}, Computers \&
  Fluids 12~(4) (1984) 255--280.

\bibitem{Leonard1985_2003}
A.~Leonard, {Computing Three-Dimensional Incompressible Flows with Vortex
  Elements}, Annual Review of Fluid Mechanics 17 (1985) 523--559.

\bibitem{Lesoinne1996}
M.~Lesoinne, C.~Farhat, {Geometric conservation laws for flow problems with
  moving boundaries and deformable meshes, and their impact on aeroelastic
  computations}, Computer Methods in Applied Mechanics and Engineering
  134~(1-2) (1996) 71--90.

\bibitem{Lighthill1963}
M.~J. Lighthill, {Introduction: boundary layer theory}, in: Laminar boundary
  layers, Oxford University Press, 1963, pp. 46--113.

\bibitem{Logg2010}
A.~Logg, G.~N. Wells, {DOLFIN}, ACM Transactions on Mathematical Software
  37~(2) (2010) 1--28.

\bibitem{Manickathan2015}
L.~Manickathan, {Hybrid Eulerian-Lagrangian vortex particle method}, Master's
  thesis, Delft University of Technology (2015).

\bibitem{Monaghan1985}
J.~Monaghan, {Extrapolating B-splines for interpolation}, Journal of
  Computational Physics 60~(2) (1985) 253--262.

\bibitem{Morgenthal2007}
G.~Morgenthal, J.~Walther, {An immersed interface method for the Vortex-In-Cell
  algorithm}, Computers \& Structures 85~(11-14) (2007) 712--726.

\bibitem{Nordmark1991}
H.~O. Nordmark, {Rezoning for higher order vortex methods}, Journal of
  Computational Physics 97~(2) (1991) 366--397.

\bibitem{Ould-Salihi2001a}
M.~L. Ould-Salihi, G.~H. Cottet, M.~{El Hamraoui}, {Blending Finite-Difference
  and Vortex Methods for Incompressible Flow Computations}, SIAM Journal on
  Scientific Computing 22~(5) (2001) 1655--1674.

\bibitem{OxleyThesis}
G.~S. Oxley, {A 2-D Hybrid Euler-Compressible Vortex Particle Method For
  Transonic Rotorcraft Flows}, Ph.D. thesis, Carleton University (2009).

\bibitem{Papadakis2014}
G.~Papadakis, S.~G. Voutsinas, {In view of accelerating CFD simulations through
  coupling with vortex particle approximations}, Journal of Physics: Conference
  Series 524~(1) (2014) 012126.

\bibitem{Ploumhans2000}
P.~Ploumhans, G.~Winckelmans, {Vortex Methods for High-Resolution Simulations
  of Viscous Flow Past Bluff Bodies of General Geometry}, Journal of
  Computational Physics 165~(2) (2000) 354--406.

\bibitem{Ploumhans2002}
P.~Ploumhans, G.~Winckelmans, J.~Salmon, A.~Leonard, M.~Warren, {Vortex Methods
  for Direct Numerical Simulation of Three-Dimensional Bluff Body Flows:
  Application to the Sphere at Re=300, 500, and 1000}, Journal of Computational
  Physics 178~(2) (2002) 427--463.

\bibitem{Rasmussen2011}
J.~T.~J. Rasmussen, G.-H. Cottet, J.~H. Walther, {A multiresolution remeshed
  Vortex-In-Cell algorithm using patches}, Journal of Computational Physics
  230~(17) (2011) 6742--6755.

\bibitem{RaviartVortex1985}
P.~Raviart, {An analysis of particle methods}, in: Numerical Methods in Fluid
  Dynamics, vol. 1127, Springer Berlin Heidelberg, 1985, pp. 243--324.

\bibitem{Renac2013}
F.~Renac, S.~G\'{e}rald, C.~Marmignon, F.~Coquel, {Fast time implicit-explicit
  discontinuous Galerkin method for the compressible Navier-Stokes equations},
  Journal of Computational Physics 251 (2013) 272--291.

\bibitem{Rosenfeld1991}
M.~Rosenfeld, D.~Kwak, M.~Vinokur, {A fractional step solution method for the
  unsteady incompressible Navier-Stokes equations in generalized coordinate
  systems}, Journal of Computational Physics 94~(1) (1991) 102--137.

\bibitem{Rosenhead1930}
L.~Rosenhead, {The Spread of Vorticity in the Wake Behind a Cylinder},
  Proceedings of the Royal Society of London A: Mathematical, Physical and
  Engineering Sciences 127~(806) (1930) 590--612.

\bibitem{Sani1994}
R.~L. Sani, P.~M. Gresho, {Resume and remarks on the open boundary condition
  minisymposium}, International Journal for Numerical Methods in Fluids 18~(10)
  (1994) 983--1008.

\bibitem{Sbalzarini2006}
I.~F. Sbalzarini, J.~H. Walther, M.~Bergdorf, S.~E. Hieber, E.~M. Kotsalis,
  P.~Koumoutsakos, {PPM - A highly efficient parallel particle-mesh library for
  the simulation of continuum systems}, Journal of Computational Physics
  215~(2) (2006) 566--588.

\bibitem{Shankar1996}
S.~Shankar, L.~van Dommelen, {A New Diffusion Procedure for Vortex Methods},
  Journal of Computational Physics 127~(1) (1996) 88--109.

\bibitem{ShielsPhD}
D.~Shiels, {Simulation of controlled bluff body flow with a viscous vortex
  method}, Ph.D. thesis, California Institute of Technology (1998).

\bibitem{speckThesis}
D.~R. Speck, {Generalized Algebraic Kernels and Multipole Expansions for
  massively parallel Vortex Particle Methods}, Ph.D. thesis, Bergische
  Universitat Wuppertal (2011).

\bibitem{Stock2010}
M.~Stock, A.~Gharakhani, C.~Stone, {Modeling Rotor Wakes with a Hybrid
  OVERFLOW-Vortex Method on a GPU Cluster}, in: 28th AIAA Applied Aerodynamics
  Conference, 2010.

\bibitem{TaylorHood1973}
C.~Taylor, P.~Hood, {A numerical solution of the Navier-Stokes equations using
  the finite element technique}, Computers \& Fluids 1~(1) (1973) 73--100.

\bibitem{Tezduyar2001}
T.~E. Tezduyar, {Finite Element Methods for Flow Problems with Moving
  Boundaries and Interfaces}, Archives of Computational Methods in Engineering
  8~(2) (2001) 83--130.

\bibitem{Wee2006}
D.~Wee, A.~F. Ghoniem, {Modified interpolation kernels for treating diffusion
  and remeshing in vortex methods}, Journal of Computational Physics 213~(1)
  (2006) 239--263.

\bibitem{winckelmans1993}
G.~Winckelmans, A.~Leonard, {Contributions to Vortex Particle Methods for the
  Computation of Three-Dimensional Incompressible Unsteady Flows}, Journal of
  Computational Physics 109~(2) (1993) 247--273.

\bibitem{winckelmansThesis}
G.~S. Winckelmans, {Topics in vortex methods for the computation of three- and
  two-dimensional incompressible unsteady flows}, Ph.D. thesis, California
  Institute of Technology, USA (1989).

\bibitem{WinckelmansGeneralVortex2004}
G.~S. Winckelmans, {Vortex methods}, in: Encyclopedia of Computational
  Mechanics, vol.~3, John Wiley \& Sons, 2004, pp. 129--153.

\bibitem{Yokota2012}
R.~Yokota, L.~A. Barba, {A tuned and scalable fast multipole method as a
  preeminent algorithm for exascale systems}, International Journal of High
  Performance Computing Applications 26~(4) (2012) 337--346.

\end{thebibliography}





\end{document}